\documentclass{article}

\usepackage{amsmath,amssymb,amsthm,amscd,mathrsfs,amsxtra,mathtools}

\usepackage{tikz}
\usetikzlibrary{calc,intersections}
\usepackage{float}

\usepackage[
  textwidth=125mm,
  textheight=185mm
]{geometry}

\usepackage{tikz-cd}
\usepackage{enumitem}
\setlist[enumerate,1]{font=\normalfont, label=(\roman*)}

\usepackage[american]{babel}
\usepackage{hyperref}
\usepackage[sort=use]{glossaries}

\renewcommand{\glsglossarymark}[1]{}

\usepackage[postpunc=;]{glossaries-extra}
\GlsXtrEnablePreLocationTag{page }{pages }

\newtheorem{theorem}{Theorem}[section]

\newtheorem{lemma}[theorem]{Lemma}

\newtheorem{proposition}[theorem]{Proposition}

\newtheorem{corollary}[theorem]{Corollary}

\theoremstyle{definition}
\newtheorem{definition}[theorem]{Definition}

\newtheorem{example}[theorem]{Example}

\newtheorem{remark}[theorem]{Remark}

\DeclareMathOperator{\cone}{cone}

\DeclareMathOperator{\Img}{Im}
\DeclareMathOperator{\Hom}{Hom}
\DeclareMathOperator{\End}{End}
\DeclareMathOperator{\Ext}{Ext}

\DeclareMathOperator{\Coh}{Coh}

\DeclareMathOperator{\cok}{cok}
\DeclareMathOperator{\supp}{supp}
\DeclareMathOperator{\Slice}{Slice}
\DeclareMathOperator{\ch}{ch}

\DeclareMathOperator{\Per}{Per}

\DeclareMathOperator{\Stab}{Stab}

\DeclareMathOperator{\pr}{pr}
\DeclareMathOperator{\Br}{Br}

\makeglossaries

\newglossaryentry{hearts}{
    name={\ensuremath{\heartsuit^b(\mathcal{D})}},
    description={the set of hearts of bounded $t$-structures on $\mathcal{D}$}
}

\newglossaryentry{verdier-quotient}{
    name={\ensuremath{\pi_{\mathcal{N}}\colon \mathcal{D}\longrightarrow \mathcal{D/N}}},
    description={the Verdier quotient functor}
}

\newglossaryentry{per-heart}{
    name={\ensuremath{\Per(\mathcal{D},\mathcal{N};\mathcal{A},\mathcal{H_N})}},
    description={the subcategory $\pi_{\mathcal{N}}^{-1}(\mathcal{A})\cap{}^\perp(\mathcal{N}^{\geq1}_{\mathcal{H_N}})\cap(\mathcal{N}^{\leq-1}_{\mathcal{H_N}})^\perp$ of $\mathcal{D}$}
}

\newglossaryentry{flat}{
    name={\ensuremath{\flat \colon \Stab^{\mathrm{w}}_\Lambda(\mathcal{D},\mathcal{N}) \longrightarrow \Stab_\Lambda(\mathcal{D}/\mathcal{N})}},
    description={the descent map}
}

\newglossaryentry{sharp}{
    name={\ensuremath{\sharp_{\mathcal{H_N}} \colon \Stab_\Lambda(\mathcal{D}/\mathcal{N}) \longrightarrow \Stab^{\mathrm{w}}_\Lambda(\mathcal{D},\mathcal{N})}},
    description={the lifting map with prescribed heart $\mathcal{H_N}$}
}

\newglossaryentry{stability-space}{
    name={\ensuremath{\Stab_\Lambda(\mathcal{D})}},
    description={the space of Bridgeland stability conditions on $\mathcal{D}$ with respect to the fixed class map to $\Lambda$}
}

\newglossaryentry{weak-stability-space}{
    name={\ensuremath{\Stab^{\mathrm{w}}_\Lambda(\mathcal{D})}},
    description={the space of weak stability conditions on $\mathcal{D}$ with respect to the fixed class map to $\Lambda$}
}

\newglossaryentry{weak-stability-stratum}{
    name={\ensuremath{\Stab^{\mathrm{w}}_\Lambda(\mathcal{D},\mathcal{N})}},
    description={the subspace of $\Stab^{\mathrm{w}}_\Lambda(\mathcal{D})$ consisting of conditions whose massless subcategory is exactly $\mathcal{N}$}
}

\newglossaryentry{massless}{
    name={\ensuremath{\mathcal{N}_\sigma}},
    description={the full subcategory of objects $E$ with $m_\sigma(E)=0$}
}

\newglossaryentry{null-heart-stratum}{
    name={\ensuremath{\mathcal{S}_{\mathcal{H_N}}}},
    description={the locus of $(Z,\mathcal{B})\in\Stab^{\mathrm{w}}_\Lambda(\mathcal{D},\mathcal{N})$ satisfying $\mathcal{B}\cap\mathcal{N}=\mathcal{H_N}$}
}

\newglossaryentry{lifting-image}{
    name={\ensuremath{\mathcal{I}_{\mathcal{H_N}}}},
    description={the image $\flat(\mathcal{S}_{\mathcal{H_N}})\subseteq\Stab_\Lambda(\mathcal{D}/\mathcal{N})$}
}

\newglossaryentry{Vwindow}{
    name={\ensuremath{\mathcal{V}_q}},
    description={the subcategory ${}^\perp(\mathcal{N}^{\geq q})\cap(\mathcal{N}^{\leq q})^\perp$}
}

\newglossaryentry{Wwindow}{
    name={\ensuremath{\mathcal{W}_q}},
    description={the subcategory ${}^\perp(\mathcal{N}^{\geq {q+1}})\cap(\mathcal{N}^{\leq {q-1}})^\perp$}
}

\newglossaryentry{dsl}{
    name={\ensuremath{d^{\mathcal{D}}_{\mathrm{sl}}(\mathcal{P},\mathcal{Q})}},
    description={the generalized metric on slicings given by the uniform difference of their largest and smallest Harder--Narasimhan phases}
}

\newglossaryentry{product-metric}{
    name={\ensuremath{d^{\mathcal{D}}((Z,\mathcal{P}),(W,\mathcal{Q}))}},
    description={the generalized metric $\max\{\|Z-W\|,d_{\mathrm{sl}}^{\mathcal{D}}(\mathcal{P},\mathcal{Q})\}$ between $(Z,\mathcal{P})$ and $(W,\mathcal{Q})$}
}

\newglossaryentry{stab-boundary}{
    name={\ensuremath{\partial\Stab_\Lambda(\mathcal{D})}},
    description={the set $\overline{\Stab_\Lambda(\mathcal{D})}\setminus\Stab_\Lambda(\mathcal{D})$, with closure in $\Hom_{\mathbb{Z}}(\Lambda,\mathbb{C})\times\Slice^{\mathrm{all}}(\mathcal{D})$}
}

\newglossaryentry{Af}{
    name={\ensuremath{\mathcal{A}_f}},
    description={the category $\Coh(Y)\cap\ker\mathbf{R}f_\ast$ for $f\colon Y\longrightarrow X$}
}

\newglossaryentry{Api}{
    name={\ensuremath{\mathcal{A}_\pi}},
    description={the category $\Coh(\widetilde{X})\cap\ker\mathbf{R}\pi_\ast$ for the minimal resolution $\pi\colon\widetilde{X}\longrightarrow X$}
}

\newglossaryentry{BrX}{
    name={\ensuremath{\Br_X}},
    description={the product $\prod_{i=1}^r\Br(\Delta_i)$ acting by spherical twists along the exceptional curves}
}

\title{Stability manifolds and resolutions of singularities}  
\author{Tzu-Yang Chou}
\date{}

\usepackage{microtype}
\begin{document}
\maketitle

\begin{abstract}
In this article, we study descent and lifting of Bridgeland stability conditions along Verdier quotients.
We give a torsion-pair criterion for lifting every stability condition on the quotient to a weak stability condition with a prescribed heart on the null category.
When this criterion holds for every heart on the null category, the space of weak stability conditions with this massless subcategory is homeomorphic to the product of the stability space of the quotient and the discrete set of hearts of bounded $t$-structures on the null category.

We apply these results to contractions of relative dimension one, resolutions of nodal varieties, and blowups of smooth projective varieties along smooth centers.
For projective surfaces with ADE singularities, we identify the discrete factor through a simply transitive braid group action.
When the minimal resolution is a K3 surface, every weak stability condition obtained by lifting a numerical stability condition from the singular surface is the limit of a continuous, quasi-convergent family of numerical Bridgeland stability conditions on the resolution with the same heart.
\end{abstract}



\section{Introduction}

In \cite{Bri07}, Bridgeland introduced stability conditions on triangulated categories and proved that their spaces admit complex manifold structures.
Recently, Li \cite{Li26} proved that the bounded derived category of every smooth projective variety over $\mathbb{C}$ admits a Bridgeland stability condition,
and this existence result has been extended to projective schemes over arbitrary fields in \cite{LLLMPSZ26}.

With this existence result, it is natural to study the relations between stability conditions on different varieties and,
more generally, on different triangulated categories.
In this article, we study the relations between stability conditions on a triangulated category and those on its Verdier localizations.

In \cite{Cho26}, we constructed a stability condition on a projective surface with an ADE singularity by descending a weak stability condition from its resolution, 
with massless subcategory given by the kernel of pushforward.

This construction motivates the question of when a stability condition on a Verdier localization arises from descent of weak stability conditions.
We establish a torsion-pair criterion to lift every stability condition on the Verdier localization along a prescribed heart on the kernel.

Our main results are the following:
\begin{theorem}\label{thm:main1}
Let $\mathcal{D}$ be a triangulated category, let $\mathcal{N}\subset \mathcal{D}$ be a thick subcategory, and let $\pi_\mathcal{N}\colon \mathcal{D}\to \mathcal{D/N}$ be the Verdier quotient.
Then there is a natural continuous descent map
\[\flat\colon \Stab_{\Lambda}^{\mathrm{w}}(\mathcal{D},\mathcal{N}) \to \Stab_{\Lambda}(\mathcal{D/N}).\]

Moreover, assume that $\Stab_\Lambda(\mathcal{D/N}) \neq \varnothing$ and let $\mathcal{H_N} \subseteq \mathcal{N}$ be a heart of a bounded $t$-structure.
Then the pairs $(\mathcal{N}_{\mathcal{H_N}}^{\leq 0},(\mathcal{N}_{\mathcal{H_N}}^{\leq 0})^\perp)$ and $({}^\perp(\mathcal{N}_{\mathcal{H_N}}^{\geq 0}),\mathcal{N}_{\mathcal{H_N}}^{\geq 0})$ are torsion pairs in $\mathcal{D}$ if and only if 
$(Z,\mathcal{A}) \longmapsto (Z \circ \pi_\mathcal{N}, \Per(\mathcal{D},\mathcal{N};\mathcal{A},\mathcal{H_N}))$,
where $\Per(\mathcal{D},\mathcal{N};\mathcal{A},\mathcal{H_N}):=\pi_{\mathcal{N}}^{-1}(\mathcal{A})
\cap{}^\perp(\mathcal{N}_{\mathcal{H_N}}^{\geq 1})
\cap(\mathcal{N}_{\mathcal{H_N}}^{\leq -1})^\perp $, defines a continuous lifting map
\[\sharp_{\mathcal{H_N}}\colon \Stab_{\Lambda}(\mathcal{D/N}) \to \Stab_{\Lambda}^{\mathrm{w}}(\mathcal{D},\mathcal{N})\]
such that $\flat \circ \sharp_{\mathcal{H_N}} = 1_{\Stab_\Lambda(\mathcal{D/N})}$.
\end{theorem}

The two torsion pairs allow us to lift truncation triangles, hearts of bounded $t$-structures, and Harder--Narasimhan filtrations from $\mathcal{D}/\mathcal{N}$.

In particular, Theorem \ref{thm:main1} implies that every stability condition on, for example, a projective variety $X$ with rational singularities admitting a relative-dimension-one resolution,
arises by descent from a weak stability condition on the resolution.

Moreover, a weak stability condition with massless subcategory $\mathcal{N}$ is uniquely determined by its descent and its heart on $\mathcal{N}$.
If the torsion-pair criterion above holds for every heart of a bounded $t$-structure on $\mathcal{N}$,
we obtain the following result,
which decomposes the space of weak stability conditions with a prescribed massless subcategory into a disjoint union of clopen copies of the stability space of the quotient.

\begin{theorem}\label{thm:main2} 
Assume that $\Stab_\Lambda(\mathcal{D/N})\neq\varnothing$.
Then the pairs $(\mathcal{N}_{\mathcal{H_N}}^{\leq 0},(\mathcal{N}_{\mathcal{H_N}}^{\leq 0})^\perp)$ and $({}^\perp(\mathcal{N}_{\mathcal{H_N}}^{\geq 0}),\mathcal{N}_{\mathcal{H_N}}^{\geq 0})$ are torsion pairs in $\mathcal{D}$  for every heart $\mathcal{H_N}$ of a bounded $t$-structure on $\mathcal{N}$
if and only if the map $\tau=(Z,\mathcal{B})\longmapsto(\mathcal{B}\cap\mathcal{N},\tau^\flat)$ gives a natural homeomorphism
\begin{align*}
\Stab^{\mathrm{w}}_\Lambda(\mathcal{D},\mathcal{N})
\simeq
\heartsuit^{b}(\mathcal{N})\times\Stab_\Lambda(\mathcal{D/N}),
\end{align*}
where $\heartsuit^{b}(\mathcal{N})$ is the set of hearts of bounded $t$-structures on $\mathcal{N}$, equipped with the discrete topology.
\end{theorem}

For each prescribed heart on $\mathcal{N}$, the lifting map is an isometry,
and the images corresponding to different hearts are separated by a slicing distance of at least $1$.

We now apply this result to morphisms $f\colon Y\longrightarrow X$ for which derived pushforward induces a Verdier localization.
In the following theorem, we take $\mathcal{D}=\mathrm{D^b}(Y)$ and $\mathcal{N}=\ker\mathbf{R}f_\ast$, and identify $\mathcal{D}/\mathcal{N}$ with $\mathrm{D^b}(X)$.

\begin{theorem}\label{thm:main3}
Let $f\colon Y\longrightarrow X$ be a projective morphism between quasi-projective schemes over an algebraically closed field $k$.
Assume that $\Stab_\Lambda(X)\neq\varnothing$.
Then the equivalent conditions of Theorem \ref{thm:main2} hold in any of the following situations:
\begin{enumerate}
\item The morphism $f$ is birational, of relative dimension one, 
$\mathbf{R}f_\ast\mathcal{O}_Y\simeq\mathcal{O}_X$, and $f(\mathrm{Exc}(f))$ is a finite set of closed points.
\item $\mathrm{char}(k) \neq 2$. $X$ is a projective variety of dimension $\geq 2$ with only ordinary double points,
and $f$ is the blowup of these singular points.
\item The variety $X$ is smooth and projective, and $f$ is the blowup along a smooth closed subvariety of codimension $\geq 2$.
\end{enumerate}
\end{theorem}

Consider the minimal resolution $\pi\colon\widetilde{X}\longrightarrow X$ of a normal projective surface $X$ over $\mathbb{C}$ with canonical singularities,
and set $\mathcal{N}:=\ker\mathbf{R}\pi_\ast$.
Theorem \ref{thm:ADE-null-stratum} shows that spherical twists induce a simply transitive action of $\prod_i\mathrm{Br}(\Delta_i)$ on $\heartsuit^b(\mathcal{N})$,
where the $\Delta_i$ are the Dynkin diagrams corresponding to the singularities of $X$.
Consequently, we obtain a homeomorphism
\begin{equation*}
\begin{aligned}
\left(\prod_i\mathrm{Br}(\Delta_i)\right)\times\Stab_\Lambda(X)
&\simeq\Stab^w_\Lambda(\widetilde{X},\mathcal{N}),\\
(\beta,\sigma)
&\longmapsto\beta\cdot\sigma^\sharp_{\Coh(\widetilde{X})\cap\mathcal{N}}.
\end{aligned}
\end{equation*}

In \cite{Cho26}, the weak stability condition on the resolution is obtained by degenerating a family of Bridgeland stability conditions with fixed heart.
We now ask whether every weak stability condition constructed above can be realized as a limit of Bridgeland stability conditions on $\widetilde{X}$, with convergence of both central charges and slicings.
We prove this when $\widetilde{X}$ is a K3 surface, while keeping the prescribed lifted heart fixed.

\begin{theorem}\label{thm:intro-BR}
If $\widetilde{X}$ is a K3 surface,
then for any $\sigma=(Z,\mathcal{A})\in\Stab_{\mathrm{num}}(X)$ and any $\mathcal{H_N}\in\heartsuit^b(\mathcal{N})$, 
there exists a continuous family, defined for all sufficiently small $\epsilon>0$,
$\sigma_\epsilon:=(Z_\epsilon,\Per(\mathrm{D^b}(\widetilde{X}),\mathcal{N};\mathcal{A},\mathcal{H_N}))\in\Stab_{\mathrm{num}}(\widetilde{X})$ converging to $\sigma^\sharp_{\mathcal{H_N}}$ as $\epsilon \longrightarrow 0$.
\end{theorem}

The support property on $X$ does not control the exceptional numerical directions on $\widetilde{X}$.
Using Serre duality on the K3 surface, we obtain a numerical estimate for suitable lifts of stable objects,
which allows us to perturb the central charge with the heart fixed and verify the Harder--Narasimhan and support properties, as well as convergence of the slicings.

Finally, we prove that the approximating families are quasi-convergent in the sense of \cite{HLJR26}; see Proposition \ref{prop:k3-quasi-convergence}.

\subsection{Related works}
In \cite{KLP26}, the authors induce $t$-structures on semiorthogonal components by constructing perverse $t$-structures on the ambient category,
while in \cite{BLMS23}, the authors give a criterion under which weak stability conditions restrict to Bridgeland stability conditions on such components.
Our constructions of hearts and stability conditions concern Verdier localization and do not require a semiorthogonal decomposition.

Jin \cite{Jin23} uses simple-minded collections to lift finite-length hearts of bounded $t$-structures from Verdier quotients.
Our lifting criterion concerns arbitrary hearts of bounded $t$-structures on the quotient, with a prescribed heart on the null category.

In \cite{BPPW22}, the authors study a partial compactification of stability manifolds by lax stability conditions.
They establish descent to the massless quotient and obtain local gluing and product descriptions of the lax closure.
For the weak stability conditions considered here, our torsion-pair criterion provides liftings of every stability condition on the quotient with a prescribed heart on the null category.
When this criterion holds for every such heart, Theorem \ref{thm:main2} gives a global product decomposition.

Bridgeland \cite{Bri09} describes a connected component of the stability space of the null category of a Kleinian resolution as a covering space,
with spherical twists acting as deck transformations.
Brav and Thomas \cite{BT11} prove the faithfulness of the associated braid group action,
while Hara and Wemyss \cite{HW26} show that every heart of a bounded $t$-structure on the null category belongs to the braid group orbit of the standard heart.

In \cite{KKO25}, the authors construct a family of stability conditions on a K3 surface with a fixed heart,
converging to a lax stability condition,
whose massless subcategory is generated by a spherical bundle.
Our K3 construction also gives convergence of central charges and slicings,
with the massless subcategory given by the kernel of pushforward.

Robotis and Vilches \cite{RV26} construct strongly quasi-convergent paths of stability conditions with a fixed heart associated to contractions of ADE configurations on K3 surfaces.
They obtain stability conditions on the singular surfaces from these paths and study the associated boundary strata in the space of augmented stability conditions.
For the same contractions, we start with any numerical stability condition on the singular surface and any prescribed heart of a bounded $t$-structure on the null category.
Theorem \ref{thm:intro-BR} realizes the weak stability condition obtained by lifting as the limit of a family of Bridgeland stability conditions on the resolution with the same heart.
This family is also quasi-convergent by Proposition \ref{prop:k3-quasi-convergence}.

\subsection{Plan of the paper}
This paper is organized as follows.
Section 2 contains preliminaries on bounded $t$-structures, torsion pairs, and stability conditions,
together with the topology on the spaces of stability conditions.

In Section 3, we study descent and lifting of (weak) stability conditions along Verdier quotients.
In \S3.1, we show that a weak stability condition descends to a Bridgeland stability condition on the quotient by its massless subcategory.
In \S3.2, we give a torsion-pair criterion for lifting hearts of bounded $t$-structures with a prescribed heart on the null category.
We then lift Harder--Narasimhan filtrations and construct weak stability conditions in \S3.3.
Finally, in \S3.4, we prove the continuity of descent and lifting and establish the global product decomposition.

Section 4 gives geometric applications of these results.
We consider contractions of relative dimension one in \S4.1, resolutions of varieties with ordinary double points in \S4.2, and blowups of smooth varieties along smooth centers in \S4.3.

In Section 5, we specialize to projective surfaces with canonical singularities.
We describe the hearts of bounded $t$-structures on the null category by braid group actions in \S5.1,
and give a corresponding description for partial resolutions.
In \S5.2, we assume that the minimal resolution is a K3 surface and realize every weak stability condition obtained by lifting a numerical stability condition as the endpoint of a quasi-convergent path of Bridgeland stability conditions with a fixed heart.

\subsection*{Acknowledgements}
I would like to express my gratitude to Hsueh-Yung Lin and Yukinobu Toda for valuable discussions and suggestions  throughout the development of this work.
Thanks to Arend Bayer and Tianle Mao for their helpful comments.

This work was supported by the National Science and Technology Council of Taiwan under grant No. 115-2811-M-002-146 and by the Yushan Fellow Program of the Ministry of Education under grant No. NTU-114V1006-5.

\paragraph{Conventions.}
For a triangulated category $\mathcal{D}$ and $S\subseteq\mathcal{D}$, we write $\langle S\rangle$ for the smallest thick subcategory of $\mathcal{D}$ containing $S$, and $\langle S\rangle_{ex}$ for its extension closure.
Let $S^\perp:=\{E\in\mathcal{D}\mid\Hom_{\mathcal{D}}(S,E)=0\}$ and ${}^\perp S:=\{E\in\mathcal{D}\mid\Hom_{\mathcal{D}}(E,S)=0\}$.
When working in a heart, we take orthogonals inside that heart.

\printglossary[title= List of Notations]

\section{Preliminaries}
In this section, we review the notions of stability conditions and bounded $t$-structures.

Let $\mathcal{D}$ be a triangulated category.

\begin{definition}
We say that a $k$-linear triangulated category $\mathcal{D}$ is proper if $\bigoplus_{i \in \mathbb{Z}} \Hom_{\mathcal{D}}(X,Y[i])$ is finite-dimensional over $k$ for all $X,Y \in \mathcal{D}$.

Equivalently, $\mathcal{D}$ is Hom-finite and $\Hom_{\mathcal{D}}(X,Y[i])=0$ for all but finitely many $i$.
\end{definition}

The most crucial examples are the derived category $\mathrm{D^b}(V)$ of a smooth, proper scheme $V$ and its full triangulated subcategories.

For a heart $\mathcal{A}$ of a bounded $t$-structure on a triangulated category $\mathcal{D}$,
we denote by $(\mathcal{D}_{\mathcal{A}}^{\leq 0},\mathcal{D}_{\mathcal{A}}^{\geq 0})$ the corresponding bounded $t$-structure.
When the heart is clear from the context, we omit it from the subscripts and write $\mathcal{D}^{\leq r}$ and $\mathcal{D}^{\geq r}$.
We also write \gls{hearts} for the set of hearts of bounded $t$-structures on $\mathcal{D}$.

\begin{definition}
A torsion pair in $\mathcal{D}$ is a pair $(\mathcal{T},\mathcal{F})$ of full additive subcategories closed under direct summands such that $\Hom_{\mathcal{D}}(\mathcal{T},\mathcal{F})=0$ and every $E\in\mathcal{D}$ admits a distinguished triangle $T\longrightarrow E\longrightarrow F\longrightarrow T[1]$ with $T\in\mathcal{T}$ and $F\in\mathcal{F}$.
\end{definition}

For a $t$-structure $(\mathcal{D}^{\leq0},\mathcal{D}^{\geq0})$, the pair $(\mathcal{D}^{\leq0},\mathcal{D}^{\geq1})$ is a torsion pair in $\mathcal{D}$.

\begin{definition}
Let $\mathcal{A}$ be an abelian category.
A torsion pair in $\mathcal{A}$ is a pair $(\mathcal{T},\mathcal{F})$ of full additive subcategories such that $\Hom(\mathcal{T},\mathcal{F})=0$ and every $A\in\mathcal{A}$ admits a short exact sequence $0\longrightarrow T\longrightarrow A\longrightarrow F\longrightarrow0$ with $T\in\mathcal{T}$ and $F\in\mathcal{F}$.
\end{definition}

\begin{theorem}[{\normalfont\cite{HRS96}}]\label{HRStilt}
Let $\mathcal{A}$ be the heart of a bounded $t$-structure in $\mathcal{D}$ with a torsion pair $(\mathcal{T}, \mathcal{F})$.
The category
\begin{equation*}
\mathcal{A}^\sharp := \{ E \in \mathcal{D} \mid H^i(E) = 0 \ {\rm if} \ i \neq 0, -1, H^0(E) \in \mathcal{T} , H^{-1}(E) \in  \mathcal{F} \}
\end{equation*}
is the heart of a bounded $t$-structure on $\mathcal{D}$.
\end{theorem}

Conversely, if $\mathcal{A},\mathcal{B}\in\heartsuit^b(\mathcal{D})$ satisfy $\mathcal{A}\subseteq\mathcal{D}_{\mathcal{B}}^{[0,1]}$, then $\mathcal{B}$ is the tilt of $\mathcal{A}$ with respect to the torsion pair $(\mathcal{A}\cap\mathcal{B},\mathcal{A}\cap\mathcal{B}[-1])$ by \cite[Lemma 1.1.2]{Pol07}.

\begin{definition}
Let $\mathcal{A}$ be an abelian category.

A group homomorphism $Z\colon \mathrm{K}(\mathcal{A}) \longrightarrow \mathbb{C}$ is said to be a weak stability function (resp. stability function) if for any nonzero object $E\in \mathcal{A}$, we have $\Im Z(E) \geq 0$ with $\Im Z(E) = 0 \Rightarrow \Re Z(E) \leq 0$ (resp. $ \Re Z(E) < 0$).

For a weak stability function $Z$ on $\mathcal{A}$ and any $E\in\mathcal{A}$, we set
\begin{equation}
\mu_Z(E):=
\begin{cases}
\frac{-\Re Z(E)}{\Im Z(E)}, & \text{if } \Im Z(E)>0, \\
+\infty, & \text{if } \Im Z(E)=0.
\end{cases}
\end{equation}

We say that an object $0\neq E \in \mathcal{A}$ is $Z$-semistable (resp. $Z$-stable) if
for every subobject $0\neq F \subsetneq E$, we have $\mu_Z (F) \leq \mu_Z (E)$ (resp. $\mu_Z (F) < \mu_Z (E)$).
The phase of a $Z$-semistable object $0\neq E\in\mathcal{A}$ is $\pi^{-1}\arg Z(E)\in(0,1]$ if $Z(E)\neq0$, and $1$ otherwise.
\end{definition}

We fix a finite rank lattice $\Lambda$ and a  group homomorphism $v\colon \mathrm{K}(\mathcal{D}) \twoheadrightarrow \Lambda$.

\begin{definition}\label{stability}
A weak stability condition (resp. Bridgeland stability condition) on $\mathcal{D}$ is a pair $\sigma = (Z,\mathcal{A})$ consisting of a group homomorphism (called the central charge of $\sigma$) $Z \colon \Lambda \longrightarrow \mathbb{C}$ and a heart $\mathcal{A}$ of a bounded $t$-structure on $\mathcal{D}$.
We write $\mu_\sigma:=\mu_Z$ and use $\sigma$-(semi)stability for $Z$-(semi)stability.
The following conditions are required:
\begin{enumerate}
    \item[(a)] The composition $Z \circ v \colon \mathrm{K}(\mathcal{A})=\mathrm{K}(\mathcal{D}) \longrightarrow \Lambda \longrightarrow \mathbb{C}$ is a weak stability function (resp. stability function) on $\mathcal{A}$.

    \item[(b)] For any $0\neq E \in \mathcal{A}$, there is a filtration
    \[0=E_0 \subset E_1\subset \dots \subset E_{n-1} \subset E_n=E,\]
    such that the quotients $A_i:=E_i/E_{i-1}$ are $\sigma$-semistable (called the Harder--Narasimhan factors) with $\mu_\sigma(A_1)> \dots >\mu_\sigma(A_n)$.
    \item[(c)] (Support property) There exists a quadratic form $Q$ on $\Lambda \otimes \mathbb{R}$ such that $Q|_{\ker Z}$ is negative definite, and $Q(v(E))\geq 0$, for all $\sigma$-semistable $E \in \mathcal{A}$.
\end{enumerate}
Any pair $\sigma=(Z,\mathcal{A})$ satisfying the conditions (a) and (b) is said to be a weak pre-stability condition (resp. pre-stability condition).
\end{definition}

For a fixed norm on $\Lambda_{\mathbb{R}}$, the support property is equivalent to the existence of $K>0$ such that $\|v(E)\|\leq K|Z(E)|$ for every semistable $E$.

We let \gls{stability-space} and \gls{weak-stability-space} be the sets of Bridgeland stability conditions and weak stability conditions on $\mathcal{D}$ with respect to the fixed class map $v\colon \mathrm{K}(\mathcal{D})\longrightarrow\Lambda$, respectively.

We write $\Slice^{\mathrm{all}}(\mathcal{D})$ for the set of all (possibly non-locally finite) slicings on $\mathcal{D}$ in the sense of \cite[Definition 3.3]{Bri07},
and we will freely pass between $(Z,\mathcal{A})$ and $(Z,\mathcal{P})$, with $\mathcal{A}=\mathcal{P}(0,1]$.

For a slicing $\mathcal{P}$, we denote the largest and smallest Harder--Narasimhan phases of a nonzero object $E$ by $\phi^\pm_{\mathcal{P}}(E)$ respectively, 
and set $\mathcal{P}(I):=\langle \mathcal{P}(\phi)\ |\ \phi\in I\rangle_{ex}$, where $I\subseteq\mathbb{R}$ is an interval.

Let $\sigma=(Z_\mathcal{B},\mathcal{B})$ be a weak pre-stability condition on
$\mathcal{D}$.

\begin{definition}
For a nonzero object $E\in\mathcal{D}$,
we define the mass of $E$ by $m_\sigma(E):=\sum_j\sum_{i=1}^{m_j}|Z_{\mathcal{B}}(E_{j,i})|$,
where $E_{j,1},\dots,E_{j,m_j}$ are the Harder--Narasimhan factors of $H^j_{\mathcal{B}}(E)$.
We also define $m_\sigma(0):=0$.

We say that an object $E\in \mathcal{D}$ is massless if $m_{\sigma}(E)=0$,
and we define the massless subcategory \gls{massless}$\subseteq \mathcal{D}$ of $\sigma$ to be the full subcategory of massless objects.
\end{definition}

\begin{lemma}\label{lem:massless}
$\mathcal{N}_\sigma$ is thick and $\mathcal{N}_\sigma=\langle \mathcal{B} \cap \ker Z_\mathcal{B}\rangle$.
\end{lemma}
\begin{proof}
The first assertion is clear.
Now let $E\in \mathcal{B}\cap \ker Z_\mathcal{B}$.
As $Z_\mathcal{B}(E)=0$, $E$ itself is $\sigma$-semistable,
and hence $m_\sigma(E)=0$.

Conversely, let $E \in \mathcal{N}_\sigma$ be a massless object.
Then each $H^j_\mathcal{B}(E)$ belongs to $\mathcal{B} \cap \ker Z_\mathcal{B}$ and hence $E \in \langle \mathcal{B} \cap \ker Z_\mathcal{B}\rangle$.
\end{proof}

For a thick subcategory $\mathcal{N}\subseteq\mathcal{D}$, we write \gls{weak-stability-stratum} for the subset of $\Stab^{\mathrm{w}}_\Lambda(\mathcal{D})$ consisting of weak stability conditions whose massless subcategory is exactly $\mathcal{N}$.

For any slicings $\mathcal{P},\mathcal{Q}$ on $\mathcal{D}$,
we recall from \cite[Section 6]{Bri07} the following generalized metric and its equivalent characterization, both of which remain valid without assuming local finiteness:
\begin{equation}\label{eq:slicing-metric}
\begin{aligned}
\gls{dsl}
&:=
\sup_{0\neq E\in\mathcal{D}}
\max\{
|\phi^+_{\mathcal{P}}(E)-\phi^+_{\mathcal{Q}}(E)|,
|\phi^-_{\mathcal{P}}(E)-\phi^-_{\mathcal{Q}}(E)|
\}\\
&=
\inf\{
\epsilon>0:
\mathcal{Q}(\phi)\subseteq
\mathcal{P}([\phi-\epsilon,\phi+\epsilon])
\text{ for every }\phi\in\mathbb{R}
\}.
\end{aligned}
\end{equation}
Fix a norm $\|\cdot\|$ on $\Hom_{\mathbb{Z}}(\Lambda,\mathbb{C})$. We equip
$\Hom_{\mathbb{Z}}(\Lambda,\mathbb{C})\times\Slice^{\mathrm{all}}(\mathcal{D})$
with the topology induced by the generalized metric
\begin{equation}\label{eq:stability-metric}
\gls{product-metric}
:=
\max\{
\|Z-W\|,
d^{\mathcal{D}}_{\mathrm{sl}}(\mathcal{P},\mathcal{Q})
\},
\end{equation}
that is, the product topology.
On $\Stab_\Lambda(\mathcal{D})$ the induced topology is the usual Bridgeland topology by \cite[Section 6]{Bri07}.

Throughout this paper, spaces of (weak) stability conditions carry the subspace topology induced by this product topology for the specified category, lattice and class map.
We denote by $\gls{stab-boundary}:=\overline{\Stab_{\Lambda}(\mathcal{D})}\setminus\Stab_{\Lambda}(\mathcal{D})$ the set of points added by taking this closure.

\section{Descent and lifting of weak stability conditions}

In this section, we study the descent and lifting of weak stability conditions.

Throughout this section, we let $\mathcal{D}$ be a triangulated category, $\mathcal{N}\subseteq\mathcal{D}$ be a thick subcategory,
and \gls{verdier-quotient} be the Verdier quotient functor.

\subsection{Descent of weak stability conditions}
We fix a finite rank lattice $\Lambda$ and a surjection $v \colon \mathrm{K}(\mathcal{D})\longrightarrow \Lambda$. 

We will use the following standard criterion for descending a bounded $t$-structure.

\begin{lemma}\label{lem:quotient-heart-criterion}
Let $\mathcal{A}\in\heartsuit^b(\mathcal{D})$.
Then $\pi_{\mathcal{N}}(\mathcal{A})\in\heartsuit^b(\mathcal{D}/\mathcal{N})$ if and only if 
$\mathcal{A}\cap\mathcal{N}$ is a Serre subcategory of $\mathcal{A}$ and $H^i_{\mathcal{A}}(\mathcal{N})\subseteq\mathcal{N}$ for every $i$.

In this case, $\mathcal{A}\cap\mathcal{N}\in\heartsuit^b(\mathcal{N})$, $\pi_{\mathcal{N}}$ is $t$-exact, and induces an equivalence $\mathcal{A}/(\mathcal{A}\cap\mathcal{N})=\pi_{\mathcal{N}}(\mathcal{A})$.
\end{lemma}

We also need the following useful lemma.
\begin{lemma}[{\normalfont\cite[Lemma 14.1.9]{Kra22}}]\label{lem:Kra22}
Let $p \colon \mathcal{B}\longrightarrow \mathcal{B}/ \ker p$ be a Serre quotient.
Then for every short exact sequence \[0\longrightarrow A'\longrightarrow p(E)\longrightarrow A''\longrightarrow 0\]
in $\mathcal{B}/\ker p$,
there exists a short exact sequence \[0\longrightarrow E'\longrightarrow E\longrightarrow E''\longrightarrow 0\] in $\mathcal{B}$ such that $p(E')= A'$ and $p(E'')= A''$.
\end{lemma}

\begin{proposition}\label{prop:weak-BS-descend}
Let $\sigma=(Z_\mathcal{B},\mathcal{B})$ be a weak pre-stability condition on
$\mathcal{D}$ with respect to $(\Lambda,v)$ and let $\mathcal{N}_\sigma\subset \mathcal{D}$ be its massless subcategory.
Then:
\begin{enumerate}
\item There is a unique pre-stability condition $\sigma^\flat=(Z,\pi_{\mathcal{N}_\sigma}(\mathcal{B}))$ on $\mathcal{D/N}_\sigma$, with respect to the quotient lattice $\Lambda/\bigl(v(\mathrm{K}(\mathcal{N}_\sigma))\bigr)^{\mathrm{sat}}$, such that $Z_{\mathcal{B}}=Z\circ\pi_{\mathcal{N}_\sigma}$;

\item if $E\in \mathcal{B}$ is $\sigma$-semistable, then $\pi_{\mathcal{N}_\sigma}(E)$ is either zero or $\sigma^\flat$-semistable;

\item for every $\sigma^\flat$-semistable object $A\in \pi_{\mathcal{N}_\sigma}(\mathcal{B})$, there exists a $\sigma$-semistable object $E\in \mathcal{B}$ such that $\pi_{\mathcal{N}_\sigma}(E)=A$;

\item Assume moreover that $\sigma$ is a weak stability condition. Then $v(K(\mathcal{N}_\sigma))=0$, and $\sigma^\flat$ is a stability condition with respect to $\Lambda$.
\end{enumerate}
\end{proposition}

\begin{proof}
We first prove that $\pi_{\mathcal{N}_\sigma}(\mathcal{B}) \in \heartsuit^b(\mathcal{D}/\mathcal{N}_\sigma)$.
By the definition of mass, we have $H^i_{\mathcal{B}}(\mathcal{N}_\sigma) \subseteq \mathcal{N}_\sigma$ for every $i$, 
and $\mathcal{B} \cap \mathcal{N}_\sigma= \mathcal{B} \cap \ker Z_{\mathcal{B}}$,
which is a Serre subcategory of $\mathcal{B}$ since $Z_{\mathcal{B}}$ is a weak stability function.
Therefore, Lemma \ref{lem:quotient-heart-criterion} applies.

For every $E\in\mathcal{N}_\sigma$, we have $Z_{\mathcal{B}}(E)=\sum_i(-1)^i Z_{\mathcal{B}}(H^i_{\mathcal{B}}(E))=0$,
so $Z_{\mathcal{B}}$ vanishes on $v(\mathrm{K}(\mathcal{N}_\sigma))$, and hence on its saturation.
Then, there is an induced map $Z\colon\Lambda/\bigl(v(\mathrm{K}(\mathcal{N}_\sigma))\bigr)^{\mathrm{sat}}\longrightarrow\mathbb{C}$
such that $Z_{\mathcal{B}}=Z\circ\pi_{\mathcal{N}_\sigma}$.

Next, we show that $Z$ is a stability function on $\pi_{\mathcal{N}_\sigma}(\mathcal{B})$. 
Given a nonzero object $A\in \pi_{\mathcal{N}_\sigma}(\mathcal{B})$,
we choose $E\in \mathcal{B}$ with $\pi_{\mathcal{N}_\sigma}(E)= A$.
Then $Z(A)=Z_\mathcal{B}(E)$ and hence lies in $\mathbb{H} \cup \{0\}$.
Assume that $Z(A)=0$.
Then $Z_{\mathcal{B}}(E)=Z(\pi_{\mathcal{N}_\sigma}(E))=0$, so $E\in \mathcal{B}\cap \ker Z_{\mathcal{B}}\subseteq \mathcal{N}_\sigma$, hence $\pi_{\mathcal{N}_\sigma}(E)=A=0$, a contradiction. 

We now prove (ii). Given a $Z_\mathcal{B}$-semistable object $E\in \mathcal{B}$ with $\pi_{\mathcal{N}_\sigma}(E) \neq 0$,
we assume that $\pi_{\mathcal{N}_\sigma}(E)$ is $Z$-unstable.
Then there is a nonzero proper subobject $A \subset \pi_{\mathcal{N}_\sigma}(E)$ in $\pi_{\mathcal{N}_\sigma}(\mathcal{B})$ with
$\mu_{Z}(A)>\mu_{Z}(\pi_{\mathcal{N}_\sigma}(E))$.
By Lemma \ref{lem:Kra22}, there exists a monomorphism $F\hookrightarrow E$ in $\mathcal{B}$ such that $\pi_{\mathcal{N}_\sigma}(F)=A$, and then we have
$\mu_{Z_\mathcal{B}}(F)=\mu_{Z}(A)>\mu_{Z}(\pi_{\mathcal{N}_\sigma}(E))=\mu_{Z_\mathcal{B}}(E)$,
contradicting the $Z_\mathcal{B}$-semistability of $E$. Thus $\pi_{\mathcal{N}_\sigma}(E)$ is $Z$-semistable.

We then return to (i).
It remains to prove the Harder--Narasimhan property for $(Z,\pi_{\mathcal{N}_\sigma}(\mathcal{B}))$. 
For any nonzero $A\in \pi_{\mathcal{N}_\sigma}(\mathcal{B})$, we choose $E\in \mathcal{B}$ with $\pi_{\mathcal{N}_\sigma}(E)= A$.
Then by (ii), applying the exact functor $\pi_{\mathcal{N}_\sigma}|_{\mathcal{B}}$ to the Harder--Narasimhan filtration of $E$ with respect to $Z_\mathcal{B}$,
we obtain the Harder--Narasimhan filtration of $A$ with respect to $Z$.

For (iii), given a $Z$-semistable object $A\in \pi_{\mathcal{N}_\sigma}(\mathcal{B})$, we choose $E\in \mathcal{B}$ with $\pi_{\mathcal{N}_\sigma}(E)= A$. 

We may assume that $E$ is $Z_\mathcal{B}$-unstable and take the Harder--Narasimhan filtration  $0=E_0\subset E_1\subset \cdots \subset E_m=E$ of $E$ with respect to $Z_\mathcal{B}$.
As $A$ is $Z$-semistable, $\pi_{\mathcal{N}_\sigma}(E_{i+1}/E_i)=0$ for all $i\neq m-1$.
But if $\pi_{\mathcal{N}_\sigma}(E_{i+1}/E_i)=0$,
then $\mu_{Z_\mathcal{B}}(E_{i+1}/E_i) = + \infty$.
Therefore, the filtration can only be a two-step filtration $0 \subseteq E_1 \subseteq E$,
and now $E/E_1$ would be a $Z_\mathcal{B}$-semistable lift of $A$.

For (iv), let $Q$ be a quadratic form on $\Lambda_{\mathbb{R}}$ giving the support property for $\sigma$. 
We prove that $v(K(\mathcal{N}_\sigma))=0$.
By Lemma \ref{lem:massless}, it suffices to show that $v(E)=0$ for every $E\in \mathcal{B}\cap \ker Z_{\mathcal{B}}$.

Let $E\in \mathcal{B}\cap\ker Z_{\mathcal{B}}$. 
If $v(E)\neq 0$, then $Q(v(E))<0$ by the negative definiteness of $Q|_{\ker Z_{\mathcal{B}}}$.
On the other hand, $\mu_{Z_{\mathcal{B}}}(E)=+\infty$, so $E$ is $Z_{\mathcal{B}}$-semistable and hence $Q(v(E))\geq 0$. 
This proves that $v(E)=0$.
Now, the same quadratic form $Q$ gives the support property for $\sigma^\flat$: $Q|_{\ker Z}=Q|_{\ker Z_\mathcal{B}}$ is negative definite,
and if $A\in \pi_{\mathcal{N}_\sigma}(\mathcal{B})$ is $Z$-semistable, then by (iii) there exists a $Z_{\mathcal{B}}$-semistable object $E\in \mathcal{B}$ with $\pi_{\mathcal{N}_\sigma}(E) = A$,
so $Q(v(A))=Q(v(E))\geq 0$.
\end{proof}

The proposition gives the descent map \gls{flat}.
For a related descent result for lax stability conditions, see \cite[Proposition 4.19]{BPPW22}.

\subsection{Lifting of bounded hearts}
The next step is to describe the canonical form of bounded hearts on $\mathcal{D}$ which descend to $\mathcal{D/N}$.
This leads to the perverse-style construction introduced below.

\begin{definition}
Given $\mathcal{A}\in\heartsuit^{b}(\mathcal{D}/\mathcal{N})$ and $\mathcal{H_N}\in\heartsuit^{b}(\mathcal{N})$, we define
\[
\mathcal{D}_{\mathcal{A},\mathcal{H_N}}^{\leq 0}
:=
\{E\in \mathcal{D} \mid \pi_\mathcal{N}(E)\in (\mathcal{D/N})_{\mathcal{A}}^{\leq 0},\ \Hom_{\mathcal{D}}(E,\mathcal{N}_{\mathcal{H_N}}^{\geq 1})=0\},
\]
\[
\mathcal{D}_{\mathcal{A},\mathcal{H_N}}^{\geq 0}
:=
\{E\in \mathcal{D} \mid \pi_\mathcal{N}(E)\in (\mathcal{D/N})_{\mathcal{A}}^{\geq 0},\ \Hom_{\mathcal{D}}(\mathcal{N}_{\mathcal{H_N}}^{\leq-1},E)=0\}.
\]
\[
\gls{per-heart}
:=
\mathcal{D}_{\mathcal{A},\mathcal{H_N}}^{\leq 0}\cap \mathcal{D}_{\mathcal{A},\mathcal{H_N}}^{\geq 0}.
\]
\end{definition}

\begin{example}
\begin{enumerate}
    \item In \cite{Bri02}, Bridgeland defines the perverse coherent sheaves ${^p}\Per(\widetilde{X}/X)$,
    which is exactly 
    \[\Per(\mathrm{D^b}(\widetilde{X}), \ker \mathbf{R}\pi_\ast; \Coh(X), \Coh(\widetilde{X}) \cap  \ker \mathbf{R}\pi_\ast [-p]).\]
    \item In \cite{Cho26}, we constructed, on the minimal resolution $\pi\colon \widetilde{X} \longrightarrow X$ of a projective surface  with a single ADE singularity,
    a heart $\mathcal{B}^0$ of a bounded $t$-structure.
    It can be seen that \[\mathcal{B}^0=\Per(\mathrm{D^b}(\widetilde{X}), \ker \mathbf{R}\pi_\ast; \Coh^0_H(X), (\Coh(\widetilde{X}) \cap  \ker \mathbf{R}\pi_\ast)[1]),\]
    where $H$ is an ample divisor on $X$ and $\Coh_H^0(X)$ denotes the tilt of $\Coh(X)$ with respect to the usual slope stability at slope $0$.
\end{enumerate}
\end{example}

If $\Per(\mathcal{D},\mathcal{N};\mathcal{A},\mathcal{H_N}) \in \heartsuit^b{\mathcal{D}}$, some good properties follow immediately.
More precisely, we have the following:

    \begin{proposition}\label{prop:formalprop}
Assume that $\mathcal{B}:=\Per(\mathcal{D},\mathcal{N};\mathcal{A},\mathcal{H_N})$ is a heart of a bounded $t$-structure on $\mathcal{D}$.
Then we have
$\mathcal{B} \cap \mathcal{N}=\mathcal{H_N},$
$ \pi_\mathcal{N}(\mathcal{B}) =\mathcal{A}$, and
$\mathcal{B}/\mathcal{H_N}=\mathcal{A} $.
\end{proposition}
\begin{proof}
The inclusion $\mathcal{H_N}\subset \mathcal{B}\cap \mathcal{N}$ is immediate. Conversely, let $E\in \mathcal{B}\cap \mathcal{N}$.
Then the definition of $\mathcal{B}$ gives $\Hom_{\mathcal{D}}(E,\mathcal{N}_{\mathcal{H_N}}^{\geq 1})=0$ and $\Hom_{\mathcal{D}}(\mathcal{N}_{\mathcal{H_N}}^{\leq -1},E)=0$.
This means exactly that $E\in \mathcal{N}_{\mathcal{H_N}}^{\leq 0}\cap \mathcal{N}_{\mathcal{H_N}}^{\geq 0}=\mathcal{H_N}$.

Next, again by the definition of $\mathcal{B}$ we have $\pi_{\mathcal{N}}(\mathcal{B})\subset \mathcal{A}$.
It follows that $\pi_{\mathcal{N}}$ is $t$-exact with respect to the bounded $t$-structures corresponding to $\mathcal{B}$ and $\mathcal{A}$.
Therefore, for every $N\in \mathcal{N}$ and every $i$,
we have $\pi_{\mathcal{N}}(H^i_{\mathcal{B}}(N))= H^i_{\mathcal{A}}(\pi_{\mathcal{N}}(N))=0$, so $H^i_{\mathcal{B}}(N)\in \mathcal{N}$.
By Lemma \ref{lem:quotient-heart-criterion},
$\mathcal{B}$ descends to $\mathcal{B}/(\mathcal{B}\cap \mathcal{N})=\mathcal{B}/\mathcal{H_N} \in \heartsuit^b(\mathcal{D/N})$.

Finally, $\pi_{\mathcal{N}}(\mathcal{B})$ and $\mathcal{A}$ are both hearts of bounded $t$-structures on $\mathcal{D/N}$ and $\pi_{\mathcal{N}}(\mathcal{B})\subset \mathcal{A}$,
so we conclude that $\mathcal{A}=\pi_{\mathcal{N}}(\mathcal{B})=\mathcal{B}/\mathcal{H_N}$.
\end{proof}

\paragraph{\textbf{Uniqueness of descended lifts.}}

The following proposition shows why this definition is reasonable.

\begin{proposition}\label{prop:uniqueness-of-descended-lifts}
Let $\mathcal{B}\in\heartsuit^{b}(\mathcal{D})$ and $\mathcal{A}\in\heartsuit^{b}(\mathcal{D}/\mathcal{N})$.
Assume that $\pi_{\mathcal{N}}(\mathcal{B})\subset \mathcal{A}$. Then $\mathcal{B}\cap \mathcal{N}$ is a heart of a bounded $t$-structure on $\mathcal{N}$, and $\mathcal{B}=\Per(\mathcal{D},\mathcal{N};\mathcal{A},\mathcal{B}\cap \mathcal{N})$.
\end{proposition}

\begin{proof}
Since $\pi_{\mathcal{N}}$ is exact and $\pi_{\mathcal{N}}(\mathcal{B})\subset \mathcal{A}$, it follows that $\pi_{\mathcal{N}}$ is t-exact with respect to the bounded $t$-structures corresponding to $\mathcal{B}$ and $\mathcal{A}$.

For $N \in \mathcal{N}$, since $\pi_{\mathcal{N}}(N)=0$,
we have $\pi_{\mathcal{N}}(H^i_{\mathcal{B}}(N))\cong H^i_{\mathcal{A}}(\pi_{\mathcal{N}}(N))=0$ for all $i$,
i.e., $H^i_{\mathcal{B}}(N)\in \mathcal{N}$ for every $i$.
It follows that $\mathcal{B}$ restricts to a heart $\mathcal{B}\cap \mathcal{N}$ of a bounded $t$-structure on $\mathcal{N}$.

We now prove that $\mathcal{B}=\Per(\mathcal{D},\mathcal{N};\mathcal{A},\mathcal{B}\cap \mathcal{N})$.
For $E\in \mathcal{B}$, we have $\pi_{\mathcal{N}}(E)\in \mathcal{A}$. Moreover, as $\mathcal{N}_{\mathcal{B}\cap \mathcal{N}}^{\geq 1}\subset \mathcal{D}_{\mathcal{B}}^{\geq 1}$ and $\mathcal{N}_{\mathcal{B}\cap \mathcal{N}}^{\leq -1}\subset \mathcal{D}_{\mathcal{B}}^{\leq -1}$, 
we have $\Hom_{\mathcal{D}}(E,\mathcal{N}_{\mathcal{B}\cap \mathcal{N}}^{\geq 1})=0$ and $\Hom_{\mathcal{D}}(\mathcal{N}_{\mathcal{B}\cap \mathcal{N}}^{\leq -1},E)=0$.
Therefore $E\in \Per(\mathcal{D},\mathcal{N};\mathcal{A},\mathcal{B}\cap \mathcal{N})$.

Conversely, let $E\in \Per(\mathcal{D},\mathcal{N};\mathcal{A},\mathcal{B}\cap \mathcal{N})$. Consider the truncation triangles
$\tau_{\mathcal{B}}^{\leq 0}E\longrightarrow E\longrightarrow \tau_{\mathcal{B}}^{\geq 1}E$
and
$\tau_{\mathcal{B}}^{\leq -1}E\longrightarrow E\longrightarrow \tau_{\mathcal{B}}^{\geq 0}E$.
Since $\pi_{\mathcal{N}}(E)\in \mathcal{A}$ and $\pi_{\mathcal{N}}$ is t-exact, we have
$\pi_{\mathcal{N}}(\tau_{\mathcal{B}}^{\geq 1}E)=0=\pi_{\mathcal{N}}(\tau_{\mathcal{B}}^{\leq -1}E)$.
Hence both $\tau_{\mathcal{B}}^{\geq 1}E$ and $\tau_{\mathcal{B}}^{\leq -1}E$ lie in $\mathcal{N}$. Since moreover
$\tau_{\mathcal{B}}^{\geq 1}E\in \mathcal{D}_{\mathcal{B}}^{\geq 1}$ and
$\tau_{\mathcal{B}}^{\leq -1}E\in \mathcal{D}_{\mathcal{B}}^{\leq -1}$,
we get
$\tau_{\mathcal{B}}^{\geq 1}E\in \mathcal{N}_{\mathcal{B}\cap \mathcal{N}}^{\geq 1}$
and
$\tau_{\mathcal{B}}^{\leq -1}E\in \mathcal{N}_{\mathcal{B}\cap \mathcal{N}}^{\leq -1}$.
This gives
$\Hom_{\mathcal{D}}(E,\tau_{\mathcal{B}}^{\geq 1}E)=0$
and
$\Hom_{\mathcal{D}}(\tau_{\mathcal{B}}^{\leq -1}E,E)=0$.
We then conclude that
$\tau_{\mathcal{B}}^{\geq 1}E=0=\tau_{\mathcal{B}}^{\leq -1}E$,
that is, $E\in \mathcal{B}$.
\end{proof}
    
\paragraph{\textbf{Existence of perverse heart lifts.}}

Our main result of this section is a useful criterion for when
$\Per(\mathcal{D},\mathcal{N};\mathcal{A},\mathcal{H_N})$ is the heart of a bounded $t$-structure.
Before proving Theorem \ref{thm:perverse-heart-lift-criterion}, we need the following preparatory lemma.

\begin{lemma}\label{lem:peeling}
Let $\mathcal{H_N}\in\heartsuit^{b}(\mathcal{N})$ and set $\gls{Wwindow}:={}^{\perp}(\mathcal{N}^{\geq q+1})\cap (\mathcal{N}^{\leq q-1})^\perp$.

If the pairs $(\mathcal{N}^{\leq 0},(\mathcal{N}^{\leq 0})^\perp)$ and $({}^\perp(\mathcal{N}^{\geq 0}),\mathcal{N}^{\geq 0})$ are torsion pairs in $\mathcal{D}$,
then for any $E\in \mathcal{D}$ and any integer $q$, there exist objects $B_0\in \mathcal{D}$ and $W\in \mathcal{W}_q$, and morphisms
$u \colon B_0\longrightarrow E$ and $v \colon B_0\longrightarrow W$ with $\cone(u) \in \mathcal{N}^{\geq q+1}$ and $\cone(v)\in \mathcal{N}^{\leq q-2}$.
In particular, $\pi_\mathcal{N}(E)\simeq\pi_\mathcal{N}(W)$ in $\mathcal{D/N}$.
\end{lemma}

\begin{proof}
By shifting, the two torsion pairs in the assumption of the lemma give, for every integer $q$, the torsion pairs $(\mathcal{N}^{\leq q},(\mathcal{N}^{\leq q})^\perp)$ and $({}^\perp(\mathcal{N}^{\geq q}),\mathcal{N}^{\geq q})$.

For $E\in \mathcal{D}$, use the second torsion pair, shifted to $q+1$, to obtain a distinguished triangle
$B_0\xlongrightarrow{u} E\longrightarrow B_1$
with
$B_0\in {}^\perp(\mathcal{N}^{\geq q+1})$
and
$B_1\in \mathcal{N}^{\geq q+1}$.
Next, apply the first torsion pair, shifted to $q-1$, to $B_0$ to obtain a distinguished  triangle
$B_2\longrightarrow B_0\xlongrightarrow{v} W$
with
$B_2\in \mathcal{N}^{\leq q-1}$
and
$W\in (\mathcal{N}^{\leq q-1})^\perp$.

It remains to show that $W\in \mathcal{W}_q={}^{\perp}(\mathcal{N}^{\geq q+1})\cap (\mathcal{N}^{\leq q-1})^\perp$.
We apply $\Hom_{\mathcal{D}}(-,T)$ to the distinguished triangle
$B_2\longrightarrow B_0\longrightarrow W\longrightarrow B_2[1]$ for some $T\in \mathcal{N}^{\geq q+1}$.
Since $\Hom_{\mathcal{D}}(B_0,T)=0$ and $\Hom_{\mathcal{D}}(B_2[1],T)=0$,
we see that $\Hom_{\mathcal{D}}(W,T)=0$.

Finally, since $\cone(u)=B_1, \cone(v)=B_2[1]\in \mathcal{N}$,
the roof $E\xlongleftarrow{u} B_0\xlongrightarrow{v} W$
defines an isomorphism
$\pi_\mathcal{N}(E)\xlongrightarrow{\sim}\pi_\mathcal{N}(W)$
in $\mathcal{D/N}$.
\end{proof}

\begin{theorem}\label{thm:perverse-heart-lift-criterion}
Let $\mathcal{H_N}\in\heartsuit^{b}(\mathcal{N})$, and assume that $\heartsuit^{b}(\mathcal{D}/\mathcal{N})\neq\varnothing$.
Then the following conditions are equivalent:

\begin{enumerate}
\item The pairs $(\mathcal{N}^{\leq 0},(\mathcal{N}^{\leq 0})^\perp)$ and $({}^\perp(\mathcal{N}^{\geq 0}),\mathcal{N}^{\geq 0})$ are torsion pairs in $\mathcal{D}$.
\item For every $\mathcal{A}\in\heartsuit^{b}(\mathcal{D}/\mathcal{N})$, $\Per(\mathcal{D},\mathcal{N};\mathcal{A},\mathcal{H_N})$ is a heart of a bounded $t$-structure on $\mathcal{D}$.
\item There exists $\mathcal{A}\in\heartsuit^{b}(\mathcal{D}/\mathcal{N})$ such that, for every $n\in\mathbb{Z}$, $\Per(\mathcal{D},\mathcal{N};\mathcal{A}[n],\mathcal{H_N})$ is a heart of a bounded $t$-structure on $\mathcal{D}$.
\end{enumerate}
\end{theorem}
\begin{proof}
(i) $\Rightarrow$ (ii):
Fix a heart $\mathcal{A}\in\heartsuit^b(\mathcal{D/N})$.
We set $\mathcal{D}_{\mathcal{A},\mathcal{H_N}}^{\leq n}:=\mathcal{D}_{\mathcal{A},\mathcal{H_N}}^{\leq 0}[-n]$,
$\mathcal{D}_{\mathcal{A},\mathcal{H_N}}^{\geq n}:=\mathcal{D}_{\mathcal{A},\mathcal{H_N}}^{\geq 0}[-n]$.
It is clear that $\mathcal{D}_{\mathcal{A},\mathcal{H_N}}^{\leq 0}\subset \mathcal{D}_{\mathcal{A},\mathcal{H_N}}^{\leq 1}$, and
$\mathcal{D}_{\mathcal{A},\mathcal{H_N}}^{\geq 1}\subset \mathcal{D}_{\mathcal{A},\mathcal{H_N}}^{\geq 0}$.

We now prove that 
$\Hom_{\mathcal{D}}(\mathcal{D}_{\mathcal{A},\mathcal{H_N}}^{\leq 0},\mathcal{D}_{\mathcal{A},\mathcal{H_N}}^{\geq 1})=0$.
Take a morphism $\phi\colon F\longrightarrow G$ with $F\in \mathcal{D}_{\mathcal{A},\mathcal{H_N}}^{\leq 0}$ and 
$G\in \mathcal{D}_{\mathcal{A},\mathcal{H_N}}^{\geq 1}$.

Since $\pi_\mathcal{N}(\phi)=0$ in $\mathcal{D/N}$, 
there exists some $K\in\mathcal{N}$ such that $\phi$ factors as $\phi_2 \circ \phi_1$ for some $\phi_2\colon K\longrightarrow G$ and $\phi_1 \colon F \longrightarrow K$.

Consider the triangle
$\tau_{\mathcal{H_N}}^{\leq 0}K\xlongrightarrow{i}K\longrightarrow \tau_{\mathcal{H_N}}^{\geq 1}K$.
As $\Hom_{\mathcal{D}}(F,\tau_{\mathcal{H_N}}^{\geq 1}K)=0$, 
$\phi_1$ factors through $i$, say
$\phi_1=i\circ \phi_3$
for some $\phi_3 \colon F\longrightarrow \tau_{\mathcal{H_N}}^{\leq 0}K$.

On the other hand, since $\Hom_{\mathcal{D}}(\tau_{\mathcal{H_N}}^{\leq 0}K,G)=0$,
we have $\phi_2 \circ i=0$,
and hence $\phi=\phi_2\circ \phi_1=\phi_2 \circ i \circ \phi_3 =0$.
This proves the orthogonality.

Next, we prove that the truncation triangle exists, that is,
for every $E\in \mathcal{D}$,
there exists a distinguished triangle $E_0 \longrightarrow E\longrightarrow E_1$
such that $E_0\in \mathcal{D}_{\mathcal{A},\mathcal{H_N}}^{\leq 0}$
and $E_1\in \mathcal{D}_{\mathcal{A},\mathcal{H_N}}^{\geq 1}$.
To begin with, we claim that the truncation triangle exists for every $E\in \mathcal{W}_0:={}^\perp(\mathcal{N}^{\geq 1})\cap (\mathcal{N}^{\leq -1})^\perp$.

Consider the truncation triangle
$\tau_\mathcal{A}^{\leq 0} \pi_\mathcal{N}(E)\xlongrightarrow{\alpha} \pi_\mathcal{N}(E)\longrightarrow \tau_\mathcal{A}^{\geq 1} \pi_\mathcal{N}(E)$.
We choose $\widetilde{A}_0 \in \mathcal{D}$ with $\pi_\mathcal{N}(\widetilde{A}_0)=\tau_\mathcal{A}^{\leq 0}\pi_\mathcal{N}(E)$
and represent $\alpha$ by the roof $\widetilde{A}_0 \xlongleftarrow{s} X \xlongrightarrow{a} E$.

Applying Lemma \ref{lem:peeling} with $q=0$ to $X$,
there are morphisms
$b\colon B\longrightarrow X$ and $c\colon B\longrightarrow W$
such that $\cone(b)\in \mathcal{N}^{\geq 1}, \cone(c)\in \mathcal{N}^{\leq -2}$, and $\pi_\mathcal{N}(W)\simeq \pi_\mathcal{N}(X)$.

Consider the composition $a\circ b\colon B\longrightarrow E$.
Since $E\in\mathcal{W}_0\subseteq(\mathcal{N}^{\leq -1})^\perp$, we have $\Hom(\cone(c)[-1],E)=0$.
Applying $\Hom(-,E)$ to the distinguished triangle $\cone(c)[-1]\longrightarrow B\xlongrightarrow{c}W\longrightarrow\cone(c)$, we obtain a morphism $f\colon W\longrightarrow E$ with $f\circ c=a\circ b$,
and hence $\pi_{\mathcal{N}}(f)=\alpha$ in $\mathcal{D/N}$, under the isomorphism $\pi_{\mathcal{N}}(s\circ b)\circ\pi_{\mathcal{N}}(c)^{-1}\colon\pi_{\mathcal{N}}(W)\xlongrightarrow{\sim}\tau_{\mathcal{A}}^{\leq 0}\pi_{\mathcal{N}}(E)$.

Therefore, we may choose a triangle
$A_0\longrightarrow E\longrightarrow A_1$
such that $A_0\in \mathcal{W}_0$,
$\pi_\mathcal{N}(A_0)= \tau_\mathcal{A}^{\leq 0} \pi_\mathcal{N}(E)$,
and
$\pi_\mathcal{N}(A_1) = \tau_\mathcal{A}^{\geq 1} \pi_\mathcal{N}(E)$.

Applying the first torsion pair in (i) to $A_1$, we obtain a distinguished triangle
$A_2 \longrightarrow A_1\longrightarrow E_1$
with $A_2\in \mathcal{N}^{\leq 0}$ and $E_1\in (\mathcal{N}^{\leq 0})^\perp$.
Applying the octahedral axiom to the composition $(E \longrightarrow A_1 \longrightarrow E_1)$,
we obtain a triangle
$E_0 \longrightarrow E\longrightarrow E_1$
where $E_0$ fits into the distinguished triangle $A_0 \longrightarrow E_0 \longrightarrow A_2$.

Since $\pi_\mathcal{N}(A_2)=0$, we have
$\pi_\mathcal{N}(E_1)= \pi_\mathcal{N}(A_1)=\tau_\mathcal{A}^{\geq 1} \pi_\mathcal{N}(E)$,
and hence
$E_1\in \mathcal{D}_{\mathcal{A},\mathcal{H_N}}^{\geq 1}$.
On the other hand, $A_0 \in \mathcal{W}_0$ and
$\pi_\mathcal{N}(A_0) =  \tau_\mathcal{A}^{\leq 0} \pi_\mathcal{N}(E)\in (\mathcal{D/N})_{\mathcal{A}}^{\leq 0}$,
so $A_0 \in \mathcal{D}_{\mathcal{A},\mathcal{H_N}}^{\leq 0}$.
Also, $A_2 \in \mathcal{N}^{\leq 0}\subset \mathcal{D}_{\mathcal{A},\mathcal{H_N}}^{\leq 0}$.
Therefore $E_0 \in \mathcal{D}_{\mathcal{A},\mathcal{H_N}}^{\leq 0}$.
This proves the claim.

Now for arbitrary $E\in \mathcal{D}$,
Lemma \ref{lem:peeling} with $q=0$ applied to $E$ gives $u\colon B_0\longrightarrow E$ and $v\colon B_0\longrightarrow W$ with $\cone(u) \in \mathcal{N}^{\geq 1}$ and $\cone(v) \in \mathcal{N}^{\leq -2}$.
By our claim above there exists a truncation triangle
$W_0\longrightarrow W\longrightarrow W_1\longrightarrow$
such that
$W_0\in \mathcal{D}_{\mathcal{A},\mathcal{H_N}}^{\leq 0}$
and
$W_1\in \mathcal{D}_{\mathcal{A},\mathcal{H_N}}^{\geq 1}$.

Applying the octahedral axiom to the composition
$(B_0\longrightarrow W\longrightarrow W_1)$,
we obtain a triangle
$E_0\longrightarrow B_0\longrightarrow W_1$
where $E_0$ fits into a distinguished triangle
$\cone(v)[-1]\longrightarrow E_0 \longrightarrow W_0$.
Applying the octahedral axiom again to the composition
$(E_0\longrightarrow B_0\longrightarrow E)$,
we obtain a triangle
$E_0\longrightarrow E\longrightarrow E_1$
where $E_1$ fits into a distinguished triangle
$W_1\longrightarrow E_1\longrightarrow \cone(u)$.

Now we have
$W_0\in \mathcal{D}_{\mathcal{A},\mathcal{H_N}}^{\leq 0}$,
$\cone(v)[-1]\in \mathcal{N}^{\leq -1}\subset \mathcal{D}_{\mathcal{A},\mathcal{H_N}}^{\leq 0}$, 
$W_1\in \mathcal{D}_{\mathcal{A},\mathcal{H_N}}^{\geq 1}$,
and
$\cone(u)\in \mathcal{N}^{\geq 1}\subset \mathcal{D}_{\mathcal{A},\mathcal{H_N}}^{\geq 1}$,
so we conclude that
$E_0\in \mathcal{D}_{\mathcal{A},\mathcal{H_N}}^{\leq 0}$
and that $E_1\in \mathcal{D}_{\mathcal{A},\mathcal{H_N}}^{\geq 1}$.
This proves the third axiom of a $t$-structure.

Finally we prove that the $t$-structure is bounded,
that is, for any $E\in\mathcal{D}$, there is an interval $[a,b]$ so that $E\in\mathcal{D}^{[a,b]}_{\mathcal{A},\mathcal{H_N}}$.

Let $E\in\mathcal{D}$. Since $\mathcal{A}$ is a heart of a bounded $t$-structure, there exists an interval $[a_0,b_0]$ such that $\pi_\mathcal{N}(E)\in (\mathcal{D/N})_{\mathcal{A}}^{[a_0,b_0]}$.

The first torsion pair gives a distinguished triangle $K_1\longrightarrow E\longrightarrow L_1$ with $K_1\in\mathcal{N}^{\leq 0}$ and $L_1\in(\mathcal{N}^{\leq 0})^\perp$.
Since the $t$-structure on $\mathcal{N}$ is bounded, we may choose $a_1\leq 0$ such that $K_1\in\mathcal{N}^{\geq a_1}$.
For any $M\in\mathcal{N}^{\leq a_1-1}$, we have $\Hom_{\mathcal{D}}(M,K_1)=0$ and $\Hom_{\mathcal{D}}(M,L_1)=0$ since $M\in\mathcal{N}^{\leq 0}$.
Applying $\Hom_{\mathcal{D}}(M,-)$ to the triangle gives $\Hom_{\mathcal{D}}(M,E)=0$. Thus $E\in(\mathcal{N}^{\leq a_1-1})^\perp$.

Similarly, the second torsion pair gives $K_2\longrightarrow E\longrightarrow L_2$ with $K_2\in{}^\perp(\mathcal{N}^{\geq 0})$ and $L_2\in\mathcal{N}^{\geq 0}$.
Choose $b_1\geq 0$ such that $L_2\in\mathcal{N}^{\leq b_1}$.
For any $M\in\mathcal{N}^{\geq b_1+1}$, both $\Hom_{\mathcal{D}}(K_2,M)$ and $\Hom_{\mathcal{D}}(L_2,M)$ vanish.
Hence $\Hom_{\mathcal{D}}(E,M)=0$, and $E\in{}^\perp(\mathcal{N}^{\geq b_1+1})$.

If we set $a:=\min\{a_0,a_1\}$ and $b:=\max\{b_0,b_1\}$,
then $\pi_\mathcal{N}(E)\in (\mathcal{D/N})_{\mathcal{A}}^{[a,b]}$, $\Hom_{\mathcal{D}}(E,\mathcal{N}^{\geq b+1})=0$, and $\Hom_{\mathcal{D}}(\mathcal{N}^{\leq a-1},E)=0$,
so $E\in\mathcal{D}_{\mathcal{A},\mathcal{H_N}}^{[a,b]}$.

(ii) $\Rightarrow$ (iii) is trivial.
Finally, for (iii) $\Rightarrow$ (i), write $\mathcal{B}_n:=\Per(\mathcal{D},\mathcal{N};\mathcal{A}[n],\mathcal{H_N})$.
We will prove that $(\mathcal{N}^{\leq 0}_{\mathcal{H_N}},(\mathcal{N}^{\leq 0}_{\mathcal{H_N}})^\perp)$ is a torsion pair in $\mathcal{D}$,
and the proof for the second pair is dual.
Given $E\in\mathcal{D}$, by boundedness of the $t$-structure with heart $\mathcal{A}$,
we may choose $n$ sufficiently large so that $\pi_{\mathcal{N}}(E)\in(\mathcal{D/N})^{\geq1}_{\mathcal{A}[n]}$.

Now consider the truncation triangle $\tau_{\mathcal{B}_n}^{\leq 0}E\longrightarrow E\longrightarrow\tau_{\mathcal{B}_n}^{\geq1}E$.
By Proposition \ref{prop:formalprop}, $\pi_{\mathcal{N}}(\tau_{\mathcal{B}_n}^{\leq 0}E)=\tau_{\mathcal{A}[n]}^{\leq 0}\pi_{\mathcal{N}}(E)=0$,
and hence $\tau_{\mathcal{B}_n}^{\leq 0}E\in \mathcal{D}^{\leq 0}_{\mathcal{B}_n} \cap \mathcal{N}=\mathcal{N}^{\leq 0}_{\mathcal{H_N}}$.
Moreover, $\mathcal{N}^{\leq 0}_{\mathcal{H_N}}\subseteq\mathcal{D}^{\leq 0}_{\mathcal{B}_n}$, so $\tau_{\mathcal{B}_n}^{\geq1}E\in (\mathcal{D}^{\leq 0}_{\mathcal{B}_n})^\perp \subseteq (\mathcal{N}^{\leq 0}_{\mathcal{H_N}})^\perp$.
This completes the proof.
\end{proof}

In such a heart there is a core part which is equivalent to the heart downstairs.

\begin{proposition}\label{prop:core-equivalence}
Assume that the equivalent conditions of Theorem \ref{thm:perverse-heart-lift-criterion} hold,
and let $\mathcal{A}\in\heartsuit^{b}(\mathcal{D}/\mathcal{N})$. Then:
\begin{enumerate}
\item $(\mathcal{H_N},\mathcal{H_N}^{\perp})$ and $({}^{\perp}\mathcal{H_N},\mathcal{H_N})$ are torsion pairs in $\Per(\mathcal{D},\mathcal{N};\mathcal{A},\mathcal{H_N})$.

\item $\pi_\mathcal{N}|_{\mathcal{H_N}^{\perp}\cap {}^{\perp}\mathcal{H_N}}\colon \mathcal{H_N}^{\perp}\cap {}^{\perp}\mathcal{H_N}\longrightarrow \mathcal{A}$ is an equivalence of categories.
\end{enumerate}
\end{proposition}

\begin{proof}
For simplicity we set $\mathcal{B}:=\Per(\mathcal{D},\mathcal{N};\mathcal{A},\mathcal{H_N})$.
By Theorem \ref{thm:perverse-heart-lift-criterion} and Proposition \ref{prop:formalprop},
$\mathcal{B}$ is a heart of a bounded $t$-structure and $\mathcal{A}=\mathcal{B}/\mathcal{H_N}$.

For (i), let $E\in\mathcal{B}$.
The first torsion pair in Theorem \ref{thm:perverse-heart-lift-criterion}(i) gives a distinguished triangle $T\longrightarrow E\longrightarrow F$ with $T\in\mathcal{N}^{\leq 0}$ and $F\in(\mathcal{N}^{\leq 0})^\perp$.
We have $\pi_\mathcal{N}(F)=\pi_\mathcal{N}(E)\in\mathcal{A}$.
For $K\in\mathcal{N}^{\geq1}$, both $\Hom_{\mathcal{D}}(E,K)$ and $\Hom_{\mathcal{D}}(T[1],K)$ vanish, so $\Hom_{\mathcal{D}}(F,K)=0$.
Hence $F\in\mathcal{B}$ by the definition of $\Per(\mathcal{D},\mathcal{N};\mathcal{A},\mathcal{H_N})$.
Since $E,F\in\mathcal{B}$, the triangle shows that $T$ has $\mathcal{B}$-cohomology only in degrees $0,1$.
On the other hand, $T\in\mathcal{N}^{\leq 0}\subseteq\mathcal{D}_{\mathcal{B}}^{\leq 0}$.
Thus $T\in\mathcal{B}\cap\mathcal{N}=\mathcal{H_N}$, and we obtain a short exact sequence $0\longrightarrow T\longrightarrow E\longrightarrow F\longrightarrow0$ in $\mathcal{B}$ with $F\in\mathcal{H_N}^{\perp}$.

For the second pair, the other ambient torsion pair gives a distinguished triangle $F\longrightarrow E\longrightarrow T$ with $F\in{}^\perp(\mathcal{N}^{\geq 0})$ and $T\in\mathcal{N}^{\geq 0}$.
Similarly, we have $\pi_\mathcal{N}(F)=\pi_\mathcal{N}(E)\in\mathcal{A}$.
For $K\in\mathcal{N}^{\leq-1}$, both $\Hom_{\mathcal{D}}(K,T[-1])$ and $\Hom_{\mathcal{D}}(K,E)$ vanish, and hence $\Hom_{\mathcal{D}}(K,F)=0$.
It follows that $F\in\mathcal{B}$.
Now $T$ has $\mathcal{B}$-cohomology only in degrees $-1,0$, while $T\in\mathcal{N}^{\geq 0}\subseteq\mathcal{D}_{\mathcal{B}}^{\geq 0}$.
Thus $T\in\mathcal{H_N}$ and we obtain the required short exact sequence in $\mathcal{B}$.

Now we prove (ii).
For $E\in \mathcal{B}$, by (i), there are short exact sequences
$0\longrightarrow T_1\longrightarrow E\longrightarrow F_1\longrightarrow 0$
and
$0\longrightarrow F_2\longrightarrow F_1\longrightarrow T_2\longrightarrow 0$
with $T_i\in \mathcal{H_N}$, $F_1\in \mathcal{H_N}^{\perp}$, and $F_2\in {}^{\perp}\mathcal{H_N}$.
Since $\mathcal{H_N}^{\perp}$ is closed under subobjects, we have $F_2\in \mathcal{H_N}^{\perp}\cap {}^{\perp}\mathcal{H_N}$.
Setting $c(E):=F_2$,
we define a functor $c\colon \mathcal{B}\longrightarrow \mathcal{H_N}^{\perp}\cap {}^{\perp}\mathcal{H_N}$ with $\pi_\mathcal{N}(c(E))=\pi_\mathcal{N}(E)$.

If $f\colon E\longrightarrow E'$ is a morphism in $\mathcal{B}$ such that $\pi_\mathcal{N}(f)$ is an isomorphism, then $\pi_\mathcal{N}(c(f))$ is also an isomorphism,
that is, both $\ker c(f)$ and $\cok c(f)$ lie in $\mathcal{H_N}$.
Since $c(E)\in \mathcal{H_N}^{\perp}$ and $c(E')\in {}^{\perp}\mathcal{H_N}$, this forces $\ker c(f)=0=\cok c(f)$, so $c(f)$ is an isomorphism.
Therefore $c$ factors through the Serre quotient $\mathcal{A}=\mathcal{B}/\mathcal{H_N}$ and gives a functor $\overline{c}\colon \mathcal{A}\longrightarrow \mathcal{H_N}^{\perp}\cap {}^{\perp}\mathcal{H_N}$ such that $\overline{c}\circ \pi_\mathcal{N}\simeq c$.

One may verify that $\overline{c}$ is a quasi-inverse of $\pi_\mathcal{N}|_{\mathcal{H_N}^{\perp}\cap {}^{\perp}\mathcal{H_N}}$.
This proves the assertion.
\end{proof}

\begin{remark}\label{lem:window-quotient}
Assume that the equivalent conditions of Theorem \ref{thm:perverse-heart-lift-criterion} hold,
and let $\mathcal{A}\in\heartsuit^{b}(\mathcal{D}/\mathcal{N})$. Then:
\begin{enumerate}
    \item Set $	\gls{Vwindow}:={}^{\perp}(\mathcal{N}^{\geq q})\cap (\mathcal{N}^{\leq q})^\perp.$ 
    Then the hypotheses of \cite[Theorem 1.1]{IY20} are satisfied,
    and we obtain that $\pi_\mathcal{N}$ induces an equivalence of additive categories $\mathcal{V}_q\xlongrightarrow{\sim}\mathcal{D/N}$.

    Note that $\mathcal{V}_0 \cap \pi^{-1}_\mathcal{N}(\mathcal{A}) = \mathcal{H_N}^{\perp}\cap {}^{\perp}\mathcal{H_N}$, and the equivalence in Proposition \ref{prop:core-equivalence} can also be viewed as the restriction of this one.
    \item $\pi_\mathcal{N}$ also induces an additive functor $\mathcal{W}_q/[\mathcal{H_N}[-q]]\longrightarrow\mathcal{D/N}$,
    which is essentially surjective by Lemma \ref{lem:peeling}.

    Let $f\in \Hom_{\mathcal{W}_0}(X,Y)$ be a morphism such that $\pi_\mathcal{N}(f)=0$.
    Then $f$ factors through some $K\in \mathcal{N}$.
    Since $X\in {}^\perp(\mathcal{N}^{\geq 1})$, $f$ factors through $\tau_{\mathcal{N}}^{\leq 0}K$.
    Consider the truncation triangle
    $\tau_{\mathcal{N}}^{\leq -1}K\longrightarrow \tau_{\mathcal{N}}^{\leq 0}K\longrightarrow H^0_{\mathcal{N}}(K)$.
    Since $Y\in (\mathcal{N}^{\leq -1})^\perp$,
    we see that $f$ factors through $H^0_{\mathcal{N}}(K)\in \mathcal{H_N}$,
    and hence $\bar{f}=0$ in $\mathcal{W}_0/[\mathcal{H_N}]$.
    That is, $\bar{\pi}_\mathcal{N}$ is faithful.

    However, such $\bar{\pi}_\mathcal{N}$ is not full in general.
\end{enumerate}
\end{remark}

\paragraph{\textbf{Verifying the torsion-pair condition.}}

We next give several sufficient criteria to verify the torsion-pair condition in Theorem \ref{thm:perverse-heart-lift-criterion}(i).

\begin{theorem}\label{thm:module-type-heart}
Assume that $\mathcal{D}$ is proper, and let $\mathcal{H_N}\in\heartsuit^{b}(\mathcal{N})$ with $\mathcal{H_N}\simeq\mathrm{mod}\text{-}\Gamma$ for some finite-dimensional $k$-algebra $\Gamma$.
Then $(\mathcal{N}^{\leq q},(\mathcal{N}^{\leq q})^\perp)$ and $({}^\perp(\mathcal{N}^{\geq q}),\mathcal{N}^{\geq q})$ are torsion pairs in $\mathcal{D}$.

\end{theorem}

Before proving Theorem \ref{thm:module-type-heart}, we first prove two preparatory lemmas under its assumptions.
In these lemmas, we fix a projective generator $P$ and an injective cogenerator $I$ in $\mathcal{H_N}\simeq\mathrm{mod}\text{-}\Gamma$.

\begin{lemma}\label{lem:detection}
For any integer $m$ and any object $E\in \mathcal{D}$, we have:
\begin{enumerate}
    \item $\Hom_{\mathcal{D}}(P[t],E)=0$  for all $t\geq m$ if and only if $ E\in (\mathcal{N}^{\leq -m})^\perp.$
    \item $\Hom_{\mathcal{D}}(E,I[t])=0$  for all  $t\leq m$ if and only if $ E\in {}^\perp(\mathcal{N}^{\geq -m})$.
\end{enumerate}
\end{lemma}
\begin{proof}
We prove (i) and then (ii) is similar.

One direction is trivial. Now assume that $\Hom_{\mathcal{D}}(P[t],E)=0$ for all $t\geq m$.
Let $S_1,\dots,S_r$ be the simple objects of $\mathcal{H_N}$.
Since $P$ is a projective generator, for each $i$ there exists a short exact sequence $0\longrightarrow K_i\longrightarrow P^{\oplus a_i}\longrightarrow S_i\longrightarrow 0$ in $\mathcal{H_N}$.

We claim that $\Hom_{\mathcal{D}}(S_i[t],E)=0$ for all $i$ and all $t\geq m$.
For $t\gg 0$ this is clear by properness of $\mathcal{D}$.

For $t \geq m$, suppose that $\Hom_{\mathcal{D}}(S_i[t+1],E)=0$ for all $i$.
Then as $\mathcal{H_N}$ is of finite length, 
we have $\Hom_{\mathcal{D}}(K_i[t+1],E)=0$.
Applying $\Hom_{\mathcal{D}}(-,E)$ to the exact triangle $P^{\oplus a_i}[t]\longrightarrow S_i[t]\longrightarrow K_i[t+1]$, we obtain an exact sequence $\Hom_{\mathcal{D}}(K_i[t+1],E)\longrightarrow \Hom_{\mathcal{D}}(S_i[t],E)\longrightarrow \Hom_{\mathcal{D}}(P^{\oplus a_i}[t],E)$,
and hence we have $\Hom_{\mathcal{D}}(S_i[t],E)=0$.

The descending induction proves the claim,
and hence $\Hom_{\mathcal{D}}(N,E)=0$ for all $N\in \mathcal{N}^{\leq -m}$, that is, $E\in (\mathcal{N}^{\leq -m})^\perp$.
\end{proof}

\begin{lemma}\label{lem:approximation}
For any integer $r$, the following statements hold.
\begin{enumerate}
\item If $E\in (\mathcal{N}^{\leq r-1})^\perp$,
then there exists a distinguished triangle
\[
C_r(E)\longrightarrow E\longrightarrow E_r'\longrightarrow C_r(E)[1]
\]
such that  $C_r(E)\in \mathcal{N}^{\leq r}$ and $E_r'\in (\mathcal{N}^{\leq r})^\perp$.

\item If $E\in {}^\perp(\mathcal{N}^{\geq r+1})$,
then there exists a distinguished  triangle
\[
E_r''\longrightarrow E\longrightarrow D_r(E)\longrightarrow E_r''[1]
\]
such that $D_r(E)\in \mathcal{N}^{\geq r}$ and 
$E_r'' \in {}^\perp(\mathcal{N}^{\geq r})$.
\end{enumerate}
\end{lemma}

\begin{proof}
We first prove (i).

Set $R:=\End_{\mathcal{H_N}}(P)$.
Since $P$ is a projective generator of $\mathcal{H_N}$, there is an exact equivalence
$\Hom_{\mathcal{H_N}}(P,-):\mathcal{H_N}\xlongrightarrow{\sim}\mathrm{mod}\text{-}R$
with exact quasi-inverse
$-\otimes_R P.$
For each $M \in \mathrm{mod}\text{-}R$, let
$\beta_M\colon M\xlongrightarrow{\sim}\Hom_{\mathcal{H_N}}(P,M\otimes_R P)$
be the natural isomorphism given by this equivalence.

Let $V:=\Hom_{\mathcal{D}}(P[-r],E)$.
Since $\mathcal{D}$ is proper, $V$ is a finite-dimensional $k$-vector space, hence a finitely generated right $R$-module.
We define a contravariant functor $\Phi_E \colon \mathrm{mod}\text{-}R \longrightarrow \mathrm{mod}\text{-}k$ by $\Phi_E(M):=\Hom_\mathcal{D}((M\otimes_R P)[-r],E)$,
and a natural transformation $\Theta \colon \Phi_E \longrightarrow \Hom_R(-,V)$ by 
$\Theta_M(h):=\Hom_{\mathcal{D}}(P[-r],h)\circ \beta_M$ for $h\in \Hom_{\mathcal{D}}((M\otimes_R P)[-r],E)$.

We now prove that $\Theta_M$ is an isomorphism.
First, for a free module $M=R^{\oplus m}$, both 
$\Phi_E(M)$ and $\Hom_R(M,V)$ are just $V^{\oplus m}$ and $\Theta_{M}$ is the identity map.

Next, we show that $\Phi_E$ is a left exact functor.
Indeed, consider a short exact sequence  $0\longrightarrow A\longrightarrow B\longrightarrow C\longrightarrow 0$ in $\mathrm{mod}\text{-}R$.
Since $-\otimes_R P$ is exact, 
there is a distinguished triangle
$(A\otimes_R P)[-r]\longrightarrow (B\otimes_R P)[-r]\longrightarrow (C\otimes_R P)[-r]$
in $\mathcal{D}$.
By our assumption,
$\Hom_{\mathcal{D}}((A\otimes_R P)[1-r],E)=0$,
and hence applying $\Hom_{\mathcal{D}}(-,E)$ to the above triangle shows that 
$0\longrightarrow \Phi_E(C)\longrightarrow \Phi_E(B)\longrightarrow \Phi_E(A)$ is exact.

For arbitrary $M$, since $R$ is a finite-dimensional $k$-algebra, every finitely generated right $R$-module is finitely presented.
As both $\Phi_E$ and $\Hom_R(-,V)$ are left exact, 
we can pass to a finite presentation of $M$ to see that $\Theta_M$ is an isomorphism.

Now define
$C_r(E):=(V\otimes_R P)[-r]$ and let
$u_r:=\Theta_V^{-1}(1_V)\colon C_r(E)\longrightarrow E$.
Note that this implies that $\Hom_\mathcal{D}(P[-r],u_r)$ is an isomorphism.
Complete $u_r$ to a distinguished triangle
\begin{equation}\label{approximation-triangle}
C_r(E)\xlongrightarrow{u_r}E\longrightarrow E_r'\longrightarrow C_r(E)[1].
\end{equation}

Since $V\otimes_R P\in \mathcal{H_N}\subset \mathcal{N}^{\leq 0}$, we have
$C_r(E)\in \mathcal{N}^{\leq r}$.
By Lemma \ref{lem:detection}, it remains to prove that
$\Hom_{\mathcal{D}}(P[t],E_r')=0$
for all $t\geq-r$.

Applying $\Hom_{\mathcal{D}}(P[t],-)$ to \eqref{approximation-triangle},
we obtain the following exact sequence:

\begin{tikzcd}
{\Hom_{\mathcal{D}}(P[t],C_r(E))} \arrow[r] & {\Hom_{\mathcal{D}}(P[t],E)} \arrow[r] & {\Hom_{\mathcal{D}}(P[t],E_r')} \arrow[lld] \\
{\Hom_{\mathcal{D}}(P[t],C_r(E)[1])} \arrow[r]                              & {\Hom_{\mathcal{D}}(P[t],E[1])}        &  
\end{tikzcd}

For $t>1-r$,
we have $\Hom_{\mathcal{D}}(P[t],C_r(E))=0=\Hom_{\mathcal{D}}(P[t],C_r(E)[1])$,
and hence $\Hom_{\mathcal{D}}(P[t],E_r') = \Hom_{\mathcal{D}}(P[t],E)=0$.

For $t=1-r$, by assumption, $\Hom_{\mathcal{D}}(P[1-r],E)=0$,
so $\Hom_{\mathcal{D}}(P[1-r],E_r')$ is the kernel of the morphism 
$\Hom_{\mathcal{D}}(P[1-r],u_r[1])= \Hom_{\mathcal{D}}(P[-r],u_r)$,
and hence $\Hom_{\mathcal{D}}(P[1-r],E_r')=0$.

For $t=-r$,  we have $\Hom_{\mathcal{D}}(P[-r],C_r(E)[1])=\Ext^1_{\mathcal{H_N}}(P,V\otimes_R P)=0$ since $P$ is projective in $\mathcal{H_N}$.
Therefore $\Hom_{\mathcal{D}}(P[-r],E_r')$ is the cokernel of the morphism 
$\Hom_{\mathcal{D}}(P[-r],u_r)$, which is $0$.

Combining the three cases, we conclude that
$\Hom_{\mathcal{D}}(P[t],E_r')=0$
for all $t\geq-r$.
This proves (i).

For (ii), apply (i) to $\mathcal{D}^{\mathrm{op}}$ and $\mathcal{N}^{\mathrm{op}}$, where $I$ is a projective generator of $\mathcal{H_N}^{\mathrm{op}}$, with the integer $-r$. 
\end{proof}

\begin{proof}[Proof of Theorem \ref{thm:module-type-heart}]
For the first torsion pair, $\Hom_{\mathcal{D}}(\mathcal{N}^{\leq 0},(\mathcal{N}^{\leq 0})^\perp)=0$ is clear.
It remains to show that every $E\in \mathcal{D}$ fits into a distinguished triangle $K\longrightarrow E\longrightarrow L$ with $K\in \mathcal{N}^{\leq 0}$ and $L\in (\mathcal{N}^{\leq 0})^\perp$.

If $E\in (\mathcal{N}^{\leq 0})^\perp$, we simply take $K:=0$ and $L:=E$. Assume now that $E\notin (\mathcal{N}^{\leq 0})^\perp$ and set $E_0:=E$.
Inductively, if $E_{a-1}\notin (\mathcal{N}^{\leq 0})^\perp$, consider the set
$\{\, r\leq 0 \mid E_{a-1}\notin (\mathcal{N}^{\leq r})^\perp \,\}$,
which admits a minimum, denoted by $r_a$, by the properness of $\mathcal{D}$ and Lemma \ref{lem:detection}(i).

By the minimality of $r_a$, we have $E_{a-1}\in (\mathcal{N}^{\leq r_a-1})^\perp$. Therefore Lemma \ref{lem:approximation}(i) yields a distinguished triangle
$C_{r_a}(E_{a-1})\longrightarrow E_{a-1}\longrightarrow E_a$
such that $C_{r_a}(E_{a-1})\in \mathcal{N}^{\leq r_a}\subset \mathcal{N}^{\leq 0}$ and $E_a\in (\mathcal{N}^{\leq r_a})^\perp$.

If the process continues, then by definition $r_{a+1}>r_a$.
Since every $r_a\leq 0$, this strictly increasing sequence terminates,
and we can take $L:=E_s \in (\mathcal{N}^{\leq 0})^\perp$ to be the terminal object.

Applying the octahedral axiom repeatedly, the above triangles glue to a distinguished triangle
$K\longrightarrow E\longrightarrow L$,
where $K$ is an extension of the objects $C_{r_a}(E_{a-1}) \in \mathcal{N}^{\leq 0}$ and in particular  $K\in \mathcal{N}^{\leq 0}$.

For the second torsion pair, apply the first part to $\mathcal{D}^{\mathrm{op}}$ and $\mathcal{N}^{\mathrm{op}}$.
Translating the resulting triangle back to $\mathcal{D}$ completes the proof.
\end{proof}

Under the additional assumptions of Theorem \ref{thm:module-type-heart}, the torsion decompositions in Proposition \ref{prop:core-equivalence} can also be constructed explicitly by the evaluation morphism associated with a projective generator and the coevaluation morphism associated with an injective cogenerator.

\begin{lemma}\label{lem:finite-amplitude-transfer}
Let $\mathcal{H_N},\mathcal{H'_N}\in\heartsuit^{b}(\mathcal{N})$.
Assume that $\mathcal{H_N}\subseteq \mathcal{N}_{\mathcal{H_N'}}^{[a,b]}$ for some integers $a,b$. Then the following are equivalent:
\begin{enumerate}
\item The pairs $(\mathcal{N}_{\mathcal{H_N}}^{\leq 0},(\mathcal{N}_{\mathcal{H_N}}^{\leq 0})^\perp)$ and $({}^\perp(\mathcal{N}_{\mathcal{H_N}}^{\geq 0}),\mathcal{N}_{\mathcal{H_N}}^{\geq 0})$ are torsion pairs in $\mathcal{D}$.
\item The pairs $(\mathcal{N}_{\mathcal{H_N'}}^{\leq 0},(\mathcal{N}_{\mathcal{H_N'}}^{\leq 0})^\perp)$ and $({}^\perp(\mathcal{N}_{\mathcal{H_N'}}^{\geq 0}),\mathcal{N}_{\mathcal{H_N'}}^{\geq 0})$ are torsion pairs in $\mathcal{D}$.
\end{enumerate}
In particular, if $\mathcal{H_N}$ is of finite length with finitely many simples,
and the pairs $(\mathcal{N}_{\mathcal{H_N}}^{\leq 0},(\mathcal{N}_{\mathcal{H_N}}^{\leq 0})^\perp)$ and $({}^\perp(\mathcal{N}_{\mathcal{H_N}}^{\geq 0}),\mathcal{N}_{\mathcal{H_N}}^{\geq 0})$ are torsion pairs in $\mathcal{D}$,
then we have the corresponding two torsion pairs for every heart of a bounded $t$-structure on $\mathcal{N}$.
\end{lemma}

\begin{proof}
Note that since the two $t$-structures are bounded, we have $\mathcal{N}_{\mathcal{H_N}}^{\leq 0}\subseteq \mathcal{N}_{\mathcal{H_N'}}^{\leq b}$
and $\mathcal{N}_{\mathcal{H_N}}^{\geq 0}\subseteq \mathcal{N}_{\mathcal{H_N'}}^{\geq a}$.
Taking right and left orthogonals of these inclusions respectively and shifting shows that $\mathcal{H_N'}\subseteq \mathcal{N}_{\mathcal{H_N}}^{[-b,-a]}$,
and hence we only need to prove that (i) implies (ii).

Assuming (i), we will prove that $(\mathcal{N}_{\mathcal{H_N'}}^{\leq 0},(\mathcal{N}_{\mathcal{H_N'}}^{\leq 0})^\perp)$ is a torsion pair in $\mathcal{D}$,
and the second torsion pair is obtained by the same argument.

For any $E\in\mathcal{D}$, by shifting the first torsion pair in (i), there is a distinguished triangle $T\longrightarrow E\longrightarrow F$
with $T\in\mathcal{N}_{\mathcal{H_N}}^{\leq -a}$ and $F\in(\mathcal{N}_{\mathcal{H_N}}^{\leq -a})^\perp$.
Taking the truncation triangle
$\tau_{\mathcal{H_N'}}^{\leq 0}T\longrightarrow T\longrightarrow\tau_{\mathcal{H_N'}}^{\geq 1}T$ 
and applying the octahedral axiom to the composition
$(\tau_{\mathcal{H_N'}}^{\leq 0}T\longrightarrow T\longrightarrow E)$,
we obtain a triangle
$\tau_{\mathcal{H_N'}}^{\leq 0}T\longrightarrow E\longrightarrow E_1$
where $E_1$ fits into a distinguished triangle
$\tau_{\mathcal{H_N'}}^{\geq 1}T\longrightarrow E_1\longrightarrow F$.

Now, for $K\in\mathcal{N}_{\mathcal{H_N'}}^{\leq 0}$, we have $\Hom_{\mathcal{D}}(K,\tau_{\mathcal{H_N'}}^{\geq 1}T)=0$ and $\Hom_{\mathcal{D}}(K,F)=0$ since $\mathcal{N}_{\mathcal{H_N'}}^{\leq 0}\subseteq\mathcal{N}_{\mathcal{H_N}}^{\leq -a}$.
Applying $\Hom_{\mathcal{D}}(K,-)$ to $\tau_{\mathcal{H_N'}}^{\geq 1}T\longrightarrow E_1\longrightarrow F$ gives $\Hom_{\mathcal{D}}(K,E_1)=0$.
Therefore, $E_1\in(\mathcal{N}_{\mathcal{H_N'}}^{\leq 0})^\perp$, and hence
$(\mathcal{N}_{\mathcal{H_N'}}^{\leq 0},(\mathcal{N}_{\mathcal{H_N'}}^{\leq 0})^\perp)$
is a torsion pair in $\mathcal{D}$.

The final assertion follows since, if $\mathcal{H_N}$ is of finite length with finitely many simples, then every heart $\mathcal{H_N}'$ of a bounded $t$-structure on $\mathcal{N}$ satisfies $\mathcal{H_N}'\subseteq\mathcal{N}_{\mathcal{H_N}}^{[a,b]}$ for some $a\leq b$ (see, for example, \cite[Corollary 4.5]{HW26}).
\end{proof}

Combined with Theorem \ref{thm:module-type-heart}, the lemma above shows that if $\mathcal{D}$ is proper and $\mathcal{N}$ admits a heart of a bounded $t$-structure of the form $\mathrm{mod}\text{-}\Gamma$ for some finite-dimensional $k$-algebra $\Gamma$, 
then the condition of Theorem \ref{thm:perverse-heart-lift-criterion}(i) holds for every heart of a bounded $t$-structure on $\mathcal{N}$.

\begin{proposition}\label{prop:admissibility-torsion-pairs}
Let $\mathcal{H_N}\in\heartsuit^{b}(\mathcal{N})$.
\begin{enumerate}
\item If $\mathcal{N}$ is right admissible in $\mathcal{D}$, then $(\mathcal{N}_{\mathcal{H_N}}^{\leq 0},(\mathcal{N}_{\mathcal{H_N}}^{\leq 0})^\perp)$ is a torsion pair in $\mathcal{D}$.
\item If $\mathcal{N}$ is left admissible in $\mathcal{D}$, then $({}^\perp(\mathcal{N}_{\mathcal{H_N}}^{\geq 0}),\mathcal{N}_{\mathcal{H_N}}^{\geq 0})$ is a torsion pair in $\mathcal{D}$.
\end{enumerate}
In particular, if $\mathcal{N}$ is an admissible component of a semiorthogonal decomposition of $\mathcal{D}$, 
then the condition of Theorem \ref{thm:perverse-heart-lift-criterion}(i) holds for every heart of a bounded $t$-structure on $\mathcal{N}$.
\end{proposition}

\begin{proof}
Let $i\colon\mathcal{N}\hookrightarrow\mathcal{D}$ be the inclusion.
We now prove (i) and the proof of (ii) is dual.

Let $i^!$ be the right adjoint of $i$. For any $E\in\mathcal{D}$, we write $E':=ii^!E$ and complete the counit to a distinguished triangle $E' \longrightarrow E\longrightarrow F$.
For every $K\in\mathcal{N}$, adjunction gives an isomorphism $\Hom_{\mathcal{D}}(K,E')\simeq\Hom_{\mathcal{D}}(K,E)$, and similarly $\Hom_{\mathcal{D}}(K,E'[1])\simeq\Hom_{\mathcal{D}}(K,E[1])$.
Therefore, applying $\Hom_{\mathcal{D}}(K,-)$ to the above triangle gives $\Hom_{\mathcal{D}}(K,F)=0$ and hence $F\in\mathcal{N}^\perp$.

Consider the truncation triangle $\tau_{\mathcal{H_N}}^{\leq 0}E'\longrightarrow E'\longrightarrow\tau_{\mathcal{H_N}}^{\geq 1}E'$ in $\mathcal{N}$. Applying the octahedral axiom to the composition $(\tau_{\mathcal{H_N}}^{\leq 0}E'\longrightarrow E'\longrightarrow E)$, we obtain a triangle $\tau_{\mathcal{H_N}}^{\leq 0}E'\longrightarrow E\longrightarrow E_1$ where $E_1$ fits into a distinguished triangle $\tau_{\mathcal{H_N}}^{\geq 1}E'\longrightarrow E_1\longrightarrow F$.

Now, for $K\in\mathcal{N}_{\mathcal{H_N}}^{\leq 0}$, we have $\Hom_{\mathcal{D}}(K,\tau_{\mathcal{H_N}}^{\geq 1}E')=0$ and $\Hom_{\mathcal{D}}(K,F)=0$.
Applying $\Hom_{\mathcal{D}}(K,-)$ to $\tau_{\mathcal{H_N}}^{\geq 1}E'\longrightarrow E_1\longrightarrow F$ gives $\Hom_{\mathcal{D}}(K,E_1)=0$.
Thus $E_1\in(\mathcal{N}_{\mathcal{H_N}}^{\leq 0})^\perp$, and hence $(\mathcal{N}_{\mathcal{H_N}}^{\leq 0},(\mathcal{N}_{\mathcal{H_N}}^{\leq 0})^\perp)$ is a torsion pair in $\mathcal{D}$.
\end{proof}

\subsection{Lifting of Harder--Narasimhan filtrations}

In this subsection, we prove that Harder--Narasimhan filtrations can be lifted along Serre quotients under some common assumptions.

For the remainder of this section, fix a finite-rank lattice $\Lambda$ and a surjective class map $v\colon\mathrm{K}(\mathcal{D}/\mathcal{N})\longrightarrow\Lambda$.
We use the induced class map $v\circ\mathrm{K}(\pi_{\mathcal{N}})$ on $\mathcal{D}$.

We fix a pre-stability condition $\sigma=(Z_\mathcal{A},\mathcal{A})$ on $\mathcal{D/N}$ and $\mathcal{B} \in \heartsuit^b(\mathcal{D})$ with $\pi_\mathcal{N}(\mathcal{B})\subseteq \mathcal{A}$,
i.e., an exact functor $\pi_\mathcal{N}\colon \mathcal{B}\longrightarrow \mathcal{B}/\mathcal{T}=\mathcal{A}$,
where $\mathcal{T}:=\mathcal{B}\cap\mathcal{N}$.
We write $Z_{\mathcal{B}}:=Z_{\mathcal{A}}\circ \pi_\mathcal{N}$ for the induced weak stability function on $\mathcal{B}$.

The following lemma shows that, under the torsion pair assumption, subobjects downstairs can be lifted inside $\mathcal{T}^{\perp}$.

\begin{lemma}\label{lem:F-lifting}
Let $E \in \mathcal{T}^\perp$ and let $0 \longrightarrow A_1 \longrightarrow \pi_\mathcal{N}(E) \longrightarrow A_2 \longrightarrow 0 $ be a short exact sequence in $\mathcal{A}$.
If $(\mathcal{T},\mathcal{T}^{\perp})$ is a torsion pair in $\mathcal{B}$,
then there exists a short exact sequence
$0\longrightarrow E_1\longrightarrow E\longrightarrow E_2\longrightarrow 0$
in $\mathcal{B}$ such that $E_i\in\mathcal{T}^{\perp}$ and $\pi_\mathcal{N}(E_i)= A_i$.
\end{lemma}

\begin{proof}
By Lemma \ref{lem:Kra22}, there is an exact sequence
$0\longrightarrow E'\longrightarrow E\longrightarrow E''\longrightarrow 0$ in $\mathcal{B}$ with $\pi_\mathcal{N}(E')=A_1$ and $\pi_\mathcal{N}(E'')=A_2$.

Now take the torsion decomposition $0\longrightarrow T\longrightarrow E'' \longrightarrow E_2 \longrightarrow 0$ of $E''$ with $T\in\mathcal{T}$ and $E_2\in\mathcal{T}^{\perp}$.
Let $E_1$ be the kernel of the composition $E\longrightarrow E''\longrightarrow E_2$. Then we have the desired short exact sequence
$0\longrightarrow E_1\longrightarrow E\longrightarrow E_2\longrightarrow 0$.
\end{proof}

\begin{theorem}\label{thm:abstract-HN-lift}
If $(\mathcal{T},\mathcal{T}^{\perp})$ is a torsion pair in $\mathcal{B}$,
then $(Z_{\mathcal{B}}, \mathcal{B})$ is a weak pre-stability condition with $(Z_{\mathcal{B}}, \mathcal{B})^\flat=\sigma$.
\end{theorem}

\begin{proof}
First, $Z_{\mathcal{B}}$ is a weak stability function on $\mathcal{B}$. Indeed, let $0\neq E\in \mathcal{B}$. If $\pi_\mathcal{N}(E)\neq 0$, then $Z_{\mathcal{B}}(E)=Z_\mathcal{A}(\pi_\mathcal{N}(E))\in \mathbb{H}$. If $\pi_\mathcal{N}(E)=0$, then $E\in \mathcal{T}$, so $Z_{\mathcal{B}}(E)=0$.

It remains to prove the Harder--Narasimhan property. 
Now for $0\neq E\in \mathcal{B}$,
we consider the short exact sequence $0\longrightarrow E_{\mathcal{T}}\longrightarrow E\longrightarrow E_{\mathcal{T}^{\perp}}\longrightarrow 0$ with $E_{\mathcal{T}}\in \mathcal{T}$ and $E_{\mathcal{T}^{\perp}}\in \mathcal{T}^{\perp}$.
Note that we may assume $E_{\mathcal{T}^{\perp}}\neq 0$.

Let $0=M_0\subset M_1\subset \cdots \subset M_r=\pi_\mathcal{N}(E_{\mathcal{T}^{\perp}})$ be the Harder--Narasimhan filtration of $\pi_\mathcal{N}(E_{\mathcal{T}^{\perp}})$ with respect to $Z_\mathcal{A}$.
Then each $A_i:=M_i/M_{i-1}$ is $Z_\mathcal{A}$-semistable and $\mu_Z(A_1)>\cdots>\mu_Z(A_r)$.

Set $G_0:=E_{\mathcal{T}^{\perp}}$. For $G_{i-1}\in \mathcal{T}^{\perp}$ with $\pi_\mathcal{N}(G_{i-1})= \pi_\mathcal{N}(E_{\mathcal{T}^{\perp}})/M_{i-1}$,
Lemma \ref{lem:F-lifting} gives an exact sequence $0\longrightarrow F_i\longrightarrow G_{i-1}\longrightarrow G_i\longrightarrow 0$ with $F_i,G_i\in \mathcal{T}^{\perp}$, $\pi_\mathcal{N}(F_i)\cong A_i$, and $\pi_\mathcal{N}(G_i)\cong \pi_\mathcal{N}(E_{\mathcal{T}^{\perp}})/M_i$.
This process terminates at $G_r=0$.

For each $0\leq i\leq r$, let $E_i:=\ker(E\longrightarrow E_{\mathcal{T}^{\perp}}\longrightarrow G_i)$.
Then $E_0=E_{\mathcal{T}}$, $E_r=E$, and $E_i/E_{i-1}\cong F_i$ for all $1\leq i\leq r$.

We claim that if $F\in \mathcal{T}^{\perp}$ and $\pi_\mathcal{N}(F)$ is $Z$-semistable in $\mathcal{A}$, then $F$ is $Z_{\mathcal{B}}$-semistable in $\mathcal{B}$. Indeed, let $0\neq G\subset F$ be a subobject. Since $\mathcal{T}^{\perp}$ is closed under subobjects, we have $G\in \mathcal{T}^{\perp}$,
and hence $\pi_\mathcal{N}(G)$ is a nonzero subobject of $\pi_\mathcal{N}(F)$ in $\mathcal{A}$.
By semistability of $\pi_\mathcal{N}(F)$, we have $\mu_Z(\pi_\mathcal{N}(G))\leq \mu_Z(\pi_\mathcal{N}(F))$, hence $\mu_{Z_{\mathcal{B}}}(G)\leq \mu_{Z_{\mathcal{B}}}(F)$.

By our claim, each $F_i$ is $Z_{\mathcal{B}}$-semistable, and $\mu_{Z_{\mathcal{B}}}(F_i)=\mu_Z(A_i)$. Therefore the slopes of $F_1,\dots,F_r$ are strictly decreasing.

If $E_{\mathcal{T}}=0$, then $0\subset E_1\subset \cdots \subset E_r=E$ is a Harder--Narasimhan filtration of $E$.
If $E_{\mathcal{T}}\neq 0$ and $\mu_{Z_{\mathcal{B}}}(F_1)=+\infty$,
then $E_1$ is an extension of $F_1$ by $E_{\mathcal{T}}$ and hence semistable,
so $0\subset E_1\subset \cdots \subset E_r=E$ is a Harder--Narasimhan filtration of $E$.
If $E_{\mathcal{T}}\neq 0$ and $\mu_{Z_{\mathcal{B}}}(F_1)<+\infty$, then $0\subset E_{\mathcal{T}}\subset E_1\subset \cdots \subset E_r=E$ is a Harder--Narasimhan filtration of $E$.

Thus $(Z_{\mathcal{B}},\mathcal{B})$ is a weak pre-stability condition,
and Proposition \ref{prop:weak-BS-descend}(i) implies that $(Z_{\mathcal{B}},\mathcal{B})^\flat=\sigma$.
\end{proof}

\begin{remark}
Conversely, if $(Z_\mathcal{B},\mathcal{B})$ admits the Harder--Narasimhan property and if $\mathcal{P}_{\mathcal{B}}(1)$ is Noetherian, then $(\mathcal{T},\mathcal{T}^\perp)$ is a torsion pair in $\mathcal{B}$.

Indeed, $\mathcal{P}_{\mathcal{B}}(1)$ is Noetherian,
so $\mathcal{T}$ is a torsion class in $\mathcal{P}_{\mathcal{B}}(1)$.
Moreover, every semistable object of finite slope belongs to $\mathcal{T}^\perp$.
Considering the Harder--Narasimhan filtration of an arbitrary object in $\mathcal{B}$, we conclude that $(\mathcal{T},\mathcal{T}^\perp)$ is a torsion pair in $\mathcal{B}$.
\end{remark}

In summary, we have the following result.

\begin{corollary}\label{cor:weak-BS-lift}
Let $\mathcal{H_N} \in \heartsuit^b(\mathcal{N})$ and assume that $\Stab_{\Lambda}(\mathcal{D/N})\neq\varnothing$. 
Then the equivalent conditions of Theorem \ref{thm:perverse-heart-lift-criterion} hold if and only if
for every $\sigma=(Z_{\mathcal{A}},\mathcal{A})\in\Stab_{\Lambda}(\mathcal{D/N})$, the pair $\sigma^\sharp_{\mathcal{H_N}}:=(Z_{\mathcal{A}}\circ\pi_{\mathcal{N}},\Per(\mathcal{D},\mathcal{N};\mathcal{A},\mathcal{H_N}))$ belongs to $\Stab^w_{\Lambda}(\mathcal{D},\mathcal{N})$.
\end{corollary}

\begin{proof}
Assume that the equivalent conditions of Theorem \ref{thm:perverse-heart-lift-criterion} hold and consider an arbitrary stability condition $\sigma=(Z_{\mathcal{A}},\mathcal{A})\in\Stab_{\Lambda}(\mathcal{D/N})$.
Then $\Per(\mathcal{D},\mathcal{N};\mathcal{A},\mathcal{H_N})$ is a heart of a bounded $t$-structure on $\mathcal{D}$.

By Proposition \ref{prop:core-equivalence}(i), $(\mathcal{H_N},\mathcal{H_N}^{\perp})$ is a torsion pair in this heart. 
Therefore Theorem \ref{thm:abstract-HN-lift} applies and shows that $\sigma^\sharp_{\mathcal{H_N}}$ is a weak pre-stability condition, with $(\sigma^\sharp_{\mathcal{H_N}})^\flat=\sigma$.

It remains to verify the support property.
Let $Q$ be a quadratic form on $\Lambda_\mathbb{R}$ giving the support property for $\sigma$.
If $E$ is $\sigma^\sharp_{\mathcal{H_N}}$-semistable,
then by Proposition \ref{prop:weak-BS-descend}(ii),
$\pi_\mathcal{N}(E)$ is either zero or $\sigma$-semistable.
If $\pi_\mathcal{N}(E)=0$, then $Q(E)=0$. If $\pi_\mathcal{N}(E)\neq 0$, then the support property of $\sigma$ gives $Q(\pi_\mathcal{N}(E))\geq 0$, hence $Q(E)\geq 0$.

Conversely, choose $\sigma=(Z_{\mathcal{A}},\mathcal{A})\in\Stab_\Lambda(\mathcal{D/N})$.
Applying the assumption to its shifts, we see that $\mathcal{B}_n:=\Per(\mathcal{D},\mathcal{N};\mathcal{A}[n],\mathcal{H_N})$ is a heart of a bounded $t$-structure for every $n\in\mathbb{Z}$,
that is, condition (iii) of Theorem \ref{thm:perverse-heart-lift-criterion} holds.
\end{proof}

For any heart $\mathcal{H_N}$ of a bounded $t$-structure on $\mathcal{N}$ satisfying the torsion-pair condition of Theorem \ref{thm:perverse-heart-lift-criterion}(i),
this gives \gls{sharp} with $\flat\circ\sharp_{\mathcal{H_N}}=\mathrm{id}$.

\subsection{Topology of descent and lifting}

We now study the continuity of $\flat$ and $\sharp_{\mathcal{H_N}}$.

\begin{proposition}\label{prop:descent-continuity}
Let $\sigma_i=(Z_i,\mathcal{P}_i)\in\Stab^\mathrm{w}_\Lambda(\mathcal{D},\mathcal{N})$ and write $\sigma_i^\flat=(Z_i,\mathcal{Q}_i)$ for $i=1,2$. Then we have
$d^{\mathcal{D/N}}_{\mathrm{sl}}(\mathcal{Q}_1,\mathcal{Q}_2)
\leq
d^{\mathcal{D}}_{\mathrm{sl}}(\mathcal{P}_1,\mathcal{P}_2)$
and
$d^{\mathcal{D/N}}(\sigma_1^\flat,\sigma_2^\flat)
\leq
d^{\mathcal{D}}(\sigma_1,\sigma_2)$.
In particular, the descent map $\flat$ is continuous.
\end{proposition}

\begin{proof}
Assume that $d^{\mathcal{D}}_{\mathrm{sl}}(\mathcal{P}_1,\mathcal{P}_2)<\infty$ and fix $\delta>d^{\mathcal{D}}_{\mathrm{sl}}(\mathcal{P}_1,\mathcal{P}_2)$.

Let $0\neq A\in\mathcal{Q}_2(\phi)$. After shifting into the heart, Proposition \ref{prop:weak-BS-descend}(iii) gives a semistable lift $E\in\mathcal{P}_2(\phi)$ with $\pi_\mathcal{N}(E)\simeq A$. By the definition of $d^{\mathcal{D}}_{\mathrm{sl}}$, every Harder--Narasimhan factor of $E$ with respect to $\sigma_1$ has phase in $(\phi-\delta,\phi+\delta)$.

Applying $\pi_\mathcal{N}$ to the Harder--Narasimhan filtration of $E$ with respect to $\sigma_1$, Proposition \ref{prop:weak-BS-descend}(ii) shows that the nonzero images form the Harder--Narasimhan filtration of $A$ with respect to $\sigma_1^\flat$, with the same phases,
and hence $A\in\mathcal{Q}_1([\phi-\delta,\phi+\delta])$.

By the characterization \eqref{eq:slicing-metric} of $d_{\mathrm{sl}}$, we obtain
$d^{\mathcal{D/N}}_{\mathrm{sl}}(\mathcal{Q}_1,\mathcal{Q}_2)\leq\delta$.
This proves the first inequality, and the second inequality follows from Proposition \ref{prop:weak-BS-descend}(i).
\end{proof}

The slicing-distance estimate in the proof above is a special case of \cite[Corollary 3.8]{BPPW22}.

Let $\mathcal{H_N}\in\heartsuit^{b}(\mathcal{N})$. We define $\gls{null-heart-stratum}:=\{(Z,\mathcal{B})\in\Stab^\mathrm{w}_\Lambda(\mathcal{D},\mathcal{N})\mid \mathcal{B}\cap\mathcal{N}=\mathcal{H_N}\}$
and $\gls{lifting-image}:=\flat(\mathcal{S}_{\mathcal{H_N}})\subseteq\Stab_\Lambda(\mathcal{D/N})$.
Note that $\sigma \in \mathcal{I}_{\mathcal{H_N}}$ exactly when $\sigma^\sharp_{\mathcal{H_N}}$ is a weak stability condition.

\begin{proposition}\label{prop:fixed-null-heart-stratum}
For every $\mathcal{H_N}\in\heartsuit^{b}(\mathcal{N})$, the descent map restricts to an isometry
$\flat\colon\mathcal{S}_{\mathcal{H_N}}\longrightarrow\mathcal{I}_{\mathcal{H_N}}$
with respect to both $d_{\mathrm{sl}}$ and $d$.
\end{proposition}

\begin{proof}
It is clear that $\flat|_{\mathcal{S}_{\mathcal{H_N}}}$ is bijective.
Surjectivity follows from the definition of $\mathcal{I}_{\mathcal{H_N}}$,
and injectivity follows from Proposition \ref{prop:uniqueness-of-descended-lifts}.

Let $\sigma_i=(Z_i,\mathcal{P}_i)\in\mathcal{S}_{\mathcal{H_N}}$ and write $\sigma_i^\flat=(Z_i,\mathcal{Q}_i)$ for $i=1,2$.
By Proposition \ref{prop:descent-continuity}, it remains to prove $d^{\mathcal{D/N}}_{\mathrm{sl}}(\mathcal{Q}_1,\mathcal{Q}_2)
\geq
d^{\mathcal{D}}_{\mathrm{sl}}(\mathcal{P}_1,\mathcal{P}_2)$. 
Assume that $d^{\mathcal{D/N}}_{\mathrm{sl}}(\mathcal{Q}_1,\mathcal{Q}_2)<\infty$ and let $\delta>d^{\mathcal{D/N}}_{\mathrm{sl}}(\mathcal{Q}_1,\mathcal{Q}_2)$. 

Let $0\neq E\in\mathcal{P}_2(\phi)$. 
If $E\in\mathcal{N}$, then $\phi\in\mathbb{Z}$ and $E[1-\phi] \in \mathcal{B}_2 \cap \mathcal{N}= \mathcal{B}_1 \cap \mathcal{N}$ so $E \in \mathcal{P}_1(\phi)$.

If $E\notin\mathcal{N}$, set $A:=\pi_\mathcal{N}(E)$. By Proposition \ref{prop:weak-BS-descend}(ii), $A\in\mathcal{Q}_2(\phi)$.
As in the proof of Proposition \ref{prop:descent-continuity}, the massive Harder--Narasimhan factors of $E$ with respect to $\sigma_1$ descend to the Harder--Narasimhan factors of $A$ with respect to $\sigma_1^\flat$ with the same phases,
and hence every massive factor of $E$ has $\sigma_1$-phase in $(\phi-\delta,\phi+\delta)$.

It remains to control the massless factors.
If there is a massless factor of phase greater than $\phi+\delta$, say the highest Harder--Narasimhan factor $M$ is massless of some phase $q>\phi$,
then the first $\sigma_1$-Harder--Narasimhan truncation gives a distinguished triangle $M\longrightarrow E\longrightarrow E'$.
Every $\sigma_1$-Harder--Narasimhan factor of $E'[-1]$ has phase smaller than $q$, so
$\Hom_{\mathcal{D}}(M,E'[-1])=0$,
and hence the morphism $M\longrightarrow E$ is nonzero.

On the other hand, since $M \in \mathcal{N}$, we similarly have  $M \in \mathcal{P}_2(q)$.
As $E \in \mathcal{P}_2(\phi)$ and $\phi<q$, we have
$\Hom_{\mathcal{D}}(M,E)=0$, a contradiction.

The case of a massless Harder--Narasimhan factor of $\sigma_1$-phase less than $\phi-\delta$ is dual.
Therefore every Harder--Narasimhan factor of $E$ with respect to $\sigma_1$ has $\sigma_1$-phase in $(\phi-\delta,\phi+\delta)$. 
This proves the assertion.
\end{proof}

The proposition together with Corollary \ref{cor:weak-BS-lift} immediately gives the following.

\begin{corollary}\label{cor:canonical-lifting-homeomorphism}
Assume that $\Stab_\Lambda(\mathcal{D/N})\neq\varnothing$.
Then $\mathcal{I}_{\mathcal{H_N}}=\Stab_\Lambda(\mathcal{D/N})$ if and only if $(\mathcal{N}^{\leq 0},(\mathcal{N}^{\leq 0})^\perp)$ and $({}^\perp(\mathcal{N}^{\geq 0}),\mathcal{N}^{\geq 0})$ are torsion pairs in $\mathcal{D}$.
If this is the case, the lifting map
$\sharp_{\mathcal{H_N}}$
is an isometry with inverse $\flat$.
\end{corollary}

\begin{proposition}\label{prop:null-heart-separation}
Let $\mathcal{H_N}\neq\mathcal{H'_N}$ be hearts of bounded $t$-structures on $\mathcal{N}$.
For any $\sigma=(Z,\mathcal{P})\in\mathcal{S}_{\mathcal{H_N}}$ and $\tau=(W,\mathcal{Q})\in\mathcal{S}_{\mathcal{H'_N}}$,
we have $d^{\mathcal{D}}_{\mathrm{sl}}(\mathcal{P},\mathcal{Q})\geq 1$.
\end{proposition}

\begin{proof}
For $E\in\mathcal{H_N}\setminus\mathcal{H'_N}$,
$E$ is $\sigma$-semistable of phase $1$,
so
$\phi^+_{\mathcal{P}}(E)=\phi^-_{\mathcal{P}}(E)=1$.

Since $E\in\mathcal{N}$ is massless, all its Harder--Narasimhan factors with respect to $\tau$ have integral phases.
If
$\phi^+_{\mathcal{Q}}(E)=\phi^-_{\mathcal{Q}}(E)=1$,
then all these factors belong to
$\mathcal{B}_\tau\cap\mathcal{N}=\mathcal{H'_N}$.
Since $\mathcal{H'_N}$ is extension-closed, this would imply $E\in\mathcal{H'_N}$, a contradiction. Thus at least one of
$\phi^+_{\mathcal{Q}}(E)$ and $\phi^-_{\mathcal{Q}}(E)$ differs from $1$ by at least $1$, and hence
$d^{\mathcal{D}}_{\mathrm{sl}}(\mathcal{P},\mathcal{Q})\geq 1$.
\end{proof}

In particular, every nonempty $\mathcal{S}_{\mathcal{H_N}}$ is open and closed in $\Stab^\mathrm{w}_\Lambda(\mathcal{D},\mathcal{N})$, and every connected subset is contained in a unique stratum.

\begin{corollary}\label{cor:product}
There is a natural homeomorphism
\begin{align*}
\Stab^{\mathrm{w}}_\Lambda(\mathcal{D},\mathcal{N})
\simeq
\coprod_{\mathcal{H_N}\in\heartsuit^{b}(\mathcal{N})}I_{\mathcal{H_N}}.
\end{align*}
Assume moreover that $\Stab_\Lambda(\mathcal{D/N})\neq\varnothing$.
Then the torsion-pair condition of Theorem \ref{thm:perverse-heart-lift-criterion}(i) holds for every $\mathcal{H_N}\in\heartsuit^{b}(\mathcal{N})$
if and only if $\tau=(Z,\mathcal{B})\longmapsto(\mathcal{B}\cap\mathcal{N},\tau^\flat)$ gives a natural homeomorphism
\begin{align*}
\Stab^{\mathrm{w}}_\Lambda(\mathcal{D},\mathcal{N})
\simeq
\heartsuit^{b}(\mathcal{N})\times\Stab_\Lambda(\mathcal{D/N}),
\end{align*}
where $\heartsuit^{b}(\mathcal{N})$ is equipped with the discrete topology.
\end{corollary}

We conclude this section with a compatibility result with autoequivalences.

\begin{proposition}\label{lem:transport}
If $\Phi \in \mathrm{Aut}(\mathcal{D})$ satisfies $\pi_{\mathcal{N}}\circ\Phi\simeq\pi_{\mathcal{N}}$, then we have:
\begin{enumerate}
\item $\Phi(\mathcal{N})=\mathcal{N}$ and $\Phi(\Per(\mathcal{D},\mathcal{N};\mathcal{A},\mathcal{H_N}))=\Per(\mathcal{D},\mathcal{N};\mathcal{A},\Phi(\mathcal{H_N}))$ for any heart $\mathcal{A}\in\heartsuit^{b}(\mathcal{D}/\mathcal{N})$.
Moreover, the torsion-pair condition of Theorem \ref{thm:perverse-heart-lift-criterion}(i) holds for $\mathcal{H_N}$ if and only if it holds for $\Phi(\mathcal{H_N})$.

\item The functor $\Phi$ acts on $\Stab^{\mathrm{w}}_\Lambda(\mathcal{D},\mathcal{N})$ by $\Phi\cdot(Z,\mathcal{B})=(Z\circ\Phi^{-1},\Phi(\mathcal{B}))$.
We have $(\Phi\cdot\tau)^\flat=\tau^\flat$ for every $\tau\in\Stab^{\mathrm{w}}_\Lambda(\mathcal{D},\mathcal{N})$ and $\Phi\cdot\sigma^\sharp_{\mathcal{H_N}}=\sigma^\sharp_{\Phi(\mathcal{H_N})}$ for every $\sigma\in I_{\mathcal{H_N}}$.

\item The action of $\Phi$ is an isometry, that is, $d^{\mathcal{D}}(\Phi\cdot\tau_1,\Phi\cdot\tau_2)=d^{\mathcal{D}}(\tau_1,\tau_2)$ for all $\tau_1,\tau_2\in\Stab^{\mathrm{w}}_\Lambda(\mathcal{D},\mathcal{N})$.
\end{enumerate}
\end{proposition}

\begin{proof}
For (i), the first assertion is trivial.
Note that for every $q\in\mathbb{Z}$ we have
$\Phi\bigl({}^\perp(\mathcal{N}_{\mathcal{H_N}}^{\geq q})\bigr)
={}^\perp(\mathcal{N}_{\Phi(\mathcal{H_N})}^{\geq q})$ and $
\Phi\bigl((\mathcal{N}_{\mathcal{H_N}}^{\leq q})^\perp\bigr)
=(\mathcal{N}_{\Phi(\mathcal{H_N})}^{\leq q})^\perp.$
Moreover, $\pi_{\mathcal{N}}\circ\Phi\simeq\pi_{\mathcal{N}}$ gives $\Phi(\pi_{\mathcal{N}}^{-1}(\mathcal{A}))=\pi_{\mathcal{N}}^{-1}(\mathcal{A})$, so
\begin{align*}
\Phi\bigl[\Per(\mathcal{D},\mathcal{N};\mathcal{A},\mathcal{H_N})\bigr]
&=\Phi\bigl[
\pi_{\mathcal{N}}^{-1}(\mathcal{A})
\cap{}^\perp(\mathcal{N}_{\mathcal{H_N}}^{\geq 1})
\cap(\mathcal{N}_{\mathcal{H_N}}^{\leq -1})^\perp
\bigr]\\
&=\pi_{\mathcal{N}}^{-1}(\mathcal{A})
\cap{}^\perp(\mathcal{N}_{\Phi(\mathcal{H_N})}^{\geq 1})
\cap(\mathcal{N}_{\Phi(\mathcal{H_N})}^{\leq -1})^\perp\\
&=\Per(\mathcal{D},\mathcal{N};\mathcal{A},\Phi(\mathcal{H_N})).
\end{align*}
Applying $\Phi$ and $\Phi^{-1}$ to the torsion decomposition triangles respectively proves the equivalence of the torsion-pair conditions.

For (ii), let $\tau=(Z,\mathcal{B})\in\Stab^{\mathrm{w}}_\Lambda(\mathcal{D},\mathcal{N})$ and set $\Phi\cdot\tau:=(Z\circ\Phi^{-1},\Phi(\mathcal{B}))$.
It is clear that $\Phi\cdot\tau$ is a weak stability condition on $\mathcal{D}$, with massless subcategory $\mathcal{N}_{\Phi\cdot\tau}=\Phi(\mathcal{N}_\tau)=\mathcal{N}$. 
Proposition \ref{prop:weak-BS-descend} then gives
\[(\Phi\cdot\tau)^\flat
=(Z,\pi_{\mathcal{N}}(\Phi(\mathcal{B})))
=(Z,\pi_{\mathcal{N}}(\mathcal{B}))
=\tau^\flat.\]

For the second assertion, given $\sigma=(W,\mathcal{A})\in I_{\mathcal{H_N}}$,
the weak stability condition $\Phi\cdot\sigma^\sharp_{\mathcal{H_N}}$ has heart $\Per(\mathcal{D},\mathcal{N};\mathcal{A},\Phi(\mathcal{H_N}))$, and its central charge is $W\circ\pi_{\mathcal{N}}\circ\Phi^{-1}=W\circ\pi_{\mathcal{N}}$,
so $\Phi\cdot\sigma^\sharp_{\mathcal{H_N}}=\sigma^\sharp_{\Phi(\mathcal{H_N})}$.

For (iii), $\Phi$ takes $\tau$-semistable objects to $\Phi\cdot\tau$-semistable objects of the same phase.
Applying $\Phi$ to a Harder--Narasimhan filtration of $E$ then gives a Harder--Narasimhan filtration of $\Phi(E)$ whose factors have the same phases,
so the action of $\Phi$ preserves $d^{\mathcal{D}}_{\mathrm{sl}}$.
Finally, $\pi_{\mathcal{N}}\circ\Phi^{-1}\simeq\pi_{\mathcal{N}}$ gives $v\circ\mathrm{K}(\pi_{\mathcal{N}})\circ\mathrm{K}(\Phi^{-1})=v\circ\mathrm{K}(\pi_{\mathcal{N}})$,
and hence $\Phi$ fixes each central charge as a map on $\Lambda$. 
This proves the assertion.
\end{proof}

\section{Geometric applications}
In this section, we give geometric applications of the results of Section 3.
In each case, we consider a Verdier localization induced by derived pushforward, and we denote the kernel by $\mathcal{N}$.

We verify the torsion-pair condition of Theorem \ref{thm:perverse-heart-lift-criterion}(i) for a chosen heart of a bounded $t$-structure on $\mathcal{N}$,
and then Corollary \ref{cor:weak-BS-lift} lifts Bridgeland stability conditions.
When moreover the condition holds for every heart of a bounded $t$-structure on $\mathcal{N}$,
Corollary \ref{cor:product} gives the corresponding product decomposition of the space of weak stability conditions.
Throughout this section, we work over an algebraically closed field $k$.

\subsection{Relative dimension one contractions}

We first study projective morphisms of relative dimension one.
We recall the Verdier localization theorem by Bondal--Kapranov--Schechtman,
and then we prove that the restriction of coherent sheaves satisfies the torsion-pair condition of Theorem \ref{thm:perverse-heart-lift-criterion}(i) and characterize when it is equivalent to modules over a finite-dimensional algebra.

Let $f\colon Y\longrightarrow X$ be a projective morphism between quasi-projective $k$-schemes of relative dimension $1$ with $\mathbf{R}f_\ast\mathcal{O}_Y\simeq\mathcal{O}_X$.

We set $\mathcal{N}:=\ker\mathbf{R}f_\ast\subset\mathrm{D^b}(Y)$ and $\gls{Af}:=\Coh(Y)\cap\mathcal{N}$.
The spectral sequence argument in \cite[Lemma 3.1]{Bri02} shows that $\mathcal{A}_f$ is a heart of a bounded $t$-structure on $\mathcal{N}$.

\begin{proposition}\cite[Theorem 2.14]{BKS18}\label{prop:BKS}
$\mathbf{R}f_\ast$ induces a Verdier localization
\[
\mathrm{D^b}(Y)/\ker \mathbf{R}f_\ast \xlongrightarrow{\sim} \mathrm{D^b}(X).
\]
\end{proposition}

\begin{theorem}\label{thm:geometric-null-heart}
The pairs $(\mathcal{N}^{\leq 0},(\mathcal{N}^{\leq 0})^\perp)$ and $({}^\perp(\mathcal{N}^{\geq 0}),\mathcal{N}^{\geq 0})$ are torsion pairs in $\mathrm{D^b}(Y)$.
\end{theorem}

\begin{proof}
We use the adjunction $\mathbf{L}f^\ast\dashv\mathbf{R}f_\ast\dashv f^!$ between $\mathrm{D_{Qcoh}}(Y)$ and $\mathrm{D_{Qcoh}}(X)$.

For $E\in\mathrm{D^b}(Y)$, consider the distinguished triangles $E'\longrightarrow E\longrightarrow f^!\mathbf{R}f_\ast E$ and $\mathbf{L}f^\ast\mathbf{R}f_\ast E\longrightarrow E\longrightarrow E''$, induced respectively by the unit and counit of the two adjunctions.
Note that by the projection formula we have $\mathbf{R}f_\ast\mathbf{L}f^\ast=1$ and taking right adjoints gives $\mathbf{R}f_\ast f^!=1$,
and hence $\mathbf{R}f_\ast E'=0=\mathbf{R}f_\ast E''$.

Note also that $\mathbf{L}f^\ast\mathbf{R}f_\ast E \in \mathrm{D}^-(Y)$ and $f^!\mathbf{R}f_\ast E \in \mathrm{D}^+(Y)$.
It follows that $E'\in \mathrm{D}^+(Y)$ and $E''\in  \mathrm{D}^-(Y)$,
and therefore $\tau^{\leq 0}E'\in\mathcal{N}^{\leq 0}$ and $\tau^{\geq 0}E''\in\mathcal{N}^{\geq 0}$.

Applying the octahedral axiom to the composition
$(\tau^{\leq 0}E'\longrightarrow E'\longrightarrow E)$,
we obtain a distinguished triangle
$\tau^{\leq 0}E'\longrightarrow E\longrightarrow F$
where $F$ fits into a distinguished triangle
$\tau^{\geq 1}E'\longrightarrow F\longrightarrow f^!\mathbf{R}f_\ast E$.

For $K\in\mathcal{N}^{\leq 0}$, we have
$\Hom(K,\tau^{\geq 1}E')=0$ by the standard t-structure, while adjunction gives
$\Hom(K,f^!\mathbf{R}f_\ast E)=\Hom(\mathbf{R}f_\ast K,\mathbf{R}f_\ast E)=0$.
Therefore $\Hom(K,F)=0$, and hence $F\in(\mathcal{N}^{\leq 0})^\perp$.
This gives the first torsion pair.

The second torsion pair is obtained similarly by applying the octahedral axiom to the composition
$(E\longrightarrow E''\longrightarrow\tau^{\geq 0}E'')$.
\end{proof}

\begin{proposition}\label{prop:VdB-null-heart-module}
Assume moreover that $f$ is birational. 
Then there exists a finite-dimensional $k$-algebra $\Gamma$ such that $\mathcal{A}_f\simeq\mathrm{mod}\text{-}\Gamma$ if and only if $f(\mathrm{Exc}(f))$ is a finite set of closed points.
\end{proposition}

\begin{proof}
Assume that $f(\mathrm{Exc}(f))=\{x_0,\dots,x_r\}$.
As $\mathcal{A}_f$ factors as the direct sum $\bigoplus_i(\mathcal{A}_f)_{x_i}$, it suffices to treat the case $f(\mathrm{Exc}(f))=\{x_0\}$.
Moreover, by restricting to an affine open neighbourhood $U \subseteq X$ of $x_0$,
we may assume that $X=\mathrm{Spec}\, R$.

By \cite[Proposition 3.2.5]{VdB04}, there exists a vector bundle $P=P_0\oplus \mathcal{O}_Y$ which is a projective generator of ${}^{-1}\Per(Y/X)$. Setting $A:=\End_Y(P)$, \cite[Corollary 3.2.8]{VdB04} shows that the functor $\Hom_Y(P,-)$ induces an exact equivalence ${}^{-1}\Per(Y/X)\xrightarrow{\sim}\mathrm{mod}\text{-}A$.

Let $i\colon \mathcal{O}_Y \longrightarrow P$ and $p\colon P \longrightarrow \mathcal{O}_Y$ be the inclusion and projection for the decomposition $P=P_0\oplus \mathcal{O}_Y$, and set $e$ to be the idempotent $i\circ p\in A$. 
This implies that
$\Hom_Y(P,E)e = \Hom_Y(\mathcal{O}_Y,E)= \Hom_X(\mathcal{O}_X,\mathbf{R}f_\ast E)$
for any object $E\in {}^{-1}\Per(Y/X)$.
Since $\mathbf{R}f_\ast E\in \Coh(X)$ and $X$ is affine, this vanishes if and only if $\mathbf{R}f_\ast E=0$.
It follows that
$\mathcal{A}_f[1]\simeq \{M\in \mathrm{mod}(A)\mid Me=0\}=\mathrm{mod}\text{-}A/AeA$.

It remains to show that $B:=A/AeA$ is a finite-dimensional $k$-algebra.
Let $F\in \mathcal{A}_f$ be the sheaf with $\Hom_Y(P,F[1])=B$.

We then have
$B =\Ext^1_Y(P,F)= H^1(Y,P^\vee\otimes F)=H^0(X,\mathbf{R}^1f_\ast(P^\vee\otimes F))$ since $P$ is a vector bundle and $X$ is affine.
As $\mathbf{R}f_\ast F=0$, we have $\supp F \subseteq f^{-1}(x_0)$ and hence $\supp \mathbf{R}f_\ast(P^\vee\otimes F)\subseteq\{x_0\}$.
Therefore $B$, viewed as an $R$-module, is finitely generated and supported at  $x_0$,
that is, $\mathfrak{m}_{x_0}^n B=0$ for some $n$ and hence $B$ is a finitely generated $R/\mathfrak{m}_{x_0}^n$-module.
Since $R/\mathfrak{m}_{x_0}=k$,
we see that $R/\mathfrak{m}_{x_0}^n$ is finite-dimensional over $k$, and hence so is $B$.

Conversely, assume that $\mathcal{A}_f\simeq\mathrm{mod}\text{-}\Gamma$ for some finite-dimensional $k$-algebra $\Gamma$.
If $f^{-1}(x)$ is zero-dimensional, then since $f$ is projective and $\mathbf{R}f_\ast\mathcal{O}_Y\simeq\mathcal{O}_X$, $f$ is an isomorphism near $x$,
so every closed point $x\in f(\mathrm{Exc}(f))$ has a one-dimensional fiber.

Let $x$ be a closed point of $f(\mathrm{Exc}(f))$, and choose a reduced irreducible one-dimensional component $C$ of $f^{-1}(x)$.
Then applying $\mathbf{R}f_\ast$ to the exact sequence $0\longrightarrow\mathcal{I}_C\longrightarrow\mathcal{O}_Y\longrightarrow\mathcal{O}_C\longrightarrow0$ gives $\mathbf{R}^1f_\ast\mathcal{O}_C=0$,
and hence $H^1(C,\mathcal{O}_C)=0$.

Choose a smooth closed point $p\in C$, and let $\mathcal{I}_p\subset\mathcal{O}_C$ be its ideal sheaf.
The natural morphism $H^0(C,\mathcal{O}_C)\longrightarrow H^0(C,\mathcal{O}_p)$ is an isomorphism,
and the exact sequence
$0\longrightarrow\mathcal{I}_p\longrightarrow\mathcal{O}_C\longrightarrow\mathcal{O}_p\longrightarrow0$ then gives $H^0(C,\mathcal{I}_p)=H^1(C,\mathcal{I}_p)=0$,
so $\mathbf{R}f_\ast\mathcal{I}_p=0$,
that is, $\mathcal{I}_p\in\mathcal{A}_f$.

We claim that $\mathcal{I}_p$ is simple in $\mathcal{A}_f$.
Given a short exact sequence in $\mathcal{A}_f$ $0\longrightarrow K\longrightarrow \mathcal{I}_p\longrightarrow E \longrightarrow 0$  with $K\neq 0$,
since the inclusion $\mathcal{A}_f\subset\Coh(Y)$ is exact, this is also an exact sequence in $\Coh(Y)$.

As $p$ is nonsingular, $\mathcal{I}_p$ is a line bundle on $C$, and therefore on the curve $C$, $K$ has rank one and $E$ has rank zero. 
Now if $E\neq 0$, then since $E$ is a nonzero zero-dimensional sheaf on $C$, we have $H^0(C,E) \neq 0$,
and hence $f_\ast E \neq 0$, contradicting $E \in \mathcal{A}_f$.
This proves that $\mathcal{I}_p$ is simple.

For distinct closed points $x\in f(\mathrm{Exc}(f))$, the simple objects constructed above have disjoint supports and are therefore mutually non-isomorphic.
Since $\mathcal{A}_f$ has only finitely many simple objects, $f(\mathrm{Exc}(f))$ must have only finitely many closed points.
\end{proof}

For the construction of a projective generator of the coherent null category over an affine base, see \cite[Proposition 2.4]{BB22}.

Note that $\mathrm{D^b}(Y)$ is not assumed to be proper in this subsection, so the condition $\mathcal{A}_f\simeq\mathrm{mod}\text{-}\Gamma$ alone does not provide the ambient torsion-pair condition in this generality.

When $f$ is birational and $f(\mathrm{Exc}(f))$ is finite, Proposition \ref{prop:VdB-null-heart-module} shows that $\mathcal{A}_f$ is of finite length with finitely many simple objects, while Theorem \ref{thm:geometric-null-heart} gives the torsion-pair condition for $\mathcal{A}_f$.
Lemma \ref{lem:finite-amplitude-transfer} then gives the same torsion-pair condition for every heart of a bounded $t$-structure on $\ker\mathbf{R}f_\ast$.
This observation applies in particular to the partial crepant resolutions considered in Section 5.1 below, where $Y$ is singular and $\mathrm{D^b}(Y)$ is not proper.

The following example shows that the geometric criterion of Theorem \ref{thm:geometric-null-heart} applies beyond the scope of the other sufficient criteria.

\begin{example}
Let $\pi\colon\mathbb{F}_2\longrightarrow\mathbb{P}(1,1,2)$ be the minimal resolution  with exceptional curve $C$ and let $E$ be an elliptic curve over $\mathbb{C}$.

Consider the morphism $\pi\times 1_E\colon Y:=\mathbb{F}_2\times E\longrightarrow X:=\mathbb{P}(1,1,2)\times E$ of relative dimension $1$ with $\mathbf{R}(\pi\times 1_E)_\ast\mathcal{O}_Y\simeq\mathcal{O}_X$. Set $\mathcal{N}:=\ker\mathbf{R}(\pi\times 1_E)_\ast$.
The results above give $\mathrm{D^b}(Y)/\mathcal{N}\simeq\mathrm{D^b}(X)$ and the torsion-pair condition of Theorem \ref{thm:perverse-heart-lift-criterion}(i).

We now show that neither of the other two sufficient criteria from Section 3 applies to this example.
First, we claim that $\mathcal{N}$ admits no heart of a bounded $t$-structure which is of finite length with finitely many simples.

Assuming the contrary,
we then have $\mathrm{K}(\mathcal{N})\simeq\mathrm{K}(\mathcal{H_N})$,
which is a free abelian group of finite rank.
On the other hand, consider the exact functor $\mathrm{D^b}(E)\longrightarrow\mathcal{N}$ given by $F\longmapsto\mathcal{O}_C(-1)\boxtimes F$.
This functor has an exact left inverse given by $G\longmapsto\mathbf{R}p_\ast(G\otimes q^\ast\mathcal{O}_{\mathbb{P}^1}(1))$, 
where $p\colon Y\longrightarrow E$ is the second projection, and $q\colon Y\longrightarrow\mathbb{P}^1$ is the composition of the first projection with the ruling of $\mathbb{F}_2$.
Consequently, $\mathrm{K}(E)$ is a direct summand of $\mathrm{K}(\mathcal{N})$, but $\mathrm{K}(E)$ is not finitely generated, a contradiction.

Finally, we claim that $\mathcal{N}$ is neither left nor right admissible.
Assume the contrary.
Then $\mathrm{D^b}(Y)/\mathcal{N}$ would be equivalent to  ${}^\perp\mathcal{N}$ or $\mathcal{N}^\perp$, respectively, and hence would be proper since $\mathrm{D^b}(Y)$ is proper,
but $\mathrm{D^b}(Y)/\mathcal{N}\simeq\mathrm{D^b}(X)$ is not proper, a contradiction.
\end{example}

\subsection{Resolutions of nodal varieties}\label{subsec:nodal-resolutions}

We next consider resolutions of varieties with ordinary double points, the simplest singularities.

In this subsection, assume that $\rm{char}$$(k)\neq 2$ and let $X$ be a projective variety over $k$ of dimension $n\geq 2$ with ordinary double points $x_1,\ldots,x_r$ and no other singularities. 
Let $\pi\colon\widetilde{X}\longrightarrow X$ be the blowup of these points, and set $\mathcal{N}:=\ker\mathbf{R}\pi_\ast\subset\mathrm{D^b}(\widetilde{X})$.
Then \cite[Corollary 5.6, Lemma 5.7]{KS24} shows that $\mathbf{R}\pi_\ast$ induces a Verdier localization $\mathrm{D^b}(\widetilde{X})/\mathcal{N}\simeq\mathrm{D^b}(X)$.

The following proposition shows that in this setting, even the standard $t$-structure need not restrict to $\mathcal{N}$.

\begin{proposition}
If $n\geq3$, then $\Coh(\widetilde{X})\cap\mathcal{N}$ is not a heart of a bounded $t$-structure on $\mathcal{N}$.
\end{proposition}
\begin{proof}
Let $i\colon E\hookrightarrow\widetilde{X}$ be the inclusion of the exceptional divisor of $\pi$.
The divisor $E$ is a disjoint union of $r$ smooth quadrics of dimension $n-1$, so $\mathbf{R}\Gamma(E,\mathcal{O}_E)\simeq k^r$.

There is a semiorthogonal decomposition $\mathrm{D^b}(E)=\langle\ker\mathbf{R}\Gamma(E,-),\mathcal{O}_E\rangle$,
and the projection $\pr$ onto $\langle\mathcal{O}_E\rangle$ is given by $\mathbf{R}\Gamma(E,-)\otimes_{k^r}\mathcal{O}_E$.
Serre duality gives $\pr(\omega_E)=\mathbf{R}\Gamma(E,\omega_E)\otimes_{k^r}\mathcal{O}_E\simeq\mathcal{O}_E[-(n-1)]$. 
Since $\omega_E\in\mathrm{D^b}(E)^{\leq 0}$, this shows that $\pr(\mathrm{D^b}(E)^{\leq 0}) \nsubseteq\mathrm{D^b}(E)^{\leq1}$.
By \cite[Lemma 4.7]{KLP26}, the standard $t$-structure on $\mathrm{D^b}(E)$ does not restrict to $\ker\mathbf{R}\Gamma(E,-)$.

We may therefore choose $F\in\ker\mathbf{R}\Gamma(E,-)$ such that $\mathbf{R}\Gamma(E,H^0(F))\neq0$.        
That is, we have $i_\ast F\in\mathcal{N}$ and $H^0(i_\ast F)=i_\ast H^0(F)\notin\mathcal{N}$,
which proves the assertion.
\end{proof}

We therefore need to construct a different heart on $\mathcal{N}$. 
To do so, we first recall the following description of the kernel by Kuznetsov and Shinder.

\begin{proposition}\cite[Theorem 5.8]{KS24}\label{prop:KS24}
There exist $E_1,\dots,E_m,K_1,\dots,K_r\in\mathcal{N}$ with the following properties:
\begin{enumerate}
\item $\mathcal{N}=\langle E_1,\dots,E_m,K_1,\dots,K_r\rangle$, and $(E_1,\dots,E_m)$ is an exceptional collection, possibly empty.
\item $\Hom(K_i,E_j[p])=0$ for all $i,j$ and $p\in\mathbb{Z}$.
\item The objects $K_1,\ldots,K_r$ are completely orthogonal and for every $i$ we have
\begin{equation*}
\Hom(K_i,K_i[p])\simeq
\begin{cases}
k & \text{if } n \text{ is even and } p\in\{0,2\},\\
k & \text{if } n \text{ is odd and } p\in\{0,3\},\\
0 & \text{otherwise}.
\end{cases}
\end{equation*}
\end{enumerate}
\end{proposition}

We can then construct the desired heart on $\mathcal{N}$.

\begin{theorem}\label{thm:module-for-nodal}
The category $\mathcal{N}$ admits a heart $\mathcal{H_N}$ of a bounded $t$-structure with $\mathcal{H_N}\simeq\mathrm{mod}\text{-}k^{r+m}$.
\end{theorem}
\begin{proof}
Let $E_1,\dots,E_m,K_1,\dots,K_r$ be as in Proposition \ref{prop:KS24}.
We will choose integers $p_1,\dots,p_m$ so that $E_1[p_1],\dots,E_m[p_m],K_1,\dots,K_r$ form a simple-minded collection in the sense of \cite[Definition 3.1]{CSPP22}.
If $m=0$, there is nothing to do.

If $m>0$, by properness we may choose $q>0$ sufficiently large so that $\Hom(E_i,E_j[s])=0=\Hom(E_i,K_j[s])$ for all $i,j$ and $s\leq 1-q$.
Now we set $p_i:=(m-i+1)q$. For $p\leq1$, we then have 
\[\Hom(E_i[p_i],E_j[p_j+p])=\Hom(E_i,E_j[p-(j-i)q])=0 \text{ whenever } i<j, \text{ and} \]
\[\Hom(E_i[p_i],K_j[p])=\Hom(E_i,K_j[p-p_i])=0 \text{ for all }i,j.\]
Combining this with Proposition \ref{prop:KS24}, we see that $E_1[p_1],\dots,E_m[p_m],K_1,\dots,K_r$ form a simple-minded collection and that $\Hom(-,-[1])$ vanishes between any two of these objects.

Now let $\mathcal{H_N}$ be the extension closure $\langle E_1[p_1],\dots,E_m[p_m],K_1,\dots,K_r\rangle_{ex}$.
The computations above show that $\mathcal{H_N}$ is actually the category consisting of finite direct sums of these objects,
and hence the functor \[\Hom(\bigoplus_i E_i[p_i] \oplus \bigoplus_j K_j, -)\colon \mathcal{H_N} \longrightarrow \mathrm{mod}\text{-}k^{r+m}\] gives an equivalence.
 
Finally, since $\widetilde{X}$ is smooth and projective, $\mathcal{N}$ is Hom-finite and Krull--Schmidt.
The construction in the proof of \cite[Proposition 3.6]{CSPP22} therefore shows that $\mathcal{H_N}$ is a heart of a bounded $t$-structure on $\mathcal{N}$.
\end{proof}

Theorem \ref{thm:module-type-heart} and Lemma \ref{lem:finite-amplitude-transfer} therefore give the torsion-pair condition of Theorem \ref{thm:perverse-heart-lift-criterion}(i) for every heart of a bounded $t$-structure on $\mathcal{N}$.

\subsection{Smooth blowups}

We now consider blowups of smooth varieties along smooth centers.

Let $X$ be a smooth projective variety over $k$,
and let $\pi\colon\widetilde{X}\longrightarrow X$ be the blowup along a smooth closed subvariety $Z\subset X$ of codimension $c\geq2$.
Set $\mathcal{N}:=\ker\mathbf{R}\pi_\ast\subset\mathrm{D^b}(\widetilde{X})$.
Then $\mathbf{R}\pi_\ast$ induces an equivalence $\mathrm{D^b}(\widetilde{X})/\mathcal{N}\simeq\mathrm{D^b}(X)$.

\begin{proposition}\label{prop:smooth-blowup-kernel-heart}
The category $\mathcal{N}$ is admissible in $\mathrm{D^b}(\widetilde{X})$ and admits a heart of a bounded $t$-structure.
\end{proposition}

\begin{proof}
Let $i\colon E\hookrightarrow\widetilde{X}$ be the exceptional divisor.
For $-c+1\leq t\leq-1$, set $\Phi_t(F):=i_\ast(\pi^\ast F\otimes\mathcal{O}_E(-tE))$.
By Orlov's blowup formula \cite{Orl93}, the functors $\Phi_t$ are fully faithful and $\mathrm{D^b}(\widetilde{X})=\langle\Phi_{-c+1}(\mathrm{D^b}(Z)),\allowbreak\dots,\allowbreak\Phi_{-1}(\mathrm{D^b}(Z)),\allowbreak\mathbf{L}\pi^\ast\mathrm{D^b}(X)\rangle$ is a semiorthogonal decomposition into admissible subcategories.
In particular, $\mathcal{N}=(\mathbf{L}\pi^\ast\mathrm{D^b}(X))^\perp=\langle\Phi_{-c+1}(\mathrm{D^b}(Z)),\allowbreak\dots,\allowbreak\Phi_{-1}(\mathrm{D^b}(Z))\rangle$ is admissible.

For $-c+1\leq t\leq-1$, set $\mathcal{N}_t:=\langle\Phi_{-c+1}(\mathrm{D^b}(Z)),\allowbreak\dots,\allowbreak\Phi_t(\mathrm{D^b}(Z))\rangle$. We construct hearts $\mathcal{H}_{\mathcal{N}_t}$ inductively, starting with $\mathcal{H}_{\mathcal{N}_{-c+1}}:=\Phi_{-c+1}(\Coh(Z))$.
The inclusion $\mathcal{N}_t\subset\mathcal{N}_{t+1}$ is admissible,
so Proposition \ref{prop:admissibility-torsion-pairs} gives the required torsion pairs in $\mathcal{N}_{t+1}$ for $\mathcal{H}_{\mathcal{N}_t}$,
and hence Theorem \ref{thm:perverse-heart-lift-criterion} shows that 
\[\mathcal{H}_{\mathcal{N}_{t+1}}:=\Per(\mathcal{N}_{t+1},\mathcal{N}_t;\Phi_{t+1}(\Coh(Z)),\mathcal{H}_{\mathcal{N}_t})\]
is a heart of a bounded $t$-structure on $\mathcal{N}_{t+1}$.
Repeating this process gives the desired heart $\mathcal{H}_{\mathcal{N}_{-1}}$ on $\mathcal{N}_{-1}=\mathcal{N}$.
\end{proof}

By Proposition \ref{prop:admissibility-torsion-pairs} again, every heart of a bounded $t$-structure on $\mathcal{N}$ satisfies the torsion-pair condition of Theorem \ref{thm:perverse-heart-lift-criterion}(i) in $\mathrm{D^b}(\widetilde{X})$.

\begin{remark}
If $\mathrm{D^b}(Z)$ admits a full exceptional collection, then Orlov's blowup formula gives a full exceptional collection $(E_1,\dots,E_s)$ in $\mathcal{N}$.
By properness and semiorthogonality, we may assume, after shifting the $E_i$, that $\Hom(E_i,E_j[p])=0$ for $i\ne j$ and all $p\leq1$.
In particular, this collection is Ext-exceptional, so \cite[Lemma 3.14]{Mac07} shows that its extension closure $\mathcal{H_N}$ is a heart of a bounded $t$-structure on $\mathcal{N}$.
The same argument in the proof of Theorem \ref{thm:module-for-nodal} then gives $\mathcal{H_N}\simeq\mathrm{mod}\text{-}k^s$.
\end{remark}

\section{Singular surfaces}

In this section we in particular focus on the surface case.

Throughout this section, we work over $\mathbb{C}$.
Let $X$ be a normal projective surface with rational singularities, and let $\pi\colon\widetilde{X}\longrightarrow X$ be its minimal resolution.
Set $\mathcal{N}:=\ker\mathbf{R}\pi_\ast$ and $\gls{Api}:=\Coh(\widetilde{X})\cap\mathcal{N}$.
By Section 4.1, $\mathbf{R}\pi_\ast$ induces an equivalence $\mathrm{D^b}(\widetilde{X})/\mathcal{N}\simeq\mathrm{D^b}(X)$, and $\mathcal{A}_\pi\simeq\mathrm{mod}\text{-}\Gamma$ for some finite-dimensional $\mathbb{C}$-algebra $\Gamma$.

Fix a finite-rank lattice $\Lambda$ and a surjective class map $v\colon\mathrm{K}(X)\longrightarrow\Lambda$. 
Then, by Corollary \ref{cor:product}, there is a natural homeomorphism \[\Stab^{\mathrm{w}}_\Lambda(\widetilde{X},\mathcal{N})\simeq\heartsuit^{b}(\mathcal{N})\times\Stab_\Lambda(X).\]

We now describe the simple objects of $\mathcal{A}_\pi$.

\begin{proposition}\label{prop:KKS}
The simple objects of $\mathcal{A}_\pi$ are exactly the sheaves $\mathcal{O}_C(-1)$, where $C$ runs over the irreducible exceptional curves of $\pi$.
\end{proposition}

\begin{proof}
Each irreducible exceptional curve $C$ is a smooth rational curve,
and the arguments in the proof of Proposition \ref{prop:VdB-null-heart-module} show that $\mathcal{O}_C(-1)$ is simple in $\mathcal{A}_\pi$.
By \cite[Lemma 2.1]{KKS22}, every object of $\mathcal{A}_\pi$ is a finite iterated extension of these sheaves,
and hence there are no other simple objects.
\end{proof}

\subsection{Null hearts and braid actions}

For the remainder of this section, assume that $X$ has canonical singularities $x_1,\ldots,x_r$. 
Let $\Delta_i$ be the Dynkin diagram corresponding to $x_i$ and let $\mathcal{N}_i\subset\mathcal{N}$ and $\mathcal{A}_{\pi,i}\subset\mathcal{A}_\pi$ be the full subcategories consisting of objects supported on $\pi^{-1}(x_i)$.

The simple objects of $\mathcal{A}_{\pi,i}$ form a $\Delta_i$-configuration of $2$-spherical objects, 
so the associated spherical twists define an action of $\Br (\Delta_i)$ on $\mathrm{D^b}(\widetilde{X})$. 

For any simple object $\mathcal{O}_C(-1)$ of $\mathcal{A}_\pi$, the spherical twist $\mathbb{T}_{\mathcal{O}_C(-1)}$ is defined by the functorial triangle $\mathbf{R}\Hom(\mathcal{O}_C(-1),E)\otimes \mathcal{O}_C(-1) \longrightarrow E \longrightarrow \mathbb{T}_{\mathcal{O}_C(-1)}(E)$.
Applying $\mathbf{R}\pi_\ast$ shows that $\mathbf{R}\pi_\ast\circ \mathbb{T}_{\mathcal{O}_C(-1)}=\mathbf{R}\pi_\ast$,
so Proposition \ref{lem:transport} applies to every element of $	\gls{BrX}:=\prod_{i=1}^r\Br(\Delta_i)$.

We first prove a finiteness property of $\mathcal{A}_\pi$.

\begin{lemma}\label{lem:ADE-null-finiteness}
\begin{enumerate}
    \item Every object $E\in\mathcal{A}_\pi$ with $\Hom(E,E)=\mathbb{C}$ is isomorphic to $\beta(\mathcal{O}_C(-1))$ for some $\beta\in\Br_X$ and some irreducible exceptional curve $C$. 
    \item $\mathcal{A}_\pi$ has only finitely many objects with $\Hom(E,E)=\mathbb{C}$ and finitely many torsion pairs.
\end{enumerate}
\end{lemma}

\begin{proof}
(i) Consider $E\in\mathcal{A}_\pi$ with $\Hom(E,E)\simeq\mathbb{C}$.
If $E$ is simple then there is nothing to prove,
so we may assume that $E$ is not simple. 

Since $[E]\neq 0$ and the ADE Cartan matrices are symmetric and positive definite with even diagonal entries,
we have $\chi(E,E)\geq 2$.
On the other hand, since $\pi$ is crepant,
Serre duality gives $\Hom_{\mathcal{N}}(F,G[p])\simeq\Hom_{\mathcal{N}}(G,F[2-p])^\vee$ for $F,G\in\mathcal{N}$,
and hence $2\leq \chi(E,E)=2-\dim\Hom(E,E[1])$.
In particular, $\chi(E,E)=2$ and  $\Hom(E,E[1])=0$.

Expanding $\chi(E,E)=2$ in the simple basis gives a simple object $S$ of $\mathcal{A}_\pi$ with $\chi(S,E)>0$.
As $[E]\neq[S]$, positive definiteness implies 
\[0<\chi([E]-[S],[E]-[S])=4-2\chi(S,E),\]
and hence $1=\chi(S,E)=\dim\Hom(S,E)-\dim\Hom(S,E[1])+\dim\Hom(E,S)$,
so $\Hom(S,E)$ and $\Hom(E,S)$ cannot both vanish,
while they cannot both be nonzero either by the simplicity of $S$.

If $\Hom(E,S)=0$, consider the spherical twist triangle \[\mathbf{R}\Hom(S,E)\otimes S\longrightarrow E\longrightarrow\mathbb{T}_S(E).\]
Its long exact sequence of $\mathcal{A}_\pi$-cohomology gives 
\[\begin{tikzcd}
0 \arrow[r] & {H^{-1}_{\mathcal{A}_\pi}(\mathbb{T}_S(E))} \arrow[r] & {\Hom(S,E)\otimes S} \arrow[r] & E \arrow[lld] \\
& {H^0_{\mathcal{A}_\pi}(\mathbb{T}_S(E))} \arrow[r]    & {\Hom(S,E[1])\otimes S} \arrow[r] & 0
\end{tikzcd}\]
together with $H^1_{\mathcal{A}_\pi}(\mathbb{T}_S(E))\simeq\Hom(S,E[2])\otimes S=0$ and $H^i_{\mathcal{A}_\pi}(\mathbb{T}_S(E))=0$ for any other $i$. 
Since $S$ is simple, the evaluation morphism $\Hom(S,E)\otimes S\longrightarrow E$ is a monomorphism,
so $H^{-1}_{\mathcal{A}_\pi}(\mathbb{T}_S(E))=0$, and hence $\mathbb{T}_S(E)\in\mathcal{A}_\pi$.

Otherwise, we have $\Hom(S,E)=0$. Then the dual argument shows that the coevaluation morphism $E\longrightarrow\Hom(E,S)^\vee\otimes S$ is an epimorphism,
and we have $\mathbb{T}_S^{-1}(E)\in\mathcal{A}_\pi$.

In either case, the resulting object again has endomorphism algebra $\mathbb{C}$ and is of class $[E]-[S]$, so its composition length is one less than that of $E$.
The assertion then follows by induction.

(ii) For $E\in\mathcal{A}_\pi$ with $\Hom(E,E)\simeq\mathbb{C}$, (i) gives $E=\beta(\mathcal{O}_C(-1))$ and $\chi(E,E)=2$.

Since $\chi$ is positive definite and $\mathrm{K}(\mathcal{A}_\pi)$ is a finite-rank lattice, there are only finitely many possibilities for $[E]\in\mathrm{K}(\mathcal{A}_\pi)$,
and hence the lengths of these objects have a common upper bound.
As there are only finitely many braid words $\beta$ of such length and finitely many choices of $C$,
there are only finitely many such $E=\beta(\mathcal{O}_C(-1))$.

Finally, let $\mathcal{T}$ be a torsion class in $\mathcal{A}_\pi$.
We prove by induction on length that every object of $\mathcal{T}$ is a finite iterated extension of objects $F\in\mathcal{T}$ with $\Hom(F,F)=\mathbb{C}$.

Consider $0\neq G\in\mathcal{T}$ with $\Hom(G,G)\neq\mathbb{C}$.
Since $\Hom(G,G)$ is finite-dimensional over $\mathbb{C}$, there exists $0\neq f \in \Hom(G,G)$ which is not an isomorphism.
$\Img f$ and $\cok f$ are quotients of $G \in \mathcal{T}$, so both belong to $\mathcal{T}$ and have strictly smaller length.
Since $G$ is an extension of $\Img f$ and $\cok f$, induction applies.

Thus $\mathcal{T}$ is the extension closure of $\{F\in\mathcal{T}\ | \ \Hom(F,F)\simeq\mathbb{C}\}$.
The finiteness of such objects in $\mathcal{A}_\pi$ gives finitely many possibilities for $\mathcal{T}$.
\end{proof}

We then need a technical lemma to classify hearts of bounded $t$-structures on $\ker \mathbf{R}\pi_\ast$.

\begin{lemma}\label{lem:ADE-cohomology-reduction}
Let $\beta_0\in\Br_X$ and set $\mathcal{B}_0:=\beta_0(\mathcal{A}_\pi)$.
Let $\mathcal{C}\subseteq\mathcal{N}_{\mathcal{B}_0}^{[a,b]}$ be a collection of objects ($a<b$).

Assume that $\Hom(E_1,E_2[p])=0$ for all $E_1,E_2\in\mathcal{C}$ and $p<0$.
Then there exists $\beta\in\Br_X$ such that $\mathcal{B}:=\beta(\mathcal{A}_\pi)$ satisfies $\mathcal{B}_0\subseteq\mathcal{N}_{\mathcal{B}}^{[0,1]}$ and $\mathcal{C}\subseteq\mathcal{N}_{\mathcal{B}}^{[a+1,b]}$.
\end{lemma}

\begin{proof}
We adapt the argument in the proofs of \cite[Theorem 4.6 and 6.11]{HW26}. Consider the set
\[
\mathcal{K}:=\left\{
\mathcal{B}\in\heartsuit^{b}(\mathcal{N})
\ \middle|\
\begin{array}{l}
\mathcal{B}=\beta(\mathcal{A}_\pi)\text{ for some }\beta\in\Br_X,\\
\mathcal{B}_0\subseteq\mathcal{N}_{\mathcal{B}}^{[0,1]},
\mathcal{C}\subseteq\mathcal{N}_{\mathcal{B}}^{[a,b]}
\end{array}
\right\}.
\]
This is a non-empty finite set by Lemma \ref{lem:ADE-null-finiteness}(ii),
so we may choose $\mathcal{B}=\beta(\mathcal{A}_\pi) \in \mathcal{K}$ with minimal $\mathcal{N}_{\mathcal{B}}^{\leq 0}$.

It suffices to prove that $H_{\mathcal{B}}^a(E)=0$ for all $E\in\mathcal{C}$.
Assume the contrary that $H_{\mathcal{B}}^a(E) \neq 0$ for some $E\in\mathcal{C}$.
Since $\mathcal{B}$ is of finite length, we can choose a simple subobject $S\subseteq H_{\mathcal{B}}^a(E)$.
In particular, $S$ is $2$-spherical, so the finite direct sums of $S$ form a torsion-free class in $\mathcal{B}$.
Let $\mathcal{B}':=\langle S[1],{^\perp}S \rangle_{ex}$ be the  corresponding tilt.
Note that this gives $\mathcal{B}\subseteq\mathcal{N}_{\mathcal{B}'}^{[0,1]}$.

We claim that $\mathcal{B}' \in \mathcal{K}$.
First, $\mathbb{T}_S(S)=S[-1]$.
Then for any simple object $ S' \neq S$ of $\mathcal{B}$, taking $\mathcal{B}$-cohomology of the triangle defining $\mathbb{T}_S^{-1}$ gives a short exact sequence 
\[0\longrightarrow\Hom(S',S[1])^\vee\otimes S\longrightarrow\mathbb{T}_S^{-1}(S')\longrightarrow S'\longrightarrow0\]
in $\mathcal{B}$, and $\Hom(\mathbb{T}_S^{-1}(S'),S)=\Hom(S',S[-1])=0$.
Therefore, $\mathbb{T}_S^{-1}(\mathcal{B})\subseteq\mathcal{B}'$ and hence equality holds.
Now, write $S\simeq\beta(\mathcal{O}_C(-1))$ for an irreducible exceptional curve $C$,
and then $\mathbb{T}_S^{-1}\circ\beta\simeq\beta\circ\mathbb{T}_{\mathcal{O}_C(-1)}^{-1}$ gives $\mathcal{B}'=\beta(\mathbb{T}_{\mathcal{O}_C(-1)}^{-1}(\mathcal{A}_\pi))$.

Next, for $F\in\mathcal{C}$, if $H_{\mathcal{B}}^b(F)\neq0$,
then by \cite[Lemma 4.2]{HW26},
we have \[\Hom(H_{\mathcal{B}}^b(F),H_{\mathcal{B}}^a(E))=\Hom(F,E[a-b])=0\] by our assumptions.
Thus $H_{\mathcal{B}}^b(F) \in {^\perp}S \subseteq \mathcal{B}'$ for every $F\in\mathcal{C}$.
Since $\mathcal{C}\subseteq\mathcal{N}_{\mathcal{B}}^{\leq b}$, this shows that $\mathcal{C}\subseteq\mathcal{N}_{\mathcal{B}'}^{\leq b}$.
Also, shifting the inclusion $\mathcal{N}_{\mathcal{B}'}^{\leq 0}\subseteq\mathcal{N}_{\mathcal{B}}^{\leq 0}$ by $[1-a]$ and taking right orthogonals gives 
\[\mathcal{N}_{\mathcal{B}'}^{\geq a}=(\mathcal{N}_{\mathcal{B}'}^{\leq a-1})^\perp\supseteq(\mathcal{N}_{\mathcal{B}}^{\leq a-1})^\perp=\mathcal{N}_{\mathcal{B}}^{\geq a}\supseteq \mathcal{C}.\]
In summary, we obtain $\mathcal{C}\subseteq\mathcal{N}_{\mathcal{B}'}^{[a,b]}$.

On the other hand, for any simple object $S_0$ of $\mathcal{B}_0$, if $H_{\mathcal{B}}^1(S_0)\neq 0$,
then by \cite[Lemma 4.2]{HW26} again,
we have \[\Hom(H_{\mathcal{B}}^1(S_0),H_{\mathcal{B}}^a(E))=\Hom(S_0,E[a-1])=0\]
since $E[a-1]\in\mathcal{N}_{\mathcal{B}_0}^{\geq 1}$.
Thus $H_{\mathcal{B}}^1(S_0) \in {^\perp}S \subseteq \mathcal{B}'$ for every simple object $S_0$ of $\mathcal{B}_0$,
so $S_0\in\mathcal{N}_{\mathcal{B}'}^{[0,1]}$
and hence $\mathcal{B}_0\subseteq\mathcal{N}_{\mathcal{B}'}^{[0,1]}$.
This proves that $\mathcal{B}'\in \mathcal{K}$.

Moreover, since $S[1] \in \mathcal{B}'$ and $S\in\mathcal{B}$,
we see that $\mathcal{N}_{\mathcal{B}'}^{\leq 0}\subsetneq\mathcal{N}_{\mathcal{B}}^{\leq 0}$, contradicting the choice of $\mathcal{B}$.
Therefore, $H_{\mathcal{B}}^a(E)=0$ for all $E\in\mathcal{C}$, and hence $\mathcal{C}\subseteq\mathcal{N}_{\mathcal{B}}^{[a+1,b]}$.
\end{proof}

\begin{theorem}\label{thm:ADE-null-stratum}
The action of $\Br_X$ on $\heartsuit^{b}(\mathcal{N})$ is simply transitive.
\end{theorem}

\begin{proof}
First, applying Lemma \ref{lem:ADE-cohomology-reduction} with $\mathcal{B}_0=\mathcal{C}=\mathcal{A}_\pi$ and $a=0$, $b=1$, we obtain $\beta_0\in\Br_X$ such that $\mathcal{A}_\pi\subseteq\beta_0(\mathcal{A}_\pi)[-1]$. 
As they are both hearts of bounded $t$-structures, we have $\mathcal{A}_\pi=\beta_0(\mathcal{A}_\pi)[-1]$,
and hence $\beta_0^n(\mathcal{A}_\pi)=\mathcal{A}_\pi[n]$ for every $n\in\mathbb{Z}$.

Now let $\mathcal{H_N}\in\heartsuit^{b}(\mathcal{N})$.
As $\mathcal{A}_\pi$ is of finite length with finitely many simple objects, we have $\mathcal{H_N}\subseteq\mathcal{N}_{\mathcal{A}_\pi}^{[a,b]}$ for some $a<b$.

Since $\Hom(F,E[p])=0$ for $E,F\in\mathcal{H_N}$ and $p<0$, applying Lemma \ref{lem:ADE-cohomology-reduction} repeatedly, with $\mathcal{C}=\mathcal{H_N}$ and starting with $\mathcal{B}_0=\mathcal{A}_\pi$, gives  some $\beta\in\Br_X$ such that $\mathcal{H_N}=\beta(\mathcal{A}_\pi)[-b]=(\beta\beta_0^{-b})(\mathcal{A}_\pi)$. 
This proves the transitivity.

We then prove that the stabilizer of $\mathcal{A}_\pi$ is trivial.
It suffices to show that $\gamma\in\Br(\Delta_j)$ satisfying $\gamma(\mathcal{A}_{\pi,j})=\mathcal{A}_{\pi,j}$ is the identity.
Write $\gamma=\beta_1^{-1}\beta_2$ with $\beta_1,\beta_2$ positive
and let $\mathcal{S}$ be the direct sum of the simple objects of $\mathcal{A}_{\pi,j}$.
We then have $\beta_1(\mathcal{S})=\beta_2(\mathcal{S})$.

Following the proof of \cite[Theorem 3.1]{BT11}, we now prove $\beta_1=\beta_2$ by induction on the word length of $\beta_1$.
If one of them is trivial, say $\beta_1=1$, then $\Hom(\mathcal{S},\beta_2(\mathcal{S})[p])=0$ for $p>2$, so \cite[Proposition 3.1(c)]{BT11} implies $\beta_2=1$.

Assume that both $\beta_1,\beta_2$ are nontrivial,
and let $p$ be the largest integer with $\Hom(\mathcal{S},\beta_1(\mathcal{S})[p])\neq 0$.
We may then choose a simple summand $\mathcal{O}_C(-1)\subseteq\mathcal{S}$
such that $\Hom(\mathcal{O}_C(-1),\beta_i(\mathcal{S})[p])\neq0$ for $i=1,2$.

By \cite[Proposition 3.1(b)]{BT11}, $\mathbb{T}_{\mathcal{O}_C(-1)}$ is a common left factor of $\beta_1$ and $\beta_2$.
Applying $\mathbb{T}_{\mathcal{O}_C(-1)}^{-1}$ to $\beta_1(\mathcal{S})=\beta_2(\mathcal{S})$ cancels this common left factor, and induction gives $\beta_1=\beta_2$.
This completes the proof.
\end{proof}

In particular,
every heart of a bounded $t$-structure on $\mathcal{N}$ is equivalent to $\mathcal{A}_\pi\simeq\mathrm{mod}\text{-}\Gamma$,
and there is a bijection $\Br_X \longrightarrow \heartsuit^{b}(\mathcal{N})$ given by  $\beta\longmapsto\beta(\mathcal{A}_\pi)$.

\begin{corollary}\label{cor:ADE-massless-stable}
For $0\neq E\in\mathcal{N}$, $E$ is $2$-spherical if and only if it is simple in some $\mathcal{H_N}\in\heartsuit^{b}(\mathcal{N})$.
\end{corollary}

\begin{proof}
Given a $2$-spherical object $E$, choose $a<b$ such that $E\in\mathcal{N}_{\mathcal{A}_\pi}^{[a,b]}$.
Since $\Hom(E,E[p])=0$ for $p<0$,
we may apply Lemma \ref{lem:ADE-cohomology-reduction} repeatedly,
with $\mathcal{C}=\{E\}$ and starting with $\mathcal{B}_0=\mathcal{A}_\pi$,
to see that $E\in\mathcal{B}$ for some $\mathcal{B}\in\heartsuit^{b}(\mathcal{N})$.

Then Theorem \ref{thm:ADE-null-stratum} gives $\beta_1\in\Br_X$ such that $\mathcal{B}=\beta_1(\mathcal{A}_\pi)$,
and Lemma \ref{lem:ADE-null-finiteness}(i) applied to $\beta_1^{-1}(E)\in\mathcal{A}_\pi$ implies that
$E\simeq\beta_2(\mathcal{O}_C(-1))$ for some $\beta_2\in\Br_X$ and some irreducible exceptional curve $C$.
Therefore, $E$ is simple in the heart $\beta_2(\mathcal{A}_\pi)$.

Conversely, if $E$ is simple in some $\mathcal{H_N}\in\heartsuit^{b}(\mathcal{N})$, then by Theorem \ref{thm:ADE-null-stratum},
we may write $\mathcal{H_N}=\beta(\mathcal{A}_\pi)$ for some $\beta\in\Br_X$,
so $E=\beta(\mathcal{O}_C(-1))$ for some irreducible exceptional curve $C$, and hence is $2$-spherical.
\end{proof}

\paragraph{Partial resolutions.}
Let $Y$ be a normal projective partial crepant resolution of $X$, and write $\pi=fg$, where $g\colon\widetilde{X}\longrightarrow Y$ and $f\colon Y\longrightarrow X$.

We have $\mathbf{R}f_\ast\mathcal{O}_Y=\mathbf{R}\pi_\ast\mathcal{O}_{\widetilde{X}}=\mathcal{O}_X$,
so Proposition \ref{prop:BKS} shows that $\mathbf{R}f_\ast$ induces a Verdier localization $\mathrm{D^b}(Y)/\ker\mathbf{R}f_\ast\simeq\mathrm{D^b}(X)$.
The discussion following Proposition \ref{prop:VdB-null-heart-module}, together with Corollary \ref{cor:product},
gives a homeomorphism \[\mathrm{Stab}^w_\Lambda(Y,\ker\mathbf{R}f_\ast)\simeq\heartsuit^{b}(\ker\mathbf{R}f_\ast)\times\mathrm{Stab}_\Lambda(X),\]
with $\heartsuit^{b}(\ker\mathbf{R}f_\ast)$ equipped with the discrete topology.

We next describe $\heartsuit^{b}(\ker\mathbf{R}f_\ast)$.
Write $\Delta=\coprod_{i=1}^r\Delta_i$ and let $\Delta_Y\subseteq\Delta$ be the full subdiagram corresponding to the curves contracted by $g$. 
By abuse of notation, we identify the vertices of $\Delta$ with the corresponding exceptional curves.

For any full subdiagram $\Delta'\subseteq\Delta$, set $\mathcal{N}_{\Delta'}:=\langle\mathcal{O}_C(-1)\mid C\in\Delta'\rangle$,
let $\mathcal{A}_{\Delta'}\subseteq\mathcal{A}_\pi$ be the extension closure of $\{\mathcal{O}_C(-1)\mid C\in\Delta'\}$,
and let $\Br(\Delta')\subseteq\Br_X$ be the subgroup generated by the spherical twists $\mathbb{T}_{\mathcal{O}_C(-1)}$ for $C\in\Delta'$.

Note that $\mathcal{A}_{\Delta'}\subseteq\mathcal{N}_{\Delta'}$ and $\mathcal{A}_{\Delta'}$ is closed under subobjects and quotients in $\mathcal{A}_\pi$.

\begin{lemma}\label{lem:subdiagram-quotients}
For any full subdiagram $\Delta'\subseteq\Delta$, the following hold.
\begin{enumerate}
\item We have $\mathcal{A}_\pi\cap\mathcal{N}_{\Delta'}=\mathcal{A}_{\Delta'}$,
and $\mathcal{A}_\pi$ descends to the heart $\mathcal{A}_\pi/\mathcal{A}_{\Delta'}$ of a bounded $t$-structure on $\mathcal{N}/\mathcal{N}_{\Delta'}$.
\item $\Br(\Delta')$ acts trivially on $\mathcal{N}/\mathcal{N}_{\Delta'}$.
\end{enumerate}
\end{lemma}

\begin{proof}
(i) Consider the thick subcategory 
$\{E\in\mathcal{N} |H^i_{\mathcal{A}_\pi}(E)\in\mathcal{A}_{\Delta'}\text{ for all }i\in\mathbb{Z}\}$,
which contains the generators of $\mathcal{N}_{\Delta'}$, and hence contains $\mathcal{N}_{\Delta'}$.
It follows that $\mathcal{A}_\pi\cap\mathcal{N}_{\Delta'}=\mathcal{A}_{\Delta'}$,
and that $H^i_{\mathcal{A}_\pi}(E)\in\mathcal{A}_{\Delta'}$ for every $E\in\mathcal{N}_{\Delta'}$ and every $i$,
and hence the assertion follows from Lemma \ref{lem:quotient-heart-criterion}.

(ii) Given $C\in\Delta'$ and $E\in\mathcal{N}$,
the spherical twist triangle shows that $\cone(E\longrightarrow\mathbb{T}_{\mathcal{O}_C(-1)}(E)) = \mathbf{R}\Hom(\mathcal{O}_C(-1),E)\otimes\mathcal{O}_C(-1)[1]$ belongs to $\mathcal{N}_{\Delta'}$.
Therefore, the induced autoequivalence of $\mathbb{T}_{\mathcal{O}_C(-1)}$ on $\mathcal{N}/\mathcal{N}_{\Delta'}$ is the identity.
\end{proof}

\begin{proposition}\label{prop:partial-resolution-braid-cosets}
There is a bijection
\[
\begin{aligned}
\heartsuit^{b}(\ker\mathbf{R}f_\ast)
&\longrightarrow
\left\{
\Br(\Delta_Y)\beta \;\middle|\;
\begin{gathered}
\beta\in\Br_X,\quad \beta\Br(\Delta')\beta^{-1}=\Br(\Delta_Y)\\
\text{for some full subdiagram }\Delta'\subseteq\Delta
\end{gathered}
\right\},\\
\mathbf{R}g_\ast(\beta(\mathcal{A}_\pi))
&\longmapsto \Br(\Delta_Y)\beta.
\end{aligned}
\]
Moreover, there are only finitely many equivalence classes in $\heartsuit^{b}(\ker\mathbf{R}f_\ast)$.
\end{proposition}

\begin{proof}
By Proposition \ref{prop:BKS}, $\mathbf{R}g_\ast$ induces an equivalence $\mathrm{D^b}(\widetilde{X})/\mathcal{N}_{\Delta_Y}=\mathrm{D^b}(Y)$, and this restricts to an equivalence $\mathcal{N}/\mathcal{N}_{\Delta_Y}=\ker\mathbf{R}f_\ast$.

By Theorem \ref{thm:geometric-null-heart} and \ref{thm:perverse-heart-lift-criterion},
every heart $\mathcal{H}$ of a bounded $t$-structure on $\ker\mathbf{R}f_\ast$ lifts to $\Per(\mathcal{N},\mathcal{N}_{\Delta_Y};\mathcal{H},\Coh(\widetilde{X})\cap\ker\mathbf{R}g_\ast)$ on $\mathcal{N}$,
and every such lift is of the form $\beta(\mathcal{A}_\pi)$ for some unique $\beta\in\Br_X$ by Theorem \ref{thm:ADE-null-stratum}. 

We claim that $\beta(\mathcal{A}_\pi)$ descends along $\mathbf{R}g_\ast$ if and only if there exists a full subdiagram $\Delta'\subseteq\Delta$ such that $\beta\Br(\Delta')\beta^{-1}=\Br(\Delta_Y)$.

If $\beta(\mathcal{A}_\pi)$ descends along $\mathbf{R}g_\ast$,
then $\mathbf{R}g_\ast$ is $t$-exact, 
so $H^i_{\beta(\mathcal{A}_\pi)}(E)\in\mathcal{N}_{\Delta_Y}$ for every $E\in\mathcal{N}_{\Delta_Y}$ and every $i$.
Therefore, $\beta(\mathcal{A}_\pi)\cap\mathcal{N}_{\Delta_Y}$ is a heart of a bounded $t$-structure on $\mathcal{N}_{\Delta_Y}$.
Moreover, $t$-exactness implies that $\mathbf{R}g_\ast$ is exact between hearts, 
so $\beta(\mathcal{A}_\pi)\cap\mathcal{N}_{\Delta_Y}$ is closed under subobjects and quotients in $\beta(\mathcal{A}_\pi)$,
and hence $\beta(\mathcal{A}_\pi)\cap\mathcal{N}_{\Delta_Y}=\beta(\mathcal{A}_{\Delta'})$ for some full subdiagram $\Delta'\subseteq\Delta$ since $\beta(\mathcal{A}_\pi)$ is of finite length.

Applying Theorem \ref{thm:ADE-null-stratum} to $g$ then shows that $\gamma\beta(\mathcal{A}_{\Delta'})=\mathcal{A}_{\Delta_Y}$ for some $\gamma\in\Br(\Delta_Y)$.
In particular, we have $\gamma\beta\Br(\Delta')\beta^{-1}\gamma^{-1}=\Br(\Delta_Y)$,
and hence $\beta\Br(\Delta')\beta^{-1}=\Br(\Delta_Y)$ since $\gamma\in\Br(\Delta_Y)$.

Conversely, if $\beta\Br(\Delta')\beta^{-1}=\Br(\Delta_Y)$ for some full subdiagram $\Delta'\subseteq\Delta$,
then for $C\in\Delta'$, $\mathbb{T}_{\beta(\mathcal{O}_C(-1))}=\beta\mathbb{T}_{\mathcal{O}_C(-1)}\beta^{-1}\in\Br(\Delta_Y)$,
so Lemma \ref{lem:subdiagram-quotients}(ii) gives $\mathbf{R}g_\ast\beta(\mathcal{O}_C(-1))=\mathbf{R}g_\ast\mathbb{T}_{\beta(\mathcal{O}_C(-1))}(\beta(\mathcal{O}_C(-1)))=\mathbf{R}g_\ast\beta(\mathcal{O}_C(-1))[-1]$.
By Lemma \ref{lem:subdiagram-quotients}(i), $\mathcal{N}/\mathcal{N}_{\Delta_Y}$ admits a bounded $t$-structure,
so this forces $\mathbf{R}g_\ast(\beta(\mathcal{O}_C(-1)))=0$,
and hence $\beta(\mathcal{N}_{\Delta'})\subseteq\mathcal{N}_{\Delta_Y}$.

Similarly, our assumption $\beta^{-1}\Br(\Delta_Y)\beta=\Br(\Delta')$ gives $\beta^{-1}(\mathcal{N}_{\Delta_Y})\subseteq\mathcal{N}_{\Delta'}$.
Thus $\beta(\mathcal{N}_{\Delta'})=\mathcal{N}_{\Delta_Y}$,
and hence $\beta$ induces an equivalence between $\mathcal{N}/\mathcal{N}_{\Delta'}$ and $\mathcal{N}/\mathcal{N}_{\Delta_Y}$.
This shows that $\beta(\mathcal{A}_\pi)$ descends along $\mathbf{R}g_\ast$ to $\beta(\mathcal{A}_\pi)/\beta(\mathcal{A}_{\Delta'})$.
This proves the claim.

It remains to determine when $\beta_1(\mathcal{A}_\pi)$ and $\beta_2(\mathcal{A}_\pi)$ descend to the same heart on $\ker \mathbf{R}f_\ast$.
For $\gamma\in\Br(\Delta_Y)$, Lemma \ref{lem:subdiagram-quotients}(ii) implies that $\mathbf{R}g_\ast\circ\gamma=\mathbf{R}g_\ast$,
so $\gamma\beta(\mathcal{A}_\pi)$ and $\beta(\mathcal{A}_\pi)$ give the same heart.

Conversely, if $\mathbf{R}g_\ast\beta_1(\mathcal{A}_\pi)=\mathbf{R}g_\ast\beta_2(\mathcal{A}_\pi) \in \heartsuit^{b}(\ker\mathbf{R}f_\ast)$,
then the intersections $\beta_1(\mathcal{A}_\pi)\cap\mathcal{N}_{\Delta_Y}$ and $\beta_2(\mathcal{A}_\pi)\cap\mathcal{N}_{\Delta_Y}$ are hearts of bounded $t$-structures on $\mathcal{N}_{\Delta_Y}$.
By Theorem \ref{thm:ADE-null-stratum} applied to $g$,
there exists some $\gamma\in\Br(\Delta_Y)$ such that $\beta_1(\mathcal{A}_\pi)\cap\mathcal{N}_{\Delta_Y}=\gamma(\beta_2(\mathcal{A}_\pi)\cap\mathcal{N}_{\Delta_Y})=\gamma\beta_2(\mathcal{A}_\pi)\cap\mathcal{N}_{\Delta_Y}$.
Since $\mathbf{R}g_\ast\circ\gamma=\mathbf{R}g_\ast$,
we have $\mathbf{R}g_\ast(\beta_1(\mathcal{A}_\pi))=\mathbf{R}g_\ast(\gamma\beta_2(\mathcal{A}_\pi))$,
and then Proposition \ref{prop:uniqueness-of-descended-lifts} shows that
\begin{align*}
\beta_1(\mathcal{A}_\pi)
&=\Per(\mathcal{N},\mathcal{N}_{\Delta_Y};\mathbf{R}g_\ast(\beta_1(\mathcal{A}_\pi)),\beta_1(\mathcal{A}_\pi)\cap\mathcal{N}_{\Delta_Y})\\
&=\Per(\mathcal{N},\mathcal{N}_{\Delta_Y};\mathbf{R}g_\ast(\gamma\beta_2(\mathcal{A}_\pi)),\gamma\beta_2(\mathcal{A}_\pi)\cap\mathcal{N}_{\Delta_Y})\\
&=\gamma\beta_2(\mathcal{A}_\pi).
\end{align*}
Theorem \ref{thm:ADE-null-stratum} then implies $\beta_1=\gamma\beta_2$.
This proves the bijection assertion.

Finally, every such heart is equivalent to $\mathcal{A}_\pi/\mathcal{A}_{\Delta'}$ for some full subdiagram $\Delta'\subseteq\Delta$.
As there are only finitely many full subdiagrams of $\Delta$, these hearts form only finitely many equivalence classes.
\end{proof}

\subsection{Boundary realization of lifting on K3 surfaces}
In this final subsection, we assume that $\widetilde{X}$ is a K3 surface.
We will prove that $\sharp_{\mathcal{H_N}}(\Stab_{\mathrm{num}}(X))\subseteq\partial\Stab_{\mathrm{num}}(\widetilde{X})$ for every heart $\mathcal{H_N}$ of a bounded $t$-structure on $\mathcal{N}$.

\begin{theorem}\label{thm:br-k3}
If $\widetilde{X}$ is a K3 surface,
then for any $\sigma=(Z,\mathcal{A})\in\Stab_{\mathrm{num}}(X)$ and any $\mathcal{H_N}\in\heartsuit^b(\mathcal{N})$, 
there exists a continuous family, defined for all sufficiently small $\epsilon>0$,
$\sigma_\epsilon:=(Z_\epsilon,\Per(\mathrm{D^b}(\widetilde{X}),\mathcal{N};\mathcal{A},\mathcal{H_N}))\in\Stab_{\mathrm{num}}(\widetilde{X})$ converging to $\sigma^\sharp_{\mathcal{H_N}}$ as $\epsilon \longrightarrow 0$.
\end{theorem}

It suffices to construct the family for $\mathcal{H_N}=\mathcal{A}_\pi[1]=\langle\mathcal{O}_{C_i}(-1)[1]\rangle_{ex}$.
The braid action of Theorem \ref{thm:ADE-null-stratum} will then prove the assertion for every heart in $\heartsuit^{b}(\mathcal{N})$.

For the numerical Grothendieck group  of the singular surface $X$, we define $\mathrm{K_{num}}(X)$ by the quotient of $\mathrm{K}(\mathrm{D^b}(X))$ by the right kernel of the Euler characteristic 
$\chi \colon \mathrm{K}(\mathrm{D^{perf}}(X))\times \mathrm{K}(\mathrm{D^b}(X))\longrightarrow\mathbb{Z}$.

We fix a numerical stability condition $\sigma=(Z,\mathcal{A})\in\Stab_{\mathrm{num}}(X)$,
denote the numerical class maps for $\widetilde{X}$ and $X$ by $v$ and $v_X$ respectively,
and fix Euclidean norms on $\mathrm{K_{num}}(\widetilde{X})_{\mathbb{R}}$ and $\mathrm{K_{num}}(X)_{\mathbb{R}}$.

For simplicity, we write $\mathcal{B}:=\Per(\mathrm{D^b}(\widetilde{X}),\mathcal{N};\mathcal{A},\mathcal{A}_\pi[1])$, $S_i:=\mathcal{O}_{C_i}(-1)[1]$,
and $\ell(T):=\mathrm{length}_{\mathcal{A}_\pi[1]}(T)$.
We also choose a quasi-inverse $C\colon\mathrm{D^b}(X)\longrightarrow\mathcal{V}_0$ of the equivalence in Remark \ref{lem:window-quotient}(i).

\begin{lemma}\label{lem:k3-core-bound}
There exists a constant $K_\sigma>0$ such that $\|v(C(A))\|\leq K_\sigma|Z(A)|$ for every $\sigma$-stable $A\in\mathcal{A}$.
\end{lemma}

\begin{proof}
Given a $\sigma$-stable object $A\in\mathcal{A}$, we consider $C(A)\in\mathcal{B}$.
Note first that $\End_{\widetilde{X}}(C(A))=\End_X(A)=\mathbb{C}$,
and by Serre duality we then have 
\[\chi(C(A),C(A))=2-\dim\Hom(C(A),C(A)[1])\leq2.\]

The restriction of $\chi$ to $(\ker\pi_\ast)_{\mathbb{R}}$ is positive definite since $\chi(S_i,S_j)=-C_i\cdot C_j$.
Let $e:=\pr_{(\ker\pi_\ast)_{\mathbb{R}}}(v(C(A)))$ be the $\chi$-orthogonal projection.
Then we obtain that $v(C(A))-e\in(\ker\pi_\ast)_{\mathbb{R}}^\perp$ and $\pi_\ast(v(C(A))-e)=v_X(A)$.

Since $\pi_\ast\colon(\ker\pi_\ast)_{\mathbb{R}}^\perp\longrightarrow\mathrm{K_{num}}(X)_{\mathbb{R}}$ is an isomorphism,
there exists some $b>0$ independent of $A$ such that $\|v(C(A))-e\|\leq b\|v_X(A)\|$.

As $\chi$ is a bilinear form on a finite-dimensional space,
we have \[|\chi(v(C(A))-e,v(C(A))-e)|\leq a\|v(C(A))-e\|^2\leq ab^2\|v_X(A)\|^2,\]
for some $a>0$ independent of $A$, and hence
\[\chi(e,e)=\chi(C(A),C(A))-\chi(v(C(A))-e,v(C(A))-e)\leq2+ab^2\|v_X(A)\|^2.\]

By the positive definiteness of $\chi$ on $(\ker\pi_\ast)_{\mathbb{R}}$, there exists $c>0$ independent of $A$ such that $\|e\|^2\leq c\chi(e,e)$.
Therefore, we now have \[\|v(C(A))\|\leq\|v(C(A))-e\|+\|e\|\leq\sqrt{2c}+b(1+\sqrt{ac})\|v_X(A)\|.\]
Since $Z(A)\neq0$, we have $v_X(A)\neq0$. Choose $M,N>0$ such that $M\leq\|v_X(A)\|$ and that $\|v_X(A)\|\leq N|Z(A)|$ for every $\sigma$-stable $A\in\mathcal{A}$.
Then the assertion holds with $K_\sigma:=N(\sqrt{2c}M^{-1}+b(1+\sqrt{ac}))$.
\end{proof}

\begin{lemma}\label{lem:k3-core-mass-bound}
There exists a constant $L_\sigma>0$ such that $\|v(F)\|\leq L_\sigma m_\sigma(\mathbf{R}\pi_\ast F)$ for every $F\in(\mathcal{A}_\pi[1])^{\perp}\subseteq\mathcal{B}$.
\end{lemma}
\begin{proof}
We first consider the Jordan--H\"older filtrations of the Harder--Narasimhan factors of $\mathbf{R}\pi_\ast F$,
and concatenate them to obtain a filtration of $\mathbf{R}\pi_\ast F$.
Then, by repeated application of Lemma \ref{lem:F-lifting}, this lifts to a finite filtration of $F$ with factors $F_j\in(\mathcal{A}_\pi[1])^{\perp}$ whose pushforwards are $\sigma$-stable.

Since $\|v(F)\|\leq\sum_j\|v(F_j)\|$ and $m_\sigma(\mathbf{R}\pi_\ast F)=\sum_j|Z(\mathbf{R}\pi_\ast F_j)|$,
it suffices to find $L_\sigma$ such that  $\|v(F)\|\leq L_\sigma|Z(\mathbf{R}\pi_\ast F)|$ for all $F\in(\mathcal{A}_\pi[1])^{\perp}$ with $\mathbf{R}\pi_\ast F$ $\sigma$-stable.

Now given $F\in(\mathcal{A}_\pi[1])^{\perp}$,
by Proposition \ref{prop:core-equivalence},
there is a short exact sequence $0\longrightarrow C(A)\longrightarrow F\longrightarrow T\longrightarrow0$,
where $A:=\mathbf{R}\pi_\ast F$ and $T\in\mathcal{A}_\pi[1]$.

Fix a projective generator $P$ of $\mathcal{A}_\pi[1]$.
Applying $\Hom(P,-)$ to the short exact sequence gives $\Hom(P,F)=0 \longrightarrow \Hom(P,T)\longrightarrow\Hom(P,C(A)[1])$.

Let $c_P:=|\chi(P,-)|$ be the operator norm of $\chi(P,-)$ on $\mathrm{K_{num}}(\widetilde{X})_{\mathbb{R}}$ and $c_{\mathcal{H_N}}:=\max_i\|[S_i]\|$.
Then:
\begin{align*}
\|v(F)\|
&\leq\|v(C(A))\|+\|v(T)\|
\\
&\leq\|v(C(A))\|+c_{\mathcal{H_N}}\ell(T)
\\
&\leq\|v(C(A))\|+c_{\mathcal{H_N}}\dim\Hom(P,T)
&&\text{($P$ is a projective generator)}\\
&\leq\|v(C(A))\|+c_{\mathcal{H_N}}\dim\Hom(P,C(A)[1])\\
&=\|v(C(A))\|-c_{\mathcal{H_N}}\chi(P,C(A))
&&\text{(Serre duality)}\\
&\leq(1+c_{\mathcal{H_N}}c_P)\|v(C(A))\|
\\
&\leq(1+c_{\mathcal{H_N}}c_P)K_\sigma|Z(A)|
&&\text{(Lemma \ref{lem:k3-core-bound})}.
\end{align*}
The assertion follows with $L_\sigma:=(1+c_{\mathcal{H_N}}c_P)K_\sigma$.
\end{proof}

\begin{proof}[Proof of Theorem \ref{thm:br-k3}]
Write $\sigma^\sharp=\sigma^\sharp_{\mathcal{A}_\pi[1]}=(Z_0,\mathcal{B})$, where $Z_0=Z\circ\mathbf{R}\pi_\ast$.
Take the exceptional real divisor $D$ with $D\cdot C_i=-1$ for every $i$.
Set $\eta:=D\cdot\ch_1$ and $Z_\epsilon:=Z_0-\epsilon\eta$ for $\epsilon>0$.
We will show that $(Z_\epsilon, \mathcal{B})\in \Stab_{\mathrm{num}}(\widetilde{X})$ for $\epsilon>0$ sufficiently small.

First, we claim that $Z_\epsilon$ is a stability function on $\mathcal{B}$ for $\epsilon>0$ sufficiently small.
For $0 \neq E \in \mathcal{B}$, we have $\Im Z_\epsilon(E)=\Im Z_0(E)\geq0$.
If $\Im Z_0(E)=0$, consider the short exact sequence
$0\longrightarrow T_E\longrightarrow E\longrightarrow F_E\longrightarrow0$
given by Proposition \ref{prop:core-equivalence}(i),
where $T_E\in\mathcal{A}_\pi[1]$ and $F_E\in(\mathcal{A}_\pi[1])^{\perp}$.
Since $\mathbf{R}\pi_\ast T_E=0$, we have
$\mathbf{R}\pi_\ast F_E=\mathbf{R}\pi_\ast E$,
which is either zero or $\sigma$-semistable of phase $1$,
and hence $m_\sigma(\mathbf{R}\pi_\ast F_E)=|Z_0(F_E)|$.

Let $L_\sigma$ be as in Lemma \ref{lem:k3-core-mass-bound}, and choose $\epsilon>0$ such that $u:=\epsilon L_\sigma|\eta|<1/4$.
Since $\eta(S_i)=1$, we have $\eta(T_E)=\ell(T_E)$.
Then
\begin{align*}
\Re Z_\epsilon(E)
&=\Re Z_0(F_E)-\epsilon\eta(F_E)-\epsilon\eta(T_E)\\
&=-|Z_0(F_E)|-\epsilon\eta(F_E)-\epsilon\ell(T_E)\\
&\leq-|Z_0(F_E)|+\epsilon|\eta|\|v(F_E)\|-\epsilon\ell(T_E)\\
&\leq-|Z_0(F_E)|+\epsilon|\eta|L_\sigma m_\sigma(\mathbf{R}\pi_\ast F_E)-\epsilon\ell(T_E)\\
&=-(1-u)|Z_0(F_E)|-\epsilon\ell(T_E)\\
&<0.
\end{align*}

To prove the Harder--Narasimhan property of $\sigma_\epsilon=(Z_\epsilon,\mathcal{B})$,
we use the criterion of \cite[Corollary 3.6]{Bay19}.
Given $E\in\mathcal{B}$, it suffices to show that the set
\[\{v(G)\mid G\subseteq E,\ \Re Z_\epsilon(G)<\max\{0,\Re Z_\epsilon(E)\}\}\] is finite.

For such a subobject $G\subseteq E$, we consider the torsion decompositions of $G$ and $E$. 
Then we have $T_G\subseteq T_E$ and $\mathbf{R}\pi_\ast G\subseteq\mathbf{R}\pi_\ast E$.
By \cite[Lemma 5.5]{Bay19}, we have
$\Re Z_0(G)\geq m_\sigma(\mathbf{R}\pi_\ast G)-m_\sigma(\mathbf{R}\pi_\ast E)+\Re Z_0(E)$.
Lemma \ref{lem:k3-core-mass-bound} gives
$|\epsilon\eta(F_G)|\leq u m_\sigma(\mathbf{R}\pi_\ast G)$, so
\begin{align*}
\Re Z_\epsilon(G)
&=\Re Z_0(G)-\epsilon\eta(F_G)-\epsilon\ell(T_G)\\
&\geq\Re Z_0(G)-u m_\sigma(\mathbf{R}\pi_\ast G)-\epsilon\ell(T_E)\\
&\geq(1-u)m_\sigma(\mathbf{R}\pi_\ast G)-m_\sigma(\mathbf{R}\pi_\ast E)+\Re Z_0(E)-\epsilon\ell(T_E).
\end{align*}
Since $\Re Z_\epsilon(G)<\max\{0,\Re Z_\epsilon(E)\}$, this bounds $m_\sigma(\mathbf{R}\pi_\ast G)$ from above independently of $G$.
Then
\[\|v(G)\|
\leq\|v(T_G)\|+\|v(F_G)\|\leq c_{\mathcal{H_N}}\ell(T_E)+L_\sigma m_\sigma(\mathbf{R}\pi_\ast G),\]
where $c_{\mathcal{H_N}}:=\max_i\|[S_i]\|$.

Therefore, $\|v(G)\|$ is bounded  independently of $G$, and discreteness of $\mathrm{K_{num}}(\widetilde{X})$ implies that the set is finite.
This proves the Harder--Narasimhan property of $(Z_\epsilon,\mathcal{B})$.

We now prove the support property of $\sigma_\epsilon$.
If $E\in\mathcal{B}$ is $Z_\epsilon$-semistable of phase $1$,
then by the inequalities
\begin{align*}
|Z_\epsilon(E)|
&\geq \epsilon\ell(T_E)+(1-u)|Z_0(F_E)|,\\
\|v(E)\|
&\leq c_{\mathcal{H_N}}\ell(T_E)+L_\sigma|Z_0(F_E)|,
\end{align*}
we obtain that $\|v(E)\|\leq K_\epsilon|Z_\epsilon(E)|$, where
$K_\epsilon:=\max\{2L_\sigma,c_{\mathcal{H_N}}/\epsilon\}$.

Now if $E\in \mathcal{B}$ is $Z_\epsilon$-semistable of phase $\phi \in (0,1)$,
then $E\in (\mathcal{A}_\pi[1])^{\perp}$ since
every nonzero object of $\mathcal{A}_\pi[1]$ has $Z_\epsilon$-phase $1$.
By the construction in the proof of Theorem \ref{thm:abstract-HN-lift},
the $\sigma$-Harder--Narasimhan filtration of $\mathbf{R}\pi_\ast E$ lifts to the $\sigma^\sharp$-Harder--Narasimhan filtration of $E$,
with factors $F_1,\ldots,F_m \in (\mathcal{A}_\pi[1])^{\perp}$.

For each $j$, Lemma \ref{lem:k3-core-mass-bound} gives
\[|Z_\epsilon(F_j)-Z_0(F_j)|=\epsilon|\eta(F_j)|\leq \epsilon|\eta|\,\|v(F_j)\|\leq \epsilon|\eta|L_\sigma m_\sigma(\mathbf{R}\pi_\ast F_j)=u|Z_0(F_j)|.\]

Taking arguments in $(0,\pi]$, we obtain
$|\arg Z_\epsilon(F_j)-\arg Z_0(F_j)|\leq\pi\delta_\epsilon$,
where $\delta_\epsilon:=\pi^{-1}\arcsin u$.
Then by $Z_\epsilon$-semistability of $E$,
$\arg Z_\epsilon(F_1)\leq\pi\phi\leq\arg Z_\epsilon(F_m)$,
and it follows that $\phi^\pm_\sigma(\mathbf{R}\pi_\ast E)\in[\phi-\delta_\epsilon,\phi+\delta_\epsilon]$.
Consequently, 
\[
\Re(e^{-i\pi\phi}Z_0(E))
=\sum_j|Z_0(F_j)| \cos(\arg Z_0(F_j)-\pi\phi)
\geq\cos(\pi\delta_\epsilon)m_\sigma(\mathbf{R}\pi_\ast E),
\]
so by Lemma \ref{lem:k3-core-mass-bound} again, we obtain that
\begin{align*}
K_\epsilon|Z_\epsilon(E)|
=K_\epsilon\Re(e^{-i\pi\phi}Z_\epsilon(E))
&\geq K_\epsilon\Re(e^{-i\pi\phi}Z_0(E))-K_\epsilon|\epsilon\eta(E)|\\
&\geq K_\epsilon(\cos(\pi\delta_\epsilon)-u)m_\sigma(\mathbf{R}\pi_\ast E)\\
&\geq\frac{K_\epsilon}{2}m_\sigma(\mathbf{R}\pi_\ast E)\\
&\geq\|v(E)\|.
\end{align*}

In summary, $\|v(E)\|\leq K_\epsilon|Z_\epsilon(E)|$ for every $Z_\epsilon$-semistable $E\in\mathcal{B}$.
This proves that $\sigma_\epsilon=(Z_\epsilon,\mathcal{B})\in\Stab_{\mathrm{num}}(\widetilde{X})$.

Now, the support property and \cite[Lemma 2.9]{Bay19} give the continuity of $\sigma_\epsilon$ for $\epsilon>0$.
Moreover, as $\Im Z_\epsilon=\Im Z_0$, every $Z_\epsilon$-semistable object of phase $1$ is also $\sigma^\sharp$-semistable of phase $1$,
so the same arguments as above show that the $\sigma^\sharp$-Harder--Narasimhan phases of every $\sigma_\epsilon$-semistable object of phase $\phi$ lie in $[\phi-\delta_\epsilon,\phi+\delta_\epsilon]$.
Then we have $\delta_\epsilon\longrightarrow0$ and $\|Z_\epsilon-Z_0\| \longrightarrow 0$ as $\epsilon\longrightarrow 0$, and hence $\sigma_\epsilon \longrightarrow \sigma^\sharp$.

Finally, for any  $\mathcal{H_N}\in\heartsuit^{b}(\mathcal{N})$,
we choose $\beta\in\Br_X$ with $\beta(\mathcal{A}_\pi[1])=\mathcal{H_N}$ given by Theorem \ref{thm:ADE-null-stratum},
and then $\beta\cdot\sigma_\epsilon$ is the desired family in $\Stab_{\mathrm{num}}(\widetilde{X})$ by Proposition \ref{lem:transport}.
\end{proof}

\begin{remark}
The K3 assumption is only used in the proof of Lemma \ref{lem:k3-core-bound},
and the conclusion of Theorem \ref{thm:br-k3} holds for any projective surface $X$ with ADE singularities and any $\sigma\in\Stab_{\mathrm{num}}(X)$ satisfying the bound in Lemma \ref{lem:k3-core-bound}.

In particular, for geometric stability conditions on $X$ constructed in \cite{Cho26}, one can verify the bound and the construction above recovers the path $\sigma_\epsilon \longrightarrow \sigma_{\widetilde{X}}$.
\end{remark}

We conclude the section with the following observation that the family constructed in the proof of Theorem \ref{thm:br-k3}, after reparametrization $\epsilon=s^{-1}$, is also quasi-convergent in the sense of \cite[Definition 2.8]{HLJR26}.

\begin{proposition}\label{prop:k3-quasi-convergence}
The family $\sigma_s=(Z_0-s^{-1}\eta,\mathcal{B})$, defined for all $s \gg 0$, is quasi-convergent.
\end{proposition}
\begin{proof}
We first show that any path of Bridgeland stability conditions converging in the metric $d$ satisfies condition (i) of \cite[Definition 2.8]{HLJR26}.

For such a path $\sigma_t=(Z_t,\mathcal{P}_t)\longrightarrow\sigma_\infty=(Z_\infty,\mathcal{P}_\infty)$ as $t\longrightarrow\infty$,
we claim that any nonzero object $E$ is $\mathcal{P}_\infty$-semistable if and only if it is limit semistable,
that is, $\phi^+_{\sigma_t}(E)-\phi^-_{\sigma_t}(E)\longrightarrow0$.

If $E$ is $\mathcal{P}_\infty$-semistable, then $\phi^\pm_{\sigma_t}(E)$ both converge to its $\mathcal{P}_\infty$-phase,
so $E$ is limit semistable.
Conversely, if $E$ is limit semistable, then
\[0\leq\phi^+_{\mathcal{P}_\infty}(E)-\phi^-_{\mathcal{P}_\infty}(E)\leq\phi^+_{\sigma_t}(E)-\phi^-_{\sigma_t}(E)+2d_{\mathrm{sl}}(\mathcal{P}_t,\mathcal{P}_\infty)\longrightarrow0.\]
Therefore, $\phi^+_{\mathcal{P}_\infty}(E)=\phi^-_{\mathcal{P}_\infty}(E)$, 
i.e., $E$ is $\mathcal{P}_\infty$-semistable.

Let $G_1,\ldots,G_n$ be the $\mathcal{P}_\infty$-Harder--Narasimhan factors of a nonzero object $G$.
Each $G_i$ is limit semistable, and its $\sigma_t$-average phase $\phi_t(G_i)$ (see \cite[Definition 2.4]{HLJR26}) satisfies $\phi^-_{\sigma_t}(G_i)\leq\phi_t(G_i)\leq\phi^+_{\sigma_t}(G_i)$, so $\phi_t(G_i)\longrightarrow\phi_{\mathcal{P}_\infty}(G_i)$.
Therefore, $\phi_t(G_i)-\phi_t(G_{i+1}) \longrightarrow \phi_{\mathcal{P}_\infty}(G_i)-\phi_{\mathcal{P}_\infty}(G_{i+1})>0$ for every $1 \leq i<n$,
and hence the $\mathcal{P}_\infty$-Harder--Narasimhan filtration is the desired filtration.
This proves condition (i), which applies in particular to $\sigma_s\longrightarrow\sigma^\sharp$.

For (ii), set $\ell_s(E):=\log m_{\sigma_s}(E)+i\pi\phi_s(E)$ and $\ell_s(E_1/E_2):=\ell_s(E_1)-\ell_s(E_2)$ as in \cite[Definition 2.4]{HLJR26}.
We must show that for every pair of limit semistable objects $E_1,E_2$,
\[\lim_{s\rightarrow\infty}\frac{\ell_s(E_1/E_2)}{1+|\ell_s(E_1/E_2)|} \text{ exists.}\]
Let $E$ be limit semistable with respect to the family $\sigma_s$.
Then $E$ is $\sigma^\sharp$-semistable.

If $E \notin \mathcal{N}$, then $Z_s(E) \longrightarrow Z_0(E) \neq 0$.
By \cite[Lemma 2.5]{HLJR26}, applied with $\sigma=\sigma_s$ and $\epsilon=\phi^+_{\sigma_s}(E)-\phi^-_{\sigma_s}(E)\longrightarrow0$, we obtain $\log(m_{\sigma_s}(E)/|Z_s(E)|)\longrightarrow0$,
and hence $m_{\sigma_s}(E)\longrightarrow|Z_0(E)|$.

If $E\in\mathcal{N}$, then $\sigma^\sharp$-semistability implies that $E$ is a shift of a nonzero object of $\mathcal{A}_\pi[1]$ and is therefore $\sigma_s$-semistable.
Thus $sm_{\sigma_s}(E)=s|Z_s(E)|=|\eta(E)|>0$ for all sufficiently large $s$.

Now let $E_1,E_2$ be limit semistable with respect to the family $\sigma_s$.
Their $\sigma_s$-average phases converge to their $\sigma^\sharp$-phases by the first part of the proof,
so $\Im\ell_s(E_1/E_2)=\pi(\phi_s(E_1)-\phi_s(E_2))$ has a finite limit.

If both $E_i$ belong to $\mathcal{N}$ or both lie outside $\mathcal{N}$,
then $m_{\sigma_s}(E_1)/m_{\sigma_s}(E_2)$ converges to a positive finite value.
Thus $\ell_s(E_1/E_2)$ has a finite limit,
and hence so does $\ell_s(E_1/E_2)/(1+|\ell_s(E_1/E_2)|)$.

If $E_1\in\mathcal{N}$ and $E_2\notin\mathcal{N}$,
then we have $m_{\sigma_s}(E_1)/m_{\sigma_s}(E_2)\longrightarrow 0 $, so $\Re\ell_s(E_1/E_2)\longrightarrow-\infty$.
Since its imaginary part is bounded, we obtain $\ell_s(E_1/E_2)/(1+|\ell_s(E_1/E_2)|)\longrightarrow-1$.

Similarly, if $E_1\notin\mathcal{N}$ and $E_2\in\mathcal{N}$, then $m_{\sigma_s}(E_1)/m_{\sigma_s}(E_2)\longrightarrow +\infty$ and $\ell_s(E_1/E_2)/(1+|\ell_s(E_1/E_2)|)\longrightarrow1$.
This verifies condition (ii) and proves quasi-convergence.
\end{proof}

Note that convergence in $d$ alone need not imply condition (ii),
as it does not control the relative rates at which masses tend to zero.

The conclusions of Proposition \ref{prop:k3-quasi-convergence} also hold for the families $\beta \cdot \sigma_s$ obtained by braid transport,
which preserves the phases and masses of corresponding objects and commutes with $\mathbf{R}\pi_\ast$.

\bibliographystyle{alpha}
\bibliography{ref}

\noindent
\textsc{Department of Mathematics, National Taiwan University, Taipei, Taiwan}

\smallskip

\noindent
\textit{Email address:} \texttt{tychou1225@ntu.edu.tw}
\end{document}